\documentclass[10pt]{article}

\usepackage[a4paper,margin=1in]{geometry}
\usepackage[T1]{fontenc}
\usepackage{lmodern}
\usepackage{amsmath,amssymb,amsthm,mathtools,mathrsfs}
\usepackage{microtype}
\usepackage[skip=0.45\baselineskip plus 2pt,indent=0pt]{parskip}
\usepackage{enumitem}
\usepackage{xcolor}
\usepackage{aliascnt}
\usepackage[colorlinks=true,linkcolor=blue,citecolor=blue,urlcolor=blue]{hyperref}
\usepackage[nameinlink,capitalize,noabbrev]{cleveref}

\numberwithin{equation}{section}

\allowdisplaybreaks
\setlist{itemsep=2pt,topsep=4pt}
\hypersetup{
  pdftitle={Nearly Sharp Bounds for Lattice Coverings by Convex Bodies},
  pdfauthor={Heng Li and Xizhi Liu},
  pdfkeywords={Lattice covering, convex body, geometry of numbers, finite-field Kakeya, random slabs, Poisson process}
}

\newtheorem{theorem}{Theorem}[section]

\newaliascnt{lemma}{theorem}
\newtheorem{lemma}[lemma]{Lemma}
\aliascntresetthe{lemma}

\newaliascnt{proposition}{theorem}
\newtheorem{proposition}[proposition]{Proposition}
\aliascntresetthe{proposition}

\newaliascnt{corollary}{theorem}
\newtheorem{corollary}[corollary]{Corollary}
\aliascntresetthe{corollary}

\theoremstyle{remark}
\newtheorem*{remark}{Remark}

\newcommand{\R}{\mathbb R}
\newcommand{\Z}{\mathbb Z}
\newcommand{\F}{\mathbb F}
\newcommand{\T}{\mathbb T}
\newcommand{\E}{\mathbb E}
\newcommand{\Prb}{\mathbb P}
\newcommand{\vol}{\operatorname{vol}}
\newcommand{\Sph}{\mathbb S}
\newcommand{\dist}{\operatorname{dist}}
\newcommand{\covrad}{\operatorname{covrad}}
\newcommand{\dd}{\,\mathrm{d}}
\newcommand{\1}{\mathbf 1}
\newcommand{\cG}{\mathcal G}
\newcommand{\cL}{\mathcal L}
\newcommand{\cP}{\mathcal P}
\newcommand{\ip}[2]{\langle #1,#2\rangle}
\DeclareMathOperator{\Var}{Var}

\title{\Large\bfseries Nearly Sharp Bounds for Lattice Coverings by Convex Bodies}
\author{%
Heng~Li\thanks{\scriptsize School of Mathematics, Shandong University, Jinan, China, and Extremal Combinatorics and Probability Group, Institute for Basic Science, Daejeon, South Korea. Email:~\texttt{heng.li@sdu.edu.cn}.}\and
Xizhi Liu\thanks{\scriptsize School of Mathematical Sciences, University of Science and Technology of China, Hefei, China. Email:~\texttt{liuxizhi@ustc.edu.cn}.}}
\date{}

\begin{document}
\maketitle

\begin{abstract}
For an $n$-dimensional convex body $K$, let $\theta_L(K)$ denote its lattice covering density, and let $\Theta_L^{\mathrm{conv}}(n)$ and $\Theta_L^{\mathrm{sym}}(n)$ be the corresponding worst-case quantities over all convex bodies and over origin-symmetric convex bodies, respectively.
Before this work, these quantities were known only to lie between a lower bound of order $n$ and an upper bound of order $n^2$, so even their polynomial order was undetermined.
We prove that there are absolute constants $c,C>0$ such that
\[
  c n\log n
  \le
  \Theta_L^{\mathrm{sym}}(n)
  \le
  \Theta_L^{\mathrm{conv}}(n)
  \le
  Cn\log n\,(\log\log n)^{10/3+o(1)}.
\]
Thus both worst-case quantities are $n\log n\,(\log n)^{o(1)}$, and the upper and lower bounds differ by a factor at most $(\log\log n)^{10/3+o(1)}$.

For the upper bound, a vertical--horizontal amplification based on weighted Boolean cubes combines covering estimates for low-codimensional sections into an exact lattice covering of an arbitrary convex body.
For the lower bound, a random-slab construction and Poisson witnesses on flat tori show, with positive probability, that the resulting body admits no lattice covering of density below $c n\log n$.
\end{abstract}



\section{Introduction}

\subsection{Lattice coverings and the extremal problem}

A \emph{convex body} in $\R^n$ is a compact convex set with nonempty interior.
If $v_1,\ldots,v_n$ is a basis of a full-rank lattice $\Lambda\subseteq\R^n$, write $\det\Lambda\coloneqq|\det(v_1,\ldots,v_n)|$ for its positive covolume.
The translates $\{K+\lambda\colon\lambda\in\Lambda\}$ form a \emph{lattice covering} if
\[
  \Lambda+K
  \coloneqq
  \bigcup_{\lambda\in\Lambda}(\lambda+K)
  =\R^n.
\]
The density of this covering is $\vol_n(K)/\det\Lambda$.
The \emph{lattice covering density} of $K$ is
\[
  \theta_L(K)
  \coloneqq
  \inf\left\{
    \frac{\vol_n(K)}{\det\Lambda}
    \colon
    \Lambda+K=\R^n
  \right\}.
\]
This quantity is invariant under invertible affine transformations of $K$.

We study the largest possible lattice covering density in a fixed dimension.
Define
\begin{align*}
  \Theta_L^{\mathrm{conv}}(n)
  &\coloneqq
  \sup\left\{
    \theta_L(K)
    \colon
    K\subseteq\R^n\text{ is a convex body}
  \right\},\\
  \Theta_L^{\mathrm{sym}}(n)
  &\coloneqq
  \sup\left\{
    \theta_L(K)
    \colon
    K\subseteq\R^n\text{ is an origin-symmetric convex body}
  \right\}.
\end{align*}
The defining classes immediately give $\Theta_L^{\mathrm{sym}}(n)\le\Theta_L^{\mathrm{conv}}(n)$.

All logarithms are natural unless a base is displayed, and all asymptotic statements involving $n$ are taken as $n\to\infty$.
In the upper-bound construction, statements involving only the auxiliary dimension $d$ are taken as $d\to\infty$.
For positive quantities $a$ and $b$, write $a\asymp b$ if $a/b$ is bounded above and below by positive absolute constants in the relevant limit, and write $a\gg b$ if $a/b\to\infty$.

Rogers~\cite{RogersBook} gives a classical treatment of packing and covering.
For the broader geometry of numbers, see Gruber and Lekkerkerker~\cite{GruberLekkerkerker1987}; Conway and Sloane~\cite{ConwaySloane1988} give an extensive account of lattices and sphere packings.

\subsection{Previous bounds}

For comparison, Rogers~\cite{Rogers1957} proved that arbitrary translates cover every $n$-dimensional convex body with density $n\log n+n\log\log n+O(n)$, while Fejes T{\'o}th~\cite{FejesToth2009} obtained essentially the same density using $O(\log n)$ cosets of one lattice.
The restriction to a single lattice is substantially harder.

For Euclidean balls, Dumer~\cite{Dumer2007} obtained the sharper translative bound $(1/2+o(1))n\log n$, and Bukh, Gao, Liu, Pikhurko, and Sun~\cite{BukhEtAl2025} attained this density with controlled covering multiplicity and identified a limitation of random periodic constructions.

For lattice coverings, Rogers~\cite{Rogers1959} proved the universal estimate $\theta_L(K)\le n^{\log_2\log n+O(1)}$.
Ordentlich, Regev, and Weiss~\cite{ORW} established the first polynomial universal bound, $\theta_L(K)\le Cn^2$.
In subsequent work, the same authors~\cite{ORWSmooth} obtained uniform multiplicity estimates at polynomial covering density and studied random Construction-A lattices.

Near-linear upper bounds with polylogarithmic losses were already known for bodies with sufficiently many orthogonal symmetries through Gritzmann~\cite{Gritzmann1985}, and for locally anti-blocking bodies and polytopes with few vertices through Schymura, Wang, and Xue~\cite{SWX}.
For Euclidean unit balls $B_2^n$, Gao, Liu, Pikhurko, and Sun~\cite{GLPS} proved
\[
  \theta_L(B_2^n)
  =
  O\bigl(n(\log n)^{\beta_{\mathrm{ball}}}\bigr),
  \quad\text{where}\quad
  \beta_{\mathrm{ball}}
  \coloneqq
  \frac12\log_2\left(\frac{8\pi e}{3\sqrt3}\right)
  =
  1.85837\ldots.
\]
In the opposite direction, Coxeter, Few, and Rogers~\cite{CFR1959} proved $\theta_L(B_2^n)\ge(e^{-3/2}+o(1))n$.
Together with the universal estimate of Ordentlich, Regev, and Weiss~\cite{ORW}, this gave the following bounds before the present work:
\[
  (e^{-3/2}+o(1))n
  \le
  \Theta_L^{\mathrm{sym}}(n)
  \le
  \Theta_L^{\mathrm{conv}}(n)
  \le
  Cn^2.
\]
These estimates left a multiplicative gap of order $n$; in particular, even the polynomial order of the worst-case lattice covering density was unknown.

\subsection{Main results}

The results of this paper narrow this polynomial gap to an iterated-logarithmic factor and place both extremal quantities at the $n\log n$ scale in the following precise sense.

\begin{theorem}\label{thm:main}
There exist absolute constants $c,C>0$ and a nonnegative sequence $\varepsilon_n\to0$ such that, for all sufficiently large $n$,
\[
  c n\log n
  \le
  \Theta_L^{\mathrm{sym}}(n)
  \le
  \Theta_L^{\mathrm{conv}}(n)
  \le
  Cn\log n\,(\log\log n)^{10/3+\varepsilon_n}.
\]
\end{theorem}

\Cref{thm:main} follows from the following two independent results.
The first gives a universal upper bound for every convex body.

\begin{theorem}\label{thm:upper}
There exist an absolute constant $C>0$ and a nonnegative sequence $\varepsilon_n\to0$ such that, for all sufficiently large $n$,
\[
  \Theta_L^{\mathrm{conv}}(n)
  \le
  Cn\log n\,(\log\log n)^{10/3+\varepsilon_n}.
\]
\end{theorem}

\begin{remark}
The exponent $10/3$ of $\log\log n$ in \cref{thm:upper} can be reduced to $\frac23\log_2(2\pi e)=2.72946\ldots$ by the two-stage binary correction developed in our companion note~\cite{LiLiuTwoStage2026}.
Our main goal in the present paper is to identify the $n\log n$ scale up to an iterated-logarithmic factor, rather than to optimize the exponent of $\log\log n$.
\end{remark}

The second constructs origin-symmetric bodies with large lattice covering density.

\begin{theorem}\label{thm:lower}
There is an absolute constant $c>0$ such that, for every sufficiently large integer $n$, there is an origin-symmetric convex body $K_n\subseteq\R^n$ such that
\[
  \theta_L(K_n)\ge c n\log n.
\]
\end{theorem}

The witness $K_n$ constructed in the proof of \cref{thm:lower} is already in John position, with John ellipsoid $B_2^n$.

Since $(\log\log n)^{10/3+o(1)}=(\log n)^{o(1)}$, \cref{thm:main} gives
\[
  \Theta_L^{\mathrm{sym}}(n)
  =
  n\log n\,(\log n)^{o(1)}
  \qquad\text{and}\qquad
  \Theta_L^{\mathrm{conv}}(n)
  =
  n\log n\,(\log n)^{o(1)}.
\]
Thus the worst-case lattice covering density is determined up to a factor $(\log\log n)^{10/3+o(1)}$.
The upper bound also confirms and strengthens the conjecture of Schymura, Wang, and Xue~\cite[Remark~1.4]{SWX} that $\theta_L(K)\le Cn(\log n)^{1+\log_2 e}$ for every $n$-dimensional convex body $K$.

For Euclidean balls, the universal upper bound asymptotically improves the factor $(\log n)^{\beta_{\mathrm{ball}}}$ obtained by Gao, Liu, Pikhurko, and Sun~\cite{GLPS} to $\log n\,(\log\log n)^{10/3+o(1)}$.

\subsection{Ideas behind the upper bound}

Choose an auxiliary dimension $d$ much smaller than $n$, write $n=m+d$, and split $\R^n=E\oplus F$ with $\dim E=m$ and $\dim F=d$.
For $y\in F$, let $K_y\coloneqq\{x\in E\colon x+y\in K\}$ be the horizontal section of $K$ at height $y$.
The vertical space $F$ is used to find a useful family of such sections above every translate, and the horizontal space $E$ is used to cover them.
The challenge is to do this when the sections have different volumes and then encode all vertical and horizontal shifts in a single lattice.

To find useful heights, we fold a narrow Gaussian into the unit cube.
Each base point determines Gaussian weights on its Boolean vertices.
The folding identity separates the total weight of the cube from the way this weight is divided among its vertices.
Entropy estimates show that no single vertex receives too much weight and that vertices with very small weights can later be removed while retaining almost all of the total.

Good base points must be available above every vertical translate.
We select a movable information layer on which the total weight and both entropy bounds are controlled, and divide the unit cube into a finite-field grid.
The high-rank Kakeya theorem of Kopparty, Lev, Saraf, and Sudan~\cite{KLSS2011} gives a linear subspace such that every grid cell can be shifted into this layer.
Construction A lifts the subspace to a lattice in $F$.
Every vertical translate then contains a Boolean cube with total normalized Gaussian weight $d^{-5/2+o(1)}$ and the required entropy bounds.

We next place $K$ in isotropic position.
The pointwise marginal theorem of Eldan and Klartag~\cite{EldanKlartag} supplies a low-dimensional subspace $F$ on which the section volumes satisfy
\[
  \vol_m(K_y)
  =
  \bigl(1+o(1)\bigr)\vol_n(K)\gamma_d(y)
\]
throughout the relevant range.
A finite discretization replaces these sections by finitely many slightly smaller convex bodies while retaining their individual volumes.
The mean-hole estimates of Rogers~\cite{Rogers1958} and Schmidt~\cite{Schmidt1958} provide one horizontal lattice that works for every retained section.
Shannon entropy controls the total load lost in pruning, while min-entropy keeps every remaining section within the range of the mean-hole estimate.

Finally, simultaneous horizontal shifts multiply the remaining uncovered proportions, and a finite-index shear incorporates all shifts into one lattice in $E\oplus F$.
A total horizontal load of order $n\log n$ makes the uncovered proportion small enough for Rogers' completion lemma~\cite{Rogers1959}, which turns the almost covering into an exact one.
The vertical construction costs $d^{5/2+o(1)}$.
The min-entropy requirement leads to
\[
  d=(\log\log n)^{4/3+o(1)},
\]
so the vertical loss is $(\log\log n)^{10/3+o(1)}$.

\subsection{Ideas behind the lower bound}

The lower-bound construction intersects a Euclidean ball and a cube with a Poisson family of slabs.
In the normalization used later, $a=2$, $b=ae^{1/n}$, and
\[
  K_\Pi
  =
  bB_2^n\cap[-1,1]^n\cap
  \bigcap_{u\in\Pi}
  \{x\in\R^n\colon|\ip{u}{x}|\le1\}.
\]
The expected volume of this body is comparable to $\kappa_n a^n$, and a fixed positive proportion of the realizations have volume of this order.
The Euclidean ball gives a uniform outer-radius cutoff, the random slabs create the obstructions to covering, and the cube places the construction directly in John position.
Indeed, the inclusions $B_2^n\subseteq K_\Pi\subseteq[-1,1]^n$, together with the cube's contact points, show that the John ellipsoid is exactly $B_2^n$.

Fix a lattice $\Lambda$ and a target $y\in\R^n/\Lambda$.
Choose a representative $\widetilde y\in\R^n$.
Each residual vector $\widetilde y-\lambda$, with $\lambda\in\Lambda$, determines a spherical cap of slab normals that remove it.
For the targets used in the proof, all relevant residuals lie in a thin shell.
The target is uncovered when every corresponding cap is hit by the Poisson process, so the problem becomes one of constructing and controlling simultaneous Poisson witnesses.

The caps belonging to one target can be strongly dependent.
We cluster nearby caps and place a measurable slot in the common intersection of each cluster.
An injective choice of one Poisson point from every slot witnesses that the target is a hole.
A hierarchical common-lens estimate controls the activity of these witnesses and their self-interactions.

A greedy packing on the flat torus supplies
\[
  T
  =
  \exp\bigl((1/2-o(1))n\log n\bigr)
\]
targets whose residual directions are separated across different targets.
The packing remains valid in the cusps of lattice space because the image of a Euclidean ball in a flat torus has Haar measure at most its Euclidean volume divided by $\det\Lambda$.
A continuous Poisson version of Janson's inequality then makes the covering probability superpolynomially small for each fixed lattice.
Finally, every lattice that could cover the random body has bounded covering radius and therefore a short basis.
A stable net over these basis matrices makes the estimate uniform over all lattices.

The same $n\log n$ scale emerges from different balances in the two proofs.
In the upper-bound argument, exponential decay in the total horizontal load means that a sufficiently large constant multiple of $n\log n$ reaches the $(n^n+1)^{-1}$ threshold required by Rogers' completion lemma.
In the lower-bound argument, the target supply above absorbs the self-interaction cost $e^{O(R)}$ of the Poisson witnesses throughout $R\le c n\log n$, where $R\coloneqq\kappa_n a^n/\det\Lambda$.
Excluding every lattice in this range forces any covering lattice to have determinant smaller than $\kappa_n a^n/(c n\log n)$, which yields the lower bound.

The remaining gap between \cref{thm:upper,thm:lower} is the iterated-logarithmic factor.

\subsection{Organization}

The proof has two technically independent parts.
Part~I proves the universal upper bound.
In \cref{sec:fc-external-results}, we fix the notation and state the external results used in that argument.
In \cref{sec:fc-weighted-cube}, we develop the weighted full-cube construction.
In \cref{sec:fc-weighted-lifting}, we transfer its weights to sections of an arbitrary convex body and perform the finite-index lifting.
In \cref{sec:fc-upper-proof}, we combine these ingredients to prove \cref{thm:upper}.

Part~II proves the worst-case lower bound.
In \cref{lb:sec:prelim}, we collect the spherical-cap, flat-torus, and Poisson tools.
In \cref{lb:sec:body}, we define the random body and prove its volume and normalization properties.
In \cref{lb:sec:fixed}, we construct good separated targets for a fixed lattice.
In \cref{lb:sec:clusters}, we develop the common-lens Poisson witness method and prove the fixed-lattice estimate.
In \cref{lb:sec:net}, we pass from one lattice to all lattices.
Finally, in \cref{lb:sec:lower-main-proof}, we prove \cref{thm:lower} and complete the proof of \cref{thm:main}.

\part{Universal upper bounds}

\section{Tools for the upper bound}\label{sec:fc-external-results}

This section fixes the quotient-torus notation and states the four principal external inputs used in the proof: the Rogers--Schmidt mean-hole estimate~\cite{Rogers1958,Schmidt1958}, Rogers' completion lemma~\cite{Rogers1959}, the pointwise Gaussian marginal theorem of Eldan and Klartag~\cite{EldanKlartag}, and the high-rank finite-field Kakeya bound of Kopparty, Lev, Saraf, and Sudan~\cite{KLSS2011}.

For a positive integer $N$, write $[N]\coloneqq\{1,\ldots,N\}$.
For $R>0$, let $B_R^r\coloneqq\{x\in\R^r\colon\lVert x\rVert_2\le R\}$.
If $\Lambda\subseteq\R^r$ is a full-rank lattice, let $\T_\Lambda\coloneqq\R^r/\Lambda$ be the \emph{quotient torus}.
Let $\mu_\Lambda$ be the normalized Lebesgue measure on $\T_\Lambda$, equivalently its Haar probability measure.
Thus $\mu_\Lambda$ is the unique translation-invariant probability measure on $\T_\Lambda$.
Let $\pi_\Lambda\colon\R^r\to\T_\Lambda$ be the quotient map.
For a Borel set $J\subseteq\R^r$, define the \emph{uncovered proportion}
\[
  \bar\rho(\Lambda+J) \coloneqq 1-\mu_\Lambda(\pi_\Lambda(J)).
\]
A full-rank lattice $L\subseteq\R^r$ is \emph{unimodular} if $\det L=1$.
Let $\cL_r$ denote the space of unimodular lattices in $\R^r$, and let $\mu_r$ be the normalized $\operatorname{SL}_r(\R)$-invariant probability measure on $\cL_r$, with the normalization used in Siegel's mean-value theorem~\cite{Siegel1945}.
Set $b_0\coloneqq\frac14\log\frac{27}{16}$ and $\eta_r\coloneqq b_0r-3\log r$.

Let $\operatorname{GL}_r(\R)$ denote the group of invertible $r\times r$ real matrices, acting as linear maps on $\R^r$.
If $A\in\operatorname{GL}_r(\R)$, then $\theta_L(A(K))=\theta_L(K)$ and $\bar\rho(A(\Lambda)+A(J))=\bar\rho(\Lambda+J)$.
Translations of the body leave both quantities unchanged.
We shall therefore use affine normalizations and determinant-one shears without further comment.

The following formulation of the mean-hole estimate is recorded by Ordentlich, Regev, and Weiss~\cite[Theorem~2.3]{ORW}.
It originates in Rogers' mean-value argument~\cite{Rogers1958}; see also Schmidt's closely related estimates~\cite{Schmidt1958}.

\begin{theorem}\label{thm:fc-RS}
There exists an absolute constant $c_{\mathrm{RS}}>0$ such that, for every sufficiently large $r$ and every Borel set $J\subseteq\R^r$ of volume $V\le\eta_r$,
\[
  \left| \int_{\cL_r}\bar\rho(L+J)\,\dd\mu_r(L) -e^{-V} \right| < c_{\mathrm{RS}}e^{-\eta_r}.
\]
\end{theorem}

The theorem applies to arbitrary Borel sets and hence to all labels for sections of varying volume used later.
More precisely, if the desired lattice determinant is $D_h$, scaling to determinant one replaces a set $J$ by $D_h^{-1/r}J$, whose volume is $\vol_r(J)/D_h$.
This normalized volume is the load used in the horizontal argument.
The upper restriction $V\le\eta_r$ is the reason for controlling the largest load at a Boolean vertex by min-entropy.

The last step of the argument uses Rogers' completion lemma~\cite[Lemma~4]{Rogers1959} to turn a lattice with sufficiently small uncovered proportion into an exact covering.

\begin{lemma}\label{lem:fc-rogers-completion}
Let $C\subseteq\R^r$ be a convex body and let $\Lambda\subseteq\R^r$ be a lattice.
If $\bar\rho(\Lambda+C) \le (r^r+1)^{-1}$, then $\Lambda+\left(1+\frac1r\right)C=\R^r$.
\end{lemma}

We next record the pointwise Gaussian marginal estimate of Eldan and Klartag~\cite[Theorem~1]{EldanKlartag}.
Its role is to compare, simultaneously at every point of a low-dimensional ball, the volume of a section of an isotropic convex body with the standard Gaussian density.

For $y\in\R^d$, write $\gamma_d(y)\coloneqq(2\pi)^{-d/2}e^{-\lVert y\rVert_2^2/2}$.
We use the pointwise, rather than merely distributional, form of the central limit theorem for low-dimensional marginals of log-concave measures.
Write $G_{n,d}$ for the Grassmannian of $d$-dimensional linear subspaces of $\R^n$, equipped with its rotation-invariant probability measure.

For a random vector $X$ with log-concave density $f_X$, choose the upper-semicontinuous log-concave representative of $f_X$.
If $F\in G_{n,d}$, $E\coloneqq F^\perp$, and $P_F$ denotes orthogonal projection onto $F$, define the \emph{$F$-marginal density} of $X$ by
\[
  \rho_F(y) \coloneqq \int_E f_X(x+y)\,\dd x, \quad\text{for }y\in F.
\]
This is the marginal operator used by Eldan and Klartag.
By Pr\'ekopa's theorem, $\rho_F$ is log-concave and hence continuous on the interior of its support.

\begin{theorem}\label{thm:fc-EK}
There exist universal constants $C,c_1,c_2,c_3,c_4>0$ with the following property.
Let $X$ be an isotropic random vector in $\R^n$ with a log-concave density, and let $1\le d\le n^{c_1}$ be an integer.
Then there is a set $\mathcal E\subseteq G_{n,d}$ of measure at least $1-Ce^{-n^{c_2}}$ such that, for every $F\in\mathcal E$, the density $\rho_F$ defined above satisfies
\begin{equation}\label{eq:fc-EK}
  \left| \frac{\rho_F(y)}{\gamma_d(y)}-1 \right| \le Cn^{-c_3}
\end{equation}
whenever $y\in F$ and $\lVert y\rVert_2\le n^{c_4}$.
\end{theorem}

Suppose that $X$ is uniform on an isotropic convex body $K$, and fix $f_X\coloneqq\1_K/\vol_n(K)$.
Let $E\coloneqq F^\perp$ and $m\coloneqq n-d$.
For $y\in F$, define the \emph{horizontal section} $K_y\coloneqq\{x\in E\colon x+y\in K\}$.
The defining integral now gives
\[
  \rho_F(y) = \frac{\vol_m(K_y)}{\vol_n(K)}
\]
for every $y\in F$.
On $\operatorname{int}(P_F(K))$, the Brunn--Minkowski inequality~\cite[Section~7.1]{Schneider2014} also shows directly that $y\mapsto\vol_m(K_y)^{1/m}$ is continuous.

The topological support of the projected measure is exactly $P_F(K)$.
The inclusion of the support in $P_F(K)$ follows from compactness.
Conversely, every relative neighborhood of a point of $P_F(K)$ meets $\operatorname{int}(P_F(K))$, where the section-volume function is positive, and therefore has positive projected measure.
In particular, $\rho_F$ is positive and continuous on $\operatorname{int}(P_F(K))$.
Thus, once $F$ is chosen as in \cref{thm:fc-EK}, the theorem compares $\vol_m(K_y)/\vol_n(K)$ with $\gamma_d(y)$ at every $y\in F$ satisfying $\lVert y\rVert_2\le n^{c_4}$.
This pointwise conclusion, rather than merely a distributional approximation, is what the vertical construction requires.

The fourth principal input collected here is a high-rank finite-field Kakeya theorem.

The polynomial method for finite-field Kakeya sets began with Dvir's rank-one theorem~\cite{Dvir2009}.
Dvir, Kopparty, Saraf, and Sudan~\cite{DKSS2013} later developed the multiplicity method and several extensions.
Here we use the rank-$k$ theorem of Kopparty, Lev, Saraf, and Sudan~\cite[Theorem~1]{KLSS2011}, stated below in the form we need.

Let $q$ be a prime power, let $1\le k\le d-1$, and let $\operatorname{Gr}_{d,k}(\F_q)$ denote the set of $k$-dimensional linear subspaces of $\F_q^d$.
A set $A\subseteq\F_q^d$ is a \emph{Kakeya set of rank $k$} if, for every $C\in\operatorname{Gr}_{d,k}(\F_q)$, the set $A$ contains a translate of $C$.

\begin{theorem}\label{thm:fc-KLSS}
If $A\subseteq\F_q^d$ is a Kakeya set of rank $k$, then
\[
  |A| \ge \left(1+(q-1)q^{-k}\right)^{-d}q^d.
\]
\end{theorem}

The form needed in the vertical construction is the following direct contrapositive.

\begin{corollary}\label{cor:fc-KLSS-cover}
Define the \emph{rank-$k$ Kakeya threshold} by
\begin{equation}\label{eq:fc-KLSS-threshold}
  \vartheta_{d,q,k} \coloneqq 1-\left(1+(q-1)q^{-k}\right)^{-d}.
\end{equation}
If $B\subseteq\F_q^d$ satisfies $|B|/q^d>\vartheta_{d,q,k}$, then there exists $C\in\operatorname{Gr}_{d,k}(\F_q)$ such that $B+C=\F_q^d$.
\end{corollary}

\begin{proof}
Suppose instead that $B+C\ne\F_q^d$ for every $C\in\operatorname{Gr}_{d,k}(\F_q)$.
For each such $C$, some affine coset of $C$ is disjoint from $B$ and hence is contained in $A\coloneqq\F_q^d\setminus B$.
Thus $A$ is a Kakeya set of rank $k$.
By \cref{thm:fc-KLSS}, $|A|\ge\left(1+(q-1)q^{-k}\right)^{-d}q^d$, contrary to the assumed lower bound on $|B|$.
\end{proof}

We shall also use the elementary estimate $\vartheta_{d,q,k}\le d(q-1)q^{-k}$, which follows from $(1+x)^{-d}\ge1-dx$ for $x\ge0$.
Consequently, whenever the matching construction produces a set $B$ with $\vartheta_{d,q,k}<|B|/q^d\le3\vartheta_{d,q,k}$, the subspace in \cref{cor:fc-KLSS-cover} has \emph{finite-field covering load} $|B||C|/q^d=(|B|/q^d)q^k\le3d(q-1)$.
The layer-to-grid matching that produces such a set $B$ is proved in \cref{prop:fc-layer-matching}.

\section{The weighted full-cube construction}\label{sec:fc-weighted-cube}

We now construct the vertical lattice and its Boolean weights.
No Boolean vertex is discarded at this stage.
The folded Gaussian identity separates the total weight of the full cube from its distribution among the vertices.
Min-entropy controls the largest vertex weight, while Shannon entropy later controls the total weight discarded when vertices with small weights are removed.

All constants in this section are absolute.
Every $o(1)$ term is a deterministic function of $d$ alone.
We let $d\to\infty$ and put
\begin{equation}\label{eq:fc-tau}
  \tau\coloneqq\frac1{2\log d}.
\end{equation}
Let $\gamma_\tau(u)\coloneqq(2\pi\tau)^{-1/2}e^{-u^2/(2\tau)}$ be the one-dimensional Gaussian density of variance $\tau$, and set $P_\tau\coloneqq\int_{-1}^1\gamma_\tau(u)\,\dd u$.

\subsection{Folded Gaussian weights and entropy}

Set
\begin{align*}
  Z\coloneqq\int_{-1}^1e^{-u^2/(2\tau)}\,\dd u=\sqrt{2\pi\tau}\,P_\tau
  \quad\text{and}\quad
  g(u)\coloneqq Z^{-1}e^{-u^2/(2\tau)}\1_{[-1,1]}(u).
\end{align*}
Thus $g$ is the density of $N(0,\tau)$ conditioned on $[-1,1]$.
Let $X$ have density $g$, and define $\mathsf E\coloneqq\1_{\{X>0\}}$ and $T\coloneqq\mathsf E-X$.
Then $\mathsf E\in\{0,1\}$, $T\in[0,1]$, and $X=\mathsf E-T$.
The density of $T$ is $f$ below.
Since a regular conditional probability is determined only almost everywhere, we choose the following continuous version of its Bernoulli parameter:
\begin{equation}\label{eq:fc-f-p}
  f(t)\coloneqq g(-t)+g(1-t), \quad\text{and}\quad p(t)\coloneqq\frac{g(1-t)}{f(t)} =\Prb\left[\mathsf E=1\mid T=t\right].
\end{equation}
For $t=(t_1,\ldots,t_d)\in[0,1]^d$, let $P_t$ be the product probability measure on $\{0,1\}^d$ given by
\[
  P_t(\varepsilon) \coloneqq \prod_{i=1}^d p(t_i)^{\varepsilon_i} (1-p(t_i))^{1-\varepsilon_i}.
\]
For a one-dimensional density $h$, write $h^{\otimes d}(x)\coloneqq\prod_{i=1}^d h(x_i)$ for its $d$-fold product density, where $x=(x_1,\ldots,x_d)\in\R^d$.

Summing the next pointwise factorization over the Boolean cube gives the folded-Gaussian identity used throughout the vertical construction.

\begin{lemma}\label{lem:fc-full-cube-identity}
For every $t\in[0,1]^d$ and $\varepsilon\in\{0,1\}^d$,
\begin{equation}\label{eq:fc-pointwise-factorization}
  g^{\otimes d}(\varepsilon-t) = f^{\otimes d}(t)P_t(\varepsilon).
\end{equation}
Consequently,
\begin{equation}\label{eq:fc-full-cube-identity}
  \sum_{\varepsilon\in\{0,1\}^d} \gamma_\tau^{\otimes d}(\varepsilon-t) = P_\tau^d f^{\otimes d}(t).
\end{equation}
Moreover,
\begin{equation}\label{eq:fc-Ptau}
  1-P_\tau = O\left(\frac1{d\sqrt{\log d}}\right), \quad\text{and}\quad P_\tau^d=1-o(1).
\end{equation}
\end{lemma}

\begin{proof}
The chosen conditional version gives $g(e-t)=f(t)\Prb\left[\mathsf E=e\mid T=t\right]$ for every $t\in[0,1]$ and $e\in\{0,1\}$.
Taking products proves \eqref{eq:fc-pointwise-factorization}, and summing over the Boolean cube and using $\sum_\varepsilon P_t(\varepsilon)=1$ gives $\sum_\varepsilon g^{\otimes d}(\varepsilon-t)=f^{\otimes d}(t)$.
Since $g(u)=P_\tau^{-1}\gamma_\tau(u)$ for $|u|\le1$, this proves \eqref{eq:fc-full-cube-identity}.

Finally, Mills' inequality~\cite{Gordon1941} at $\tau^{-1/2}=\sqrt{2\log d}$ gives the first estimate in \eqref{eq:fc-Ptau}; the second follows from $d(1-P_\tau)=O((\log d)^{-1/2})$.
\end{proof}

The next estimates determine the entropy scale used first to choose the finite-field rank and later to choose the vertical dimension.

The \emph{binary entropy function} is $h_2(u)\coloneqq-u\log u-(1-u)\log(1-u)$ for $u\in[0,1]$.
Define the \emph{one-coordinate Shannon entropy} and \emph{one-coordinate min-entropy}, respectively, by
\begin{equation}\label{eq:fc-one-coordinate-entropies}
  h_{\mathrm{Ber}}(t)\coloneqq h_2(p(t)), \quad\text{and}\quad h_{\mathrm{Ber},\infty}(t) \coloneqq -\log\max\{p(t),1-p(t)\}.
\end{equation}
For $T$ with density $f$, put $s\coloneqq\E h_{\mathrm{Ber}}(T)$ and $s_\infty\coloneqq\E h_{\mathrm{Ber},\infty}(T)$.

\begin{lemma}\label{lem:fc-entropy-asymptotics}
As $d\to\infty$,
\begin{align}
  s &= (2\sqrt\pi+o(1)) \frac{d^{-1/4}}{\sqrt{\log d}}, \label{eq:fc-s-asymptotic}\\
  s_\infty &= \left(\frac{2\pi-4}{\sqrt\pi}+o(1)\right) \frac{d^{-1/4}}{\sqrt{\log d}}. \notag
\end{align}
In particular,
\begin{equation}\label{eq:fc-entropy-scales}
  ds\asymp ds_\infty \asymp \frac{d^{3/4}}{\sqrt{\log d}} \gg\log d.
\end{equation}
\end{lemma}

\begin{proof}
Put $t\coloneqq\frac12+\tau y$.
Directly from \eqref{eq:fc-f-p},
\begin{equation}\label{eq:fc-logistic}
  p(t)=\frac{e^y}{1+e^y}, \quad\text{and}\quad h_{\mathrm{Ber},\infty}(t)=\log(1+e^{-|y|}).
\end{equation}
The same substitution, including the Jacobian $\dd t=\tau\,\dd y$, gives
\[
  f\left(\frac12+\tau y\right)\tau = \frac{2\tau e^{-1/(8\tau)}}{Z} e^{-\tau y^2/2}\cosh(y/2).
\]
For $\psi(y)\coloneqq h_2(e^y/(1+e^y))$ and $\psi_\infty(y)\coloneqq\log(1+e^{-|y|})$, we therefore have
\begin{equation}\label{eq:fc-entropy-integrals}
  \E\psi_\star\left(\frac{T-1/2}{\tau}\right)
  =
  \frac{2\tau e^{-1/(8\tau)}}{Z}
  \int_{-1/(2\tau)}^{1/(2\tau)}
  e^{-\tau y^2/2}\cosh(y/2)\psi_\star(y)\,\dd y,
\end{equation}
for $\psi_\star\in\{\psi,\psi_\infty\}$.
Extend each integrand by zero outside the interval $[-1/(2\tau),1/(2\tau)]$.
The resulting functions on $\R$ are bounded by the fixed integrable function $C(1+|y|)e^{-|y|/2}$, uniformly in $\tau$.
Dominated convergence therefore reduces the leading constants to two one-dimensional integrals, which we now evaluate explicitly.

Both $\psi$ and $\psi_\infty$ are even.
For $y\ge0$, set $x\coloneqq e^{-y}$.
Then
\[
  \psi(y)
  =
  \log(1+x)-\frac{x\log x}{1+x},
  \quad\text{and}\quad
  \psi_\infty(y)=\log(1+x),
  \quad\text{and}\quad
  \cosh(y/2)=\frac{1+x}{2\sqrt x}.
\]
It follows that
\begin{align*}
  \int_{\R}\cosh(y/2)\psi(y)\,\dd y = A+B+C
  \quad\text{and}\quad
  \int_{\R}\cosh(y/2)\psi_\infty(y)\,\dd y = A+B,
\end{align*}
where
\[
  A\coloneqq\int_0^1x^{-3/2}\log(1+x)\,\dd x,
  \quad
  B\coloneqq\int_0^1x^{-1/2}\log(1+x)\,\dd x,
  \quad\text{and}\quad
  C\coloneqq-\int_0^1x^{-1/2}\log x\,\dd x.
\]
For $A$, perform integration by parts first on $[\varepsilon,1]$ and then let $\varepsilon\downarrow0$; the boundary term at $\varepsilon$ tends to zero.
This justification, followed by $x=u^2$ for $B$, gives
\begin{align*}
  A &= -2\log 2+2\int_0^1\frac{x^{-1/2}}{1+x}\,\dd x = \pi-2\log 2,\\
  B &= 2\int_0^1\log(1+u^2)\,\dd u
  = 2\left[ u\log(1+u^2)-2u+2\arctan u \right]_{u=0}^{u=1} = 2\log 2-4+\pi.
\end{align*}
Finally, differentiating $\int_0^1x^{a-1}\,\dd x=a^{-1}$ at $a=1/2$ gives $C=4$.
Consequently,
\begin{align*}
  \int_{\R}\cosh(y/2)\psi(y)\,\dd y =A+B+C=2\pi
  \quad\text{and}\quad
  \int_{\R}\cosh(y/2)\psi_\infty(y)\,\dd y =A+B=2\pi-4.
\end{align*}
Since $Z=\sqrt{2\pi\tau}(1+o(1))$, $e^{-1/(8\tau)}=d^{-1/4}$, and $\sqrt{2\tau/\pi}=1/\sqrt{\pi\log d}$, substituting these values into \eqref{eq:fc-entropy-integrals} proves the result.
\end{proof}

For a finite probability distribution $q$, its \emph{Shannon entropy} and \emph{min-entropy} are $H(q)\coloneqq-\sum_x q(x)\log q(x)$ and $H_\infty(q)\coloneqq-\log\max_x q(x)$, respectively.
For the product distribution above,
\[
  H(P_t)=\sum_{i=1}^d h_{\mathrm{Ber}}(t_i), \quad\text{and}\quad H_\infty(P_t)= \sum_{i=1}^d h_{\mathrm{Ber},\infty}(t_i).
\]

We will use these entropies pointwise on most base points.
The next lemma shows that both remain close to their means with probability at least $1-O(d^{-30})$.

\begin{lemma}\label{lem:fc-entropy-concentration}
There is an absolute constant $C_{\mathrm{ent}}>0$ such that, for all sufficiently large $d$,
\begin{align*}
  \Prb\left[ H(P_T)>ds+C_{\mathrm{ent}}\sqrt{ds\log d} \right] \le d^{-30}
  \quad\text{and}\quad
  \Prb\left[ H_\infty(P_T) < ds_\infty-C_{\mathrm{ent}}\sqrt{ds_\infty\log d} \right] \le d^{-30},
\end{align*}
where $T\coloneqq(T_1,\ldots,T_d)$ has density $f^{\otimes d}$.
\end{lemma}

\begin{proof}
Both functions in \eqref{eq:fc-one-coordinate-entropies} take values in $[0,\log 2]$.
Hence $\Var(h_{\mathrm{Ber}}(T))\le\E h_{\mathrm{Ber}}(T)^2\le(\log 2)s$, and similarly $\Var(h_{\mathrm{Ber},\infty}(T))\le(\log 2)s_\infty$.
For $s_\star\in\{s,s_\infty\}$, put $t_\star\coloneqq C_{\mathrm{ent}}\sqrt{ds_\star\log d}$.
Bernstein's inequality bounds the relevant upper or lower tail by
\[
  \exp\left( -\frac{t_\star^2}{2d(\log 2)s_\star+\frac23(\log 2)t_\star} \right).
\]
By \eqref{eq:fc-entropy-scales}, $ds_\star\gg\log d$, so $t_\star=o(ds_\star)$ and the exponent is at most $-cC_{\mathrm{ent}}^2\log d$.
A sufficiently large absolute choice of $C_{\mathrm{ent}}$ makes both probabilities at most $d^{-30}$.
\end{proof}

\subsection{A movable information layer}

For a probability density $\varphi$, its \emph{differential entropy} is $h(\varphi)\coloneqq-\int \varphi(x)\log \varphi(x)\,\dd x$ whenever this integral is finite.
Set $a(t)\coloneqq-\log f(t)$, so that $h(f)=\E a(T)$.

The Berry--Esseen theorem is a quantitative form of the central limit theorem: when the variance stays bounded away from zero and the third absolute moment is bounded, the distribution function of a normalized sum of $d$ independent copies differs from the standard normal distribution function by $O(d^{-1/2})$~\cite[Chapter~V]{Petrov1975}.
We use this estimate to show that a fixed-width information interval has probability of order $d^{-1/2}$, uniformly when its location is shifted by $O(\log d)$.
The next lemma records the moment bounds needed for this application.

\begin{lemma}\label{lem:fc-center-moments}
One has
\begin{align}
  h(f) &= \frac12\log\frac{\pi e}{\log d} +O\left(\frac{d^{-1/4}}{\sqrt{\log d}}\right), \label{eq:fc-hf}\\
  \Var(a(T))&=\frac12+o(1), \label{eq:fc-a-variance}\\
  \E|a(T)-h(f)|^3&\le C. \label{eq:fc-a-third}
\end{align}
\end{lemma}

\begin{proof}
We first bound the moments of the conditional-information term, then compare $a(T)$ with $U^2/2$ to obtain the variance and third moment, and finally compute the mean by the mixed entropy chain rule.

For $t\in[0,1]$ and $e\in\{0,1\}$, put $i(t,e)\coloneqq-\log\Prb\left[\mathsf E=e\mid T=t\right]$ and $j(t,e)\coloneqq-\log g(e-t)$.
The pointwise factorization gives $a(T)+i(T,\mathsf E)=j(T,\mathsf E)$, while the Gaussian form of $g$ gives
\begin{equation}\label{eq:fc-a-j-i}
  a(T) = \log Z+\frac{X^2}{2\tau}-i(T,\mathsf E).
\end{equation}

We first quantify the moments of the conditional-information term.
Let $Y\coloneqq(T-1/2)/\tau$, $I\coloneqq i(T,\mathsf E)$, and $p_y\coloneqq e^y/(1+e^y)$.
For $q\in\{1,2,3\}$, define $m_q(y)\coloneqq\E[I^q\mid Y=y]$.
For $y\ge0$, the logistic formula \eqref{eq:fc-logistic} gives
\begin{align*}
  m_q(y) = p_y\bigl(\log(1+e^{-y})\bigr)^q +(1-p_y)\bigl(\log(1+e^y)\bigr)^q
  \le e^{-qy}+e^{-y}(y+\log 2)^q \le C_q(1+y^q)e^{-y}.
\end{align*}
Interchanging the two Bernoulli outcomes shows that $m_q$ is even, so
\begin{equation}\label{eq:fc-information-moment-pointwise}
  m_q(y)\le C_q(1+|y|^q)e^{-|y|} \quad\text{for }y\in\R.
\end{equation}
The density computation in \eqref{eq:fc-entropy-integrals} now gives the exact representation
\[
  \E I^q
  =
  c_\tau
  \int_{-1/(2\tau)}^{1/(2\tau)}
    e^{-\tau y^2/2}\cosh(y/2)m_q(y)\,\dd y,
  \quad\text{where}\quad
  c_\tau\coloneqq\frac{2\tau e^{-1/(8\tau)}}{Z}.
\]
By \eqref{eq:fc-information-moment-pointwise}, the integral is at most
\[
  C_q\int_{\R}(1+|y|^q)e^{-|y|/2}\,\dd y<\infty.
\]
Moreover,
\[
  c_\tau
  \asymp
  \frac{d^{-1/4}}{\sqrt{\log d}}
  \asymp s
\]
by \eqref{eq:fc-Ptau}, \eqref{eq:fc-tau}, and \eqref{eq:fc-s-asymptotic}.
Hence
\begin{equation}\label{eq:fc-information-moments}
  \E I^q=O(s) \quad\text{for }q=1,2,3.
\end{equation}
For $q=1$, this expectation is the conditional entropy $H(\mathsf E\mid T)=s$.

We next record the required truncated-Gaussian moment bounds.
If $G\sim N(0,1)$, then $U\coloneqq X/\sqrt\tau$ has the law of $G$ conditioned on $|G|\le\alpha_d$, where $\alpha_d\coloneqq\tau^{-1/2}=\sqrt{2\log d}$.
For every fixed even $q\ge2$, integration by parts yields
\[
  \begin{aligned}
    \int_{\alpha}^{\infty}x^q\gamma_1(x)\,\dd x
    =
    \alpha^{q-1}\gamma_1(\alpha)
    +(q-1)\int_{\alpha}^{\infty}x^{q-2}\gamma_1(x)\,\dd x.
  \end{aligned}
\]
Induction in $q$, starting from Mills' inequality~\cite{Gordon1941}, therefore gives
\[
  \E\left[ |G|^q\1_{\{|G|>\alpha_d\}} \right] = O_q\left( \frac{(\log d)^{(q-1)/2}}d \right) \quad\text{for }q=2,4,6.
\]
Since $\Prb\left[|G|\le\alpha_d\right]=P_\tau=1-o(1)$, division by $P_\tau$ in the conditional moments gives
\begin{align}
  \E U^2 &= 1+O\left(\frac{\sqrt{\log d}}d\right), \label{eq:fc-U-second}\\
  \E U^4 &= 3+O\left(\frac{(\log d)^{3/2}}d\right), \label{eq:fc-U-fourth}\\
  \E U^6 &= 15+O\left(\frac{(\log d)^{5/2}}d\right). \notag
\end{align}

Set $V\coloneqq U^2/2$.
The additive constant in \eqref{eq:fc-a-j-i} does not affect variance, so $\Var(a(T))=\Var(V-I)$.
Equations \eqref{eq:fc-U-fourth} and \eqref{eq:fc-information-moments} imply
\[
  \lVert V-\E V\rVert_2=O(1), \quad\text{and}\quad \lVert I-\E I\rVert_2 \le2\lVert I\rVert_2 = O(\sqrt s).
\]
Expanding the variance and applying Cauchy--Schwarz to the covariance term, we obtain
\begin{align*}
  \left| \Var(V-I)-\Var(V) \right| &\le \Var(I)+2|\operatorname{Cov}(V,I)|\\
  &\le \lVert I-\E I\rVert_2^2 +2\lVert V-\E V\rVert_2 \lVert I-\E I\rVert_2
  = O(s)+O(\sqrt s) = o(1).
\end{align*}
On the other hand, \eqref{eq:fc-U-second} and \eqref{eq:fc-U-fourth} give $\Var(V)=\frac14(\E U^4-(\E U^2)^2)=\frac12+o(1)$.
Consequently, $\Var(a(T))=\frac12+o(1)$, as asserted in \eqref{eq:fc-a-variance}.

The sixth-moment estimate and \eqref{eq:fc-information-moments} give
\begin{align*}
  \lVert V-\E V\rVert_3 \le2\lVert V\rVert_3 = 2\left(\frac{\E|U|^6}{8}\right)^{1/3} = O(1)
  \quad\text{and}\quad
  \lVert I-\E I\rVert_3 \le2\lVert I\rVert_3 = O(s^{1/3}).
\end{align*}
Minkowski's inequality applied to $a(T)-h(f)=(V-\E V)-(I-\E I)$ now yields
\[
  \bigl(\E|a(T)-h(f)|^3\bigr)^{1/3} \le \lVert V-\E V\rVert_3 + \lVert I-\E I\rVert_3 = O(1),
\]
which proves \eqref{eq:fc-a-third}.

It remains to compute the mean.
By \eqref{eq:fc-Ptau}, $|\log P_\tau|=O(1/(d\sqrt{\log d}))$, while \eqref{eq:fc-U-second} controls the quadratic term.
Thus
\[
  h(g)
  =\log Z+\frac{\E U^2}{2}
  =
  \frac12\log(2\pi\tau)+\log P_\tau
  +\frac{\E U^2}{2}
  =
  \frac12\log(2\pi e\tau)
  {}+O\left(\frac{\sqrt{\log d}}d\right)
  =
  \frac12\log(2\pi e\tau)+o(s),
\]
where the last equality uses \eqref{eq:fc-s-asymptotic}.
Away from the null point $X=0$, the map $X\mapsto(T,\mathsf E)$ is one-to-one on each of the two half-intervals and has Jacobian one.
Splitting the entropy integral over these two branches gives
\[
  h(g)
  =
  -\int_0^1 f(t)
  \bigl[
    (1-p(t))\log\bigl(f(t)(1-p(t))\bigr)
    +p(t)\log\bigl(f(t)p(t)\bigr)
  \bigr]\,\dd t
  =
  h(f)+\int_0^1f(t)h_2(p(t))\,\dd t.
\]
This is the mixed discrete--continuous entropy chain rule in the present setting, so $h(g)=h(f)+H(\mathsf E\mid T)=h(f)+s$, and hence
\[
  h(f) = \frac12\log(2\pi e\tau)-s+o(s) = \frac12\log\frac{\pi e}{\log d} +O\left(\frac{d^{-1/4}}{\sqrt{\log d}}\right),
\]
which is \eqref{eq:fc-hf}.
\end{proof}

For the finite-field step, we need to move the information window by $O(\log d)$ without changing its order of size.
The moment bounds above and the Berry--Esseen theorem give this uniformly in the offset.

\begin{proposition}\label{prop:fc-movable-layer}
There are absolute constants $W,c_0,C_0,R_0>0$ such that the following holds for all sufficiently large $d$.
For every $r\in[-R_0\log d,R_0\log d]$, let the \emph{information layer} $\cG(r)$ be the set of all $t\in[0,1]^d$ satisfying
\begin{align}
  dh(f)+r &\le -\log f^{\otimes d}(t) \le dh(f)+r+W, \label{eq:fc-layer-information}\\
  H_\infty(P_t) &\ge ds_\infty-C_{\mathrm{ent}}\sqrt{ds_\infty\log d}, \label{eq:fc-layer-min}\\
  H(P_t) &\le ds+C_{\mathrm{ent}}\sqrt{ds\log d}. \label{eq:fc-layer-Shannon}
\end{align}
Then
\begin{equation}\label{eq:fc-layer-probability}
  \frac{c_0}{\sqrt d} \le \int_{\cG(r)}f^{\otimes d}(t)\,\dd t \le \frac{C_0}{\sqrt d},
\end{equation}
and, writing $v(r)\coloneqq\vol_d(\cG(r))$,
\begin{equation}\label{eq:fc-layer-volume}
  c_0d^{-1/2}e^{dh(f)+r} \le v(r) \le C_0d^{-1/2}e^{dh(f)+r}.
\end{equation}
\end{proposition}

\begin{proof}
Fix $R_0\coloneqq3$.
Let $A_d\coloneqq\sum_{i=1}^d a(T_i)$ and $\sigma_d^2\coloneqq\Var(a(T))$.
Although the distribution of $a(T)$ depends on $d$, the Berry--Esseen bound~\cite[Chapter~V]{Petrov1975} and \cref{lem:fc-center-moments} give
\[
  \sup_x
  \left|
    \Prb\left[
      \frac{A_d-dh(f)}{\sigma_d\sqrt d}\le x
    \right]-\Phi_{\mathrm{N}}(x)
  \right|
  \le
  \frac{C\,\E|a(T)-h(f)|^3}
       {\sigma_d^3\sqrt d}
  \le
  \frac{C_{\mathrm{BE}}}{\sqrt d}.
\]
Here $\Phi_{\mathrm{N}}$ is the standard normal distribution function.
The uniform third-moment bound and $\sigma_d^2\to1/2$ show that $C_{\mathrm{BE}}$ is absolute.
Since $|r|\le R_0\log d=o(\sqrt d)$, the standard normal density is $\gamma_1(0)+o(1)$ uniformly on the interval
\[
  \left[ \frac r{\sigma_d\sqrt d}, \frac{r+W}{\sigma_d\sqrt d} \right].
\]
Thus its probability is $(c_W+o(1))/\sqrt d$, where $c_W\coloneqq W/\sqrt\pi$.
Passing from a distribution function to an interval uses two endpoints, so the Berry--Esseen error is at most $2C_{\mathrm{BE}}/\sqrt d$.
Choose $W$ so that $c_W>4C_{\mathrm{BE}}$.
This gives uniform upper and lower bounds of order $d^{-1/2}$ for the event \eqref{eq:fc-layer-information}.
Removing the two exceptional events in \cref{lem:fc-entropy-concentration} changes its probability by at most $2d^{-30}$.
This proves \eqref{eq:fc-layer-probability}.

On $\cG(r)$, the reciprocal density $1/f^{\otimes d}(t)$ lies between $e^{dh(f)+r}$ and $e^{dh(f)+r+W}$.
Integrating this reciprocal with respect to $f^{\otimes d}(t)\,\dd t$ over $\cG(r)$ proves \eqref{eq:fc-layer-volume}, after changing the absolute constants.
\end{proof}

\subsection{Matching a layer to the Kakeya threshold}

We next match the continuous layer, the grid size, and the discrete rank.
The strict inequality required by \cref{cor:fc-KLSS-cover} must also survive integer rounding.

Bertrand's postulate states that, for every integer $N\ge2$, there is a prime strictly between $N$ and $2N$.
Applying it with $N=\lceil d(\log d)^2\rceil$, and using $2\lceil d(\log d)^2\rceil\le3d(\log d)^2$ for all sufficiently large $d$, fix once and for all a prime $Q=Q(d)$ satisfying
\begin{equation}\label{eq:fc-Q}
  d(\log d)^2\le Q\le3d(\log d)^2.
\end{equation}
Partition $[0,1)^d$ into the \emph{half-open grid cells}
\[
  \mathcal Q_b \coloneqq \frac bQ+[0,Q^{-1})^d, \quad\text{for }b\in\F_Q^d.
\]
Here each $b\in\F_Q^d$ is identified with its standard representative in $\{0,\ldots,Q-1\}^d$.

\begin{proposition}\label{prop:fc-layer-matching}
There are absolute constants $c,C>0$ such that, for every sufficiently large $d$ and the prime $Q=Q(d)$ fixed in \eqref{eq:fc-Q}, there are $r\in[-2\log Q,2\log Q]$, $k\in\{2,\ldots,d-1\}$, and $B\subseteq\F_Q^d$ such that every cell $\mathcal Q_b$ with $b\in B$ meets $\cG(r)$.
With $\beta\coloneqq|B|/Q^d$ and $\vartheta_k\coloneqq\vartheta_{d,Q,k}$, one has
\begin{equation}\label{eq:fc-matching-threshold}
  \vartheta_k<\beta\le3\vartheta_k
\end{equation}
and
\begin{equation}\label{eq:fc-matching-volume}
  c\,v(r)\le\beta\le C\,v(r).
\end{equation}
Consequently, there is a $k$-dimensional linear subspace $C_d\le\F_Q^d$ such that
\begin{equation}\label{eq:fc-B-plus-C}
  B+C_d=\F_Q^d
\end{equation}
and
\begin{equation}\label{eq:fc-finite-load}
  c\,dQ \le \beta Q^k = \frac{|B||C_d|}{Q^d} \le C\,dQ.
\end{equation}
\end{proposition}

\begin{proof}
We first choose the rank $k$, then shift the information layer to the corresponding Kakeya threshold, select enough intersecting grid cells, and finally apply \cref{cor:fc-KLSS-cover}.

Put $v_*\coloneqq d^{-1/2}e^{dh(f)}$ and $x_k\coloneqq d(Q-1)Q^{-k}$.
By \eqref{eq:fc-hf}, $\log v_*=-\frac d2\log\log d+O(d)$, whereas $x_2=d(Q-1)Q^{-2}\ge c/(\log d)^2$ and $x_{d-1}\le\exp(-(1+o(1))d\log d)$.
Thus, for all sufficiently large $d$, there is $k\in\{2,\ldots,d-1\}$ such that
\begin{equation}\label{eq:fc-k-choice}
  x_k\le v_*<Qx_k.
\end{equation}
In particular, $x_k\to0$.
Moreover, $\log(Qv_*)=-\frac d2\log\log d+O(d)\to-\infty$, so $Qx_k\le Qv_*\to0$.

Set $u_k\coloneqq(Q-1)Q^{-k}=x_k/d$.
Since $du_k=x_k\to0$, the definition \eqref{eq:fc-KLSS-threshold} and the expansion $1-(1+u_k)^{-d}=du_k(1+O(du_k))$, uniform for the selected $k$, give
\begin{equation}\label{eq:fc-theta-x}
  \vartheta_k=(1+o(1))x_k.
\end{equation}
Choose a sufficiently large absolute constant $M$, depending only on the lower constant $c_0$ in \eqref{eq:fc-layer-volume}, and put
\begin{equation}\label{eq:fc-r-choice}
  r\coloneqq \log\frac{M\vartheta_k}{v_*}.
\end{equation}
Equations \eqref{eq:fc-k-choice} and \eqref{eq:fc-theta-x} imply $-\log Q+O(1)\le r\le O(1)$, so $r\in[-2\log Q,2\log Q]$ for large $d$.
Since $\log Q=(1+o(1))\log d$, this lies in the range of \cref{prop:fc-movable-layer}, where $R_0=3$.
By \eqref{eq:fc-layer-volume} and \eqref{eq:fc-r-choice},
\begin{equation}\label{eq:fc-v-theta}
  c_0M\vartheta_k \le v(r) \le C_0M\vartheta_k.
\end{equation}
Fix $M$ so that $c_0M\ge8$.

Let $S\coloneqq\{b\in\F_Q^d\colon\mathcal Q_b\cap\cG(r)\ne\varnothing\}$.
The cells are disjoint and cover $[0,1)^d$.
Since $\cG(r)\setminus[0,1)^d$ lies in the boundary of the cube and hence has zero $d$-dimensional volume, the union of the cells indexed by $S$ has volume at least $v(r)$.
Therefore
\begin{equation}\label{eq:fc-intersecting-cells}
  \frac{|S|}{Q^d}\ge v(r)\ge8\vartheta_k.
\end{equation}
We also have $Q^{-d}\le\vartheta_k$ for all sufficiently large $d$, uniformly for the selected $k$.
Indeed, \eqref{eq:fc-theta-x} and $k\le d-1$ give $\vartheta_k Q^d\ge\frac12d(Q-1)Q^{d-k}\ge1$.
Choose any subset $B\subseteq S$ of cardinality
\begin{equation}\label{eq:fc-B-rounding}
  |B|=\left\lceil2\vartheta_kQ^d\right\rceil.
\end{equation}
The choice is possible by \eqref{eq:fc-intersecting-cells}, and the last observation gives
\[
  2\vartheta_k \le\beta \le2\vartheta_k+Q^{-d} \le3\vartheta_k.
\]
This proves the strict inequality and all rounding assertions in \eqref{eq:fc-matching-threshold}.
Combining this estimate with \eqref{eq:fc-v-theta} proves \eqref{eq:fc-matching-volume}.

Now apply \cref{cor:fc-KLSS-cover} to obtain a $k$-dimensional subspace $C_d$ satisfying \eqref{eq:fc-B-plus-C}.
Finally, \eqref{eq:fc-theta-x} gives $\beta Q^k\asymp x_kQ^k=d(Q-1)\asymp dQ$, which proves \eqref{eq:fc-finite-load}.
\end{proof}

\subsection{Construction A and the scaled weighted construction}

We lift the covering subspace $C_d$ to a Euclidean lattice so that every torus point can be moved into a selected grid cell.
The stability lemma then transfers the information-layer estimates from the chosen representative of that cell to the actual point.

For every $b\in B$, choose a representative
\begin{equation}\label{eq:fc-cell-representative}
  t_b\in\cG(r)\cap\mathcal Q_b.
\end{equation}
Let $\operatorname{red}_Q\colon\Z^d\to\F_Q^d$ be coordinatewise reduction modulo $Q$.
Define the \emph{unscaled Construction-A lattice} by
\[
  \Lambda_0
  \coloneqq
  \frac1Q\operatorname{red}_Q^{-1}(C_d)
  =
  \left\{
    \frac zQ\colon
    z\in\Z^d,\ \operatorname{red}_Q(z)\in C_d
  \right\}.
\]
Equivalently, if $\widehat c\in\{0,\ldots,Q-1\}^d$ denotes the standard integer representative of $c\in C_d$, then
\[
  \Lambda_0 = \Z^d+ \left\{\frac{\widehat c}{Q}\colon c\in C_d\right\}.
\]
Changing the integer representative of $c$ only adds an element of $\Z^d$, so the lattice is independent of this choice.
Since $\Lambda_0/\Z^d\cong C_d$, its covolume is
\begin{equation}\label{eq:fc-det-zero}
  \det\Lambda_0=|C_d|^{-1}=Q^{-k}.
\end{equation}
This is the usual Construction-A lift of a linear code, as used, for example, by Erez, Litsyn, and Zamir~\cite{ErezLitsynZamir2005}.

To pass from a selected point of a good grid cell to every point in that cell, we need the following stability estimate.

\begin{lemma}\label{lem:fc-cell-stability}
If $t,t'\in[0,1]^d$ and $\lVert t-t'\rVert_\infty\le Q^{-1}$, then
\begin{align}
  \left| \log f^{\otimes d}(t)-\log f^{\otimes d}(t') \right| &\le\Delta_d, \label{eq:fc-f-stability}\\
  |H(P_t)-H(P_{t'})| + |H_\infty(P_t)-H_\infty(P_{t'})| &\le C\Delta_d, \label{eq:fc-entropy-stability}
\end{align}
where
\begin{equation}\label{eq:fc-Delta}
  \Delta_d\coloneqq\frac{d}{Q\tau} = O\left(\frac1{\log d}\right).
\end{equation}
\end{lemma}

\begin{proof}
The two Gaussian summands defining $f$ give $|(\log f)'(u)|\le\tau^{-1}$ for $0<u<1$, with the corresponding one-sided bounds at the endpoints.
Thus the mean-value theorem proves \eqref{eq:fc-f-stability}.

In the logistic coordinate $y=(t-1/2)/\tau$, the function $h_2(e^y/(1+e^y))$ has bounded derivative.
The function $\log(1+e^{-|y|})$ is not differentiable at zero, but its two one-sided derivatives are bounded in absolute value; hence it is globally Lipschitz.
Applying these Lipschitz bounds coordinatewise proves \eqref{eq:fc-entropy-stability}.
The estimate \eqref{eq:fc-Delta} follows from \eqref{eq:fc-tau} and \eqref{eq:fc-Q}.
\end{proof}

Passing to standard-Gaussian coordinates scales $\Lambda_0$ by $\tau^{-1/2}$ and hence multiplies its covolume by $\tau^{-d/2}$.
Define the \emph{scaled vertical lattice} and the scaling map by
\begin{equation}\label{eq:fc-scaled-lattice}
  \Lambda_v \coloneqq \tau^{-1/2}\Lambda_0, \quad\text{and}\quad \Psi_d\coloneqq\tau^{-1/2}\operatorname{Id}_d.
\end{equation}
Then $\Psi_d(\Z^d)\subseteq\Lambda_v$ is a finite-index full-rank sublattice, and
\begin{equation}\label{eq:fc-scaled-determinant}
  \det\Lambda_v = \tau^{-d/2}\det\Lambda_0 = \tau^{-d/2}Q^{-k}.
\end{equation}

The next proposition is the output of the vertical construction.
It places a full Boolean pattern above every vertical translate and records its total Gaussian weight, largest normalized weight, and entropy.

\begin{proposition}\label{prop:fc-weighted-vertical}
For all sufficiently large $d$, the lattice $\Lambda_v$ defined in \eqref{eq:fc-scaled-lattice} has the following property.
For every $y\in\R^d$, there are $\lambda_y\in\Lambda_v$ and $\widetilde t_y\in[0,1)^d$ such that
\begin{equation}\label{eq:fc-vertex-position}
  y-\lambda_y-\Psi_d(\varepsilon) = \frac{\widetilde t_y-\varepsilon}{\sqrt\tau} \qquad \text{for every }\varepsilon\in\{0,1\}^d.
\end{equation}
Define the \emph{vertex weights} and their \emph{total weight} by
\[
  \mathsf w_{y,\varepsilon}
  \coloneqq
  \det\Lambda_v\,
  \gamma_d\left(
    y-\lambda_y-\Psi_d(\varepsilon)
  \right),
  \quad\text{and}\quad
  \mathsf W_y
  \coloneqq
  \sum_{\varepsilon\in\{0,1\}^d}
  \mathsf w_{y,\varepsilon}.
\]
Uniformly in $y$,
\begin{align}
  \mathsf W_y &= d^{-5/2+o(1)}, \label{eq:fc-total-weight}\\
  \frac{\max_\varepsilon\mathsf w_{y,\varepsilon}}
       {\mathsf W_y}
  &\le
  \exp\left(
    -ds_\infty
    +C\sqrt{ds_\infty\log d}
    +o(\log d)
  \right),
  \label{eq:fc-max-atom}\\
  -\sum_\varepsilon
  \frac{\mathsf w_{y,\varepsilon}}{\mathsf W_y}
  \log
  \frac{\mathsf w_{y,\varepsilon}}{\mathsf W_y}
  &\le
  ds+C\sqrt{ds\log d}+o(\log d).
  \label{eq:fc-weight-entropy}
\end{align}
Here the first estimate is uniform in the following explicit sense: there is an absolute constant $C_\xi>0$ such that one may take $\xi_d\coloneqq C_\xi\log\log d/\log d$ and obtain $d^{-5/2-\xi_d}\le\mathsf W_y\le d^{-5/2+\xi_d}$ for every $y$.
All points in \eqref{eq:fc-vertex-position} have norm at most
\begin{equation}\label{eq:fc-vertex-radius}
  \sqrt{\frac d\tau} = \sqrt{2d\log d}.
\end{equation}
\end{proposition}

\begin{proof}
We first use $B+C_d=\F_Q^d$ to move the given point into a selected cell, then derive exact formulas for the Boolean weights, and finally transfer the information-layer bounds by cell stability.

Fix $y\in\R^d$ and put $w\coloneqq\sqrt\tau\,y$.
Choose $z_0\in\Z^d$ such that $w-z_0\in[0,1)^d$, and set $w_0\coloneqq w-z_0$.
Let $z\in\F_Q^d$ be the index of the cell containing $w_0$, let $\widehat z$ be its standard integer representative, and write $w_0=\widehat z/Q+u$ with $u\in[0,Q^{-1})^d$.
By \eqref{eq:fc-B-plus-C}, write $z=b+c$ in $\F_Q^d$, where $b\in B$ and $c\in C_d$.
If $\widehat b$ and $\widehat c$ are the standard integer representatives, there is a $\kappa\in\Z^d$ such that $\widehat z=\widehat b+\widehat c-Q\kappa$.
Then
\begin{equation}\label{eq:fc-carry}
  \widetilde t_y \coloneqq w_0-\frac{\widehat c}Q+\kappa = \frac{\widehat b}Q+u \in\mathcal Q_b.
\end{equation}
Thus, with
\[
  \lambda_0 \coloneqq z_0-\kappa+\frac{\widehat c}Q \in\Lambda_0, \quad\text{and}\quad \lambda_y\coloneqq\tau^{-1/2}\lambda_0,
\]
one has $\widetilde t_y=w-\lambda_0$, which proves \eqref{eq:fc-vertex-position}.
Equation \eqref{eq:fc-vertex-radius} follows because every coordinate of $\widetilde t_y-\varepsilon$ lies in $[-1,1]$.

The determinant scaling in \eqref{eq:fc-scaled-determinant} gives the exact identity
\begin{align}
  \mathsf w_{y,\varepsilon} &= \det\Lambda_0\, \gamma_\tau^{\otimes d} (\widetilde t_y-\varepsilon), \notag\\
  \mathsf W_y &= \det\Lambda_0\,P_\tau^d f^{\otimes d}(\widetilde t_y), \label{eq:fc-exact-total-weight}\\
  \frac{\mathsf w_{y,\varepsilon}}{\mathsf W_y} &= P_{\widetilde t_y}(\varepsilon). \label{eq:fc-normalized-law}
\end{align}
Here the second and third lines follow from \cref{lem:fc-full-cube-identity}.
Thus the factor $\tau^{d/2}$ arising from the Gaussian scaling is incorporated exactly through \eqref{eq:fc-scaled-determinant}.

By \eqref{eq:fc-cell-representative} and \eqref{eq:fc-carry}, $\lVert\widetilde t_y-t_b\rVert_\infty\le Q^{-1}$.
Equations \eqref{eq:fc-layer-information} and \eqref{eq:fc-f-stability} therefore give
\[
  f^{\otimes d}(\widetilde t_y) = \exp\bigl( -dh(f)-r+O(W+\Delta_d) \bigr).
\]
On the other hand, \eqref{eq:fc-B-rounding}, \eqref{eq:fc-r-choice}, and $v_*=d^{-1/2}e^{dh(f)}$ give
\begin{equation}\label{eq:fc-beta-f}
  \beta f^{\otimes d}(\widetilde t_y) = d^{-1/2}\exp(O(W+\Delta_d)).
\end{equation}
Combining \eqref{eq:fc-finite-load}, \eqref{eq:fc-det-zero}, \eqref{eq:fc-beta-f}, and \eqref{eq:fc-exact-total-weight}, we obtain, uniformly in $y$,
\[
  \mathsf W_y = P_\tau^dQ^{-k}f^{\otimes d}(\widetilde t_y) = d^{-3/2}Q^{-1}e^{O(1)} = d^{-5/2}(\log d)^{-2}e^{O(1)}.
\]
Here we used $P_\tau^d=1-o(1)$ and $Q\asymp d(\log d)^2$.
This proves \eqref{eq:fc-total-weight}, and choosing $C_\xi$ sufficiently large gives the explicit two-sided bound in the statement.

Finally, \eqref{eq:fc-layer-min}, \eqref{eq:fc-layer-Shannon}, and \cref{lem:fc-cell-stability} yield
\begin{align*}
  H_\infty(P_{\widetilde t_y}) \ge ds_\infty -C_{\mathrm{ent}}\sqrt{ds_\infty\log d} -C\Delta_d
  \quad\text{and}\quad
  H(P_{\widetilde t_y}) \le ds+C_{\mathrm{ent}}\sqrt{ds\log d} +C\Delta_d.
\end{align*}
Since $\Delta_d=o(\log d)$, the exact normalized law \eqref{eq:fc-normalized-law} proves \eqref{eq:fc-max-atom} and \eqref{eq:fc-weight-entropy}.
\end{proof}

\section{Sections of varying volume and weighted lifting}\label{sec:fc-weighted-lifting}

The vertical construction assigns a different Gaussian weight to each Boolean vertex.
Equalizing the corresponding horizontal sections would lose this information.
In this section we retain the volume of each section, find one horizontal lattice for all loads in a prescribed window, amplify label-dependent hole bounds, and insert the resulting shifts into a finite-index shear.
The determinant of the final lattice is the product of the horizontal determinant and that of the scaled vertical lattice.

\subsection{Discretizing sections of varying volume}

Let $K\subseteq\R^n$ be a convex body.
After an invertible affine transformation, assume that the uniform probability measure on $K$ is isotropic.
Fix an integer $d<n$, put $m\coloneqq n-d$, and let $\R^n=E\oplus F$ be an orthogonal decomposition, with $\dim E=m$ and $\dim F=d$, supplied by \cref{thm:fc-EK}.
For $y\in F$, write $K_y\coloneqq\{x\in E\colon x+y\in K\}$.

For the uniform measure on $K$, the definition of $\rho_F$ gives $\rho_F(y)=\vol_m(K_y)/\vol_n(K)$.
For $F$ supplied by \cref{thm:fc-EK} and all sufficiently large $n$, so that $Cn^{-c_3}<1$, the lower bound in \eqref{eq:fc-EK} is positive at every point of the open ball $\{\lVert y\rVert_2<n^{c_4}\}$.
Each such point therefore belongs to the support of the marginal, which is $P_F(K)$.
Thus the entire open ball is contained in $P_F(K)$, and, because the ball is open, each of its points lies in $\operatorname{int}(P_F(K))$.
Consequently, throughout this range,
\begin{equation}\label{eq:fc-pointwise-section-volume}
  \vol_m(K_y) = \bigl(1+O(n^{-c_3})\bigr)\vol_n(K)\gamma_d(y).
\end{equation}

The horizontal sections vary continuously with $y$, but the later union bound requires a finite family.
The next proposition replaces all sections in the relevant ball by at most $(Cm^2)^d$ slightly smaller model sections while preserving their volumes up to a factor $1+O(m^{-1})$.

\begin{proposition}\label{prop:fc-variable-section-discretization}
There is an absolute constant $C>0$ with the following property.
Let $m\coloneqq n-d$, suppose that $m\ge2$ and $d\le n^{c_1}$, and let $R>0$ satisfy
\begin{equation}\label{eq:fc-section-discretization-range}
  2dR<n^{c_4}.
\end{equation}
For all sufficiently large $n$, every convex body $K\subseteq\R^n$ has an affine image $\widetilde K$ and an orthogonal decomposition $\R^n=E\oplus F$, with $\dim E=m$ and $\dim F=d$.
For $s>0$, write $B_s^F\coloneqq\{y\in F\colon\lVert y\rVert_2\le s\}$, and use $\gamma_d(y)\coloneqq(2\pi)^{-d/2}e^{-\lVert y\rVert_2^2/2}$ for the standard Gaussian density on $F$.
Define $\widetilde K_y\coloneqq\{x\in E\colon x+y\in\widetilde K\}$.
There are convex bodies $J_1,\ldots,J_N\subseteq E$, with $N\le(Cm^2)^d$, such that the following holds.
For every $y\in B_R^F$, there is an $a\in[N]$ satisfying
\begin{align}
  J_a&\subseteq \widetilde K_y, \label{eq:fc-variable-label-inclusion}\\
  \left(1-\frac1{m^2}\right)^{2m}
  \vol_m(\widetilde K_y)
  &\le
  \vol_m(J_a)
  \le
  \vol_m(\widetilde K_y).
  \label{eq:fc-variable-label-comparison}
\end{align}
In particular, the label $a$ can be chosen so that, uniformly in $y$,
\begin{equation}\label{eq:fc-variable-label-gaussian}
  \vol_m(J_a) = \left(1+O\left(\frac1m+n^{-c_3}\right)\right) \vol_n(\widetilde K)\gamma_d(y).
\end{equation}
\end{proposition}

\begin{proof}
We first center all sections over a containing simplex by a determinant-one shear, then choose a fine net in the target ball, and finally contract the net sections so that they lie in every nearby section.

Let $K_0$ be an affine image of $K$ whose uniform probability measure is isotropic.
Choose $F$ for $K_0$ from \cref{thm:fc-EK}, and set $E\coloneqq F^\perp$.
Let $\Sigma\subseteq F$ be a regular $d$-simplex centered at the origin and having inradius $2R$.
The circumradius of a regular $d$-simplex is $d$ times its inradius, so $B_{2R}^F\subseteq\Sigma\subseteq B_{2dR}^F$.
Since $2dR<n^{c_4}$, the support conclusion above gives $\Sigma\subseteq\operatorname{int}(P_F(K_0))$.
Moreover, \eqref{eq:fc-pointwise-section-volume} shows that every vertex section has positive $m$-dimensional volume.

Let $v_0,\ldots,v_d$ be the vertices of $\Sigma$, and choose $x_i\in (K_0)_{v_i}$.
Write
\[
  y=\sum_{i=0}^d\alpha_i(y)v_i,
  \quad\text{where}\quad
  \sum_{i=0}^d\alpha_i(y)=1,
\]
where the $\alpha_i(y)$ are the affine barycentric coordinates of $y$, and define $x_0(y)\coloneqq\sum_{i=0}^d\alpha_i(y)x_i$.
For $y\in\Sigma$, all the coefficients are nonnegative, and convexity gives $x_0(y)\in (K_0)_y$.

Apply the affine shear $\mathcal S(x,y)\coloneqq(x-x_0(y),y)$, and set $\widetilde K\coloneqq\mathcal S(K_0)$.
The linear part of $\mathcal S$ is block triangular with determinant one, and the shear preserves the vertical marginal and every section volume.
We therefore have
\begin{equation}\label{eq:fc-centered-sections}
  0\in \widetilde K_y \qquad\text{for every }y\in\Sigma.
\end{equation}

Choose an $R/m^2$-net $\mathcal Y\coloneqq\{y_1,\ldots,y_N\}\subseteq B_R^F$.
The standard packing estimate gives
\[
  N \le \left(1+\frac{2R}{R/m^2}\right)^d = (1+2m^2)^d \le (Cm^2)^d.
\]
Fix $y\in B_R^F$, and choose $y_a\in\mathcal Y$ with $\lVert y-y_a\rVert_2\le R/m^2$.
If $y=y_a$, the two inclusions below hold with $\sigma\coloneqq0$.
Otherwise, set
\[
  u\coloneqq\frac{y-y_a}{\lVert y-y_a\rVert_2},
  \quad\text{and}\quad
  \sigma\coloneqq\frac{\lVert y-y_a\rVert_2}{R}
  \le\frac1{m^2},
\]
and define $z_+\coloneqq y_a+Ru$ and $z_-\coloneqq y-Ru$.
Both $z_+$ and $z_-$ belong to $B_{2R}^F\subseteq\Sigma$, and
\[
  y=(1-\sigma)y_a+\sigma z_+, \quad\text{and}\quad y_a=(1-\sigma)y+\sigma z_-.
\]
Convexity of $\widetilde K$, together with \eqref{eq:fc-centered-sections}, gives, in either case,
\[
  (1-\sigma)\widetilde K_{y_a}\subseteq \widetilde K_y,
  \quad\text{and}\quad
  (1-\sigma)\widetilde K_y\subseteq \widetilde K_{y_a}.
\]

For every $a\in[N]$, define
\begin{equation}\label{eq:fc-variable-label-definition}
  J_a\coloneqq \left(1-\frac1{m^2}\right)\widetilde K_{y_a},
\end{equation}
where the contraction is about the origin.
By \eqref{eq:fc-centered-sections}, $0\in\widetilde K_{y_a}$.
For any convex set $C$ containing the origin, $rC\subseteq sC$ whenever $0\le r\le s$.
Since $1-m^{-2}\le1-\sigma$, we therefore have $J_a\subseteq(1-\sigma)\widetilde K_{y_a}\subseteq\widetilde K_y$, which proves \eqref{eq:fc-variable-label-inclusion}.
The second inclusion gives
\[
  \vol_m(\widetilde K_{y_a}) \ge (1-\sigma)^m\vol_m(\widetilde K_y) \ge \left(1-\frac1{m^2}\right)^m\vol_m(\widetilde K_y).
\]
Multiplying by the volume contraction in \eqref{eq:fc-variable-label-definition} proves the lower bound in \eqref{eq:fc-variable-label-comparison}; its upper bound follows from \eqref{eq:fc-variable-label-inclusion}.
Finally, $(1-m^{-2})^{2m}=1+O(m^{-1})$, and \eqref{eq:fc-pointwise-section-volume} proves \eqref{eq:fc-variable-label-gaussian}.
\end{proof}

\subsection{A common horizontal lattice and simultaneous shifts}

The next lemma applies the Rogers--Schmidt estimate to all model sections at once.
The sections may have different volumes; only their individual loads enter the union bound.

\begin{lemma}\label{lem:fc-varying-common-lattice}
Let $m$ be sufficiently large, let $D_h>0$, and let $J_a\subseteq\R^m$, $a\in\mathcal A$, be finitely many compact sets of positive volume.
Define the \emph{load} of $J_a$ by $u_a\coloneqq\vol_m(J_a)/D_h$, and suppose $u_a\le\eta_m$ for every $a\in\mathcal A$.
Let $0<c<1$.
If
\begin{equation}\label{eq:fc-varying-union-condition}
  (1+c_{\mathrm{RS}}) \sum_{a\in\mathcal A}e^{-(1-c)u_a}<1,
\end{equation}
then there exists a lattice $\Lambda_h\subseteq\R^m$, with $\det\Lambda_h=D_h$, such that
\begin{equation}\label{eq:fc-varying-hole}
  \bar\rho(\Lambda_h+J_a)\le e^{-cu_a} \qquad\text{for every }a\in\mathcal A.
\end{equation}
In particular, if $|\mathcal A|\le N$ and $T\le u_a\le\eta_m$ for every $a\in\mathcal A$, then the conclusion with $c=1/2$ holds whenever
\begin{equation}\label{eq:fc-varying-threshold}
  T\ge2\log\bigl(2(1+c_{\mathrm{RS}})N\bigr).
\end{equation}
\end{lemma}

\begin{proof}
For $a\in\mathcal A$, the rescaled set $\widetilde J_a\coloneqq D_h^{-1/m}J_a$ has volume $u_a$.
For a Haar--Siegel random unimodular lattice $L$, \cref{thm:fc-RS} and $u_a\le\eta_m$ give
\[
  \E_L\bar\rho(L+\widetilde J_a) \le e^{-u_a}+c_{\mathrm{RS}}e^{-\eta_m} \le (1+c_{\mathrm{RS}})e^{-u_a}.
\]
Markov's inequality therefore gives
\[
  \Prb_L\left[ \bar\rho(L+\widetilde J_a)>e^{-cu_a} \right] \le (1+c_{\mathrm{RS}})e^{-(1-c)u_a}.
\]
Condition \eqref{eq:fc-varying-union-condition} and a union bound produce one $L$ for which none of these events occurs.
Taking $\Lambda_h\coloneqq D_h^{1/m}L$ proves \eqref{eq:fc-varying-hole}.
Under \eqref{eq:fc-varying-threshold}, the left-hand side of \eqref{eq:fc-varying-union-condition}, with $c=1/2$, is at most $1/2$.
\end{proof}

The next result chooses shifts that make all relevant labeled Boolean intersections small at the same time.
It will be applied to the horizontal quotient torus.

Let $\T$ be a compact abelian group with Haar probability measure $\mu$.
Let $U_1,\ldots,U_N\subseteq\T$ be measurable, and suppose $\mu(U_a)\le\delta_a\le1$, where $\delta_a$ may depend on $a$.

For $I\subseteq[d]$, a \emph{labeled pattern} consists of a nonempty set $P\subseteq\{0,1\}^I$ and a \emph{label map} $\ell\colon P\to[N]$.
Given a labeled pattern and $z_1,\ldots,z_d\in\T$, define
\[
  \mathcal U_{I,P,\ell}(z) \coloneqq \bigcap_{\omega\in P} \left( U_{\ell(\omega)} -\sum_{i\in I}\omega_i z_i \right).
\]
For $j\in I$ and $e\in\{0,1\}$, let $\iota_{j,e}\colon\{0,1\}^{I\setminus\{j\}}\to\{0,1\}^I$ insert the coordinate $e$ in position $j$, and define the \emph{$e$-slice}
\[
  P_{j,e} \coloneqq \bigl\{ \omega'\in\{0,1\}^{I\setminus\{j\}}\colon \iota_{j,e}(\omega')\in P \bigr\}.
\]
When $P_{j,e}\ne\varnothing$, its \emph{induced label map} is
\[
  \ell_{j,e}(\omega') \coloneqq \ell\bigl(\iota_{j,e}(\omega')\bigr), \quad\text{for }\omega'\in P_{j,e}.
\]
Thus a slice retains the original labels of its surviving Boolean vertices.
A family $\mathscr C$ of triples $(I,P,\ell)$ is \emph{slice-closed} if $(I\setminus\{j\},P_{j,e},\ell_{j,e})\in\mathscr C$ whenever $(I,P,\ell)\in\mathscr C$, $j\in I$, and $P_{j,e}\ne\varnothing$.
This closure condition is exactly what the coordinate-by-coordinate induction requires: after one Boolean coordinate is fixed, every nonempty slice remains among the patterns controlled at the preceding stage.

\begin{theorem}\label{thm:fc-weighted-amplification}
Let $\mathscr C$ be a finite slice-closed family of triples $(I,P,\ell)$, allowing both the coordinate set $I$ and the pattern $P$ to vary.
If $L\ge2|\mathscr C|$, then there are shifts $z_1,\ldots,z_d\in\T$ such that, simultaneously for every $(I,P,\ell)\in\mathscr C$,
\begin{equation}\label{eq:fc-weighted-amplification-bound}
  \mu\bigl(\mathcal U_{I,P,\ell}(z)\bigr) \le L^{|P|-1} \prod_{\omega\in P}\delta_{\ell(\omega)}.
\end{equation}
\end{theorem}

\begin{proof}
If one of the factors on the right of \eqref{eq:fc-weighted-amplification-bound} is zero, then the intersection on the left is contained in a translate of a null set, so the desired inequality holds for every choice of shifts.
We may therefore consider only triples for which all relevant factors are positive.

For $0\le q\le d$, let $\mathscr C_{\le q}\coloneqq\{(I,P,\ell)\in\mathscr C\colon I\subseteq[q]\}$.
We choose the shifts inductively so that \eqref{eq:fc-weighted-amplification-bound} holds for every triple in $\mathscr C_{\le q}$.
For $q=0$, every nonempty pattern is a singleton, and the assertion is exactly $\mu(U_a)\le\delta_a$.

Suppose that $z_1,\ldots,z_{q-1}$ have been chosen for $\mathscr C_{\le q-1}$.
This property will not be affected by the choice of $z_q$, because none of the earlier intersections uses the $q$-th coordinate.

Consider a triple $(I,P,\ell)\in\mathscr C_{\le q}$ with $q\in I$, and put $I'\coloneqq I\setminus\{q\}$.
By slice closure, each nonempty $q$-slice belongs to $\mathscr C_{\le q-1}$.
If only the $e$-slice is nonempty, then
\[
  \mathcal U_{I,P,\ell}(z) = \mathcal U_{I',P_{q,e},\ell_{q,e}}(z)-e z_q.
\]
The two sets have the same measure, $|P_{q,e}|=|P|$, and the induction hypothesis gives the required bound without any loss.

It remains to treat triples for which both slices are nonempty.
For such a triple, set $A\coloneqq\mathcal U_{I',P_{q,0},\ell_{q,0}}(z)$ and $B\coloneqq\mathcal U_{I',P_{q,1},\ell_{q,1}}(z)$.
Then
\[
  \mathcal U_{I,P,\ell}(z) = A\cap(B-z_q).
\]
For Haar-uniform $z_q$, Fubini's theorem and translation invariance give
\[
  \E_{z_q}\mu\bigl(A\cap(B-z_q)\bigr) = \mu(A)\mu(B).
\]
The induction hypothesis bounds the last product by
\[
  L^{|P_{q,0}|-1}L^{|P_{q,1}|-1} \prod_{\omega\in P}\delta_{\ell(\omega)} = L^{|P|-2} \prod_{\omega\in P}\delta_{\ell(\omega)}.
\]
Markov's inequality shows that the target bound fails for this triple with probability at most $1/L$.
At most $|\mathscr C|\le L/2$ triples in $\mathscr C_{\le q}$ have two nonempty $q$-slices, so the union of their failure events has probability at most $1/2$.
We may therefore choose one $z_q$ that works for all of them.
This completes the induction.
\end{proof}

If $\mathscr A$ is a family of labeled patterns in $\{0,1\}^d$, regard each $(P,\ell)\in\mathscr A$ as the triple $([d],P,\ell)$.
Its \emph{slice closure}, denoted by $\operatorname{scl}(\mathscr A)$, is the smallest slice-closed family containing these triples.
The crude family containing every labeled nonempty subpattern on every surviving coordinate set has size at most
\[
  \sum_{I\subseteq[d]}(N+1)^{2^{|I|}} \le 2^d(N+1)^{2^d}.
\]
Consequently, with $\mathscr C\coloneqq\operatorname{scl}(\mathscr A)$ and $L\coloneqq2|\mathscr C|$,
\[
  \log L \le \log 2+d\log 2+2^d\log(N+1).
\]

\subsection{Finite-index lifting and entropy pruning}

We now build the simultaneous shifts into a single lattice.
The finite-index hypothesis is the algebraic bridge: after choosing a basis map $B$ for $\Lambda_v$, the inclusion $\Psi(\Z^d)\subseteq\Lambda_v$ yields a nonsingular integer matrix $M$ such that $\Psi=B\circ M$.
Since $M^{\mathsf T}$ is onto on the horizontal torus, the prescribed shifts can be lifted to shear parameters without changing the product of the horizontal and vertical determinants.
The next theorem formalizes this construction.

\begin{theorem}\label{thm:fc-weighted-lifting}
Let $\R^n=E\oplus F$ be an orthogonal decomposition with $\dim E=m$ and $\dim F=d$.
Let $K\subseteq E\oplus F$ be compact, and define its \emph{horizontal sections} by $K_y\coloneqq\{x\in E\colon x+y\in K\}$.
Suppose that:
\begin{enumerate}[label=(\roman*)]
\item $J_1,\ldots,J_N\subseteq E$ are compact sets, and $\Lambda_h\subseteq E$ is a lattice. For each $a\in[N]$, let $\delta_a\in[0,1]$ satisfy
  \[
    \bar\rho(\Lambda_h+J_a)\le\delta_a
    \quad\text{for every }a\in[N];
  \]
\item $\Lambda_v\subseteq F$ is a lattice, $\Psi\colon\R^d\to F$ is a linear isomorphism, and $\Psi(\Z^d)\subseteq\Lambda_v$ is a finite-index sublattice;
\item $\mathscr A$ is a finite family of pairs $(P,\ell)$ with $\varnothing\ne P\subseteq\{0,1\}^d$ and $\ell\colon P\to[N]$; the sets $P$ may have different cardinalities. For every $w\in F$, there are $\lambda\in\Lambda_v$ and $(P,\ell)\in\mathscr A$ satisfying
  \begin{equation}\label{eq:fc-weighted-witness}
    J_{\ell(\omega)}
    \subseteq
    K_{w-\lambda-\Psi(\omega)}
    \qquad\text{for every }\omega\in P.
  \end{equation}
\end{enumerate}
Let $\mathscr C\coloneqq\operatorname{scl}(\mathscr A)$ and $L\coloneqq2|\mathscr C|$.
Then there exists a lattice $\Lambda\subseteq E\oplus F$ such that
\begin{equation}\label{eq:fc-weighted-lifting-determinant}
  \det\Lambda
  =
  \det\Lambda_h\det\Lambda_v
\end{equation}
and
\begin{equation}\label{eq:fc-weighted-lifting-hole}
  \bar\rho(\Lambda+K) \le \max_{(P,\ell)\in\mathscr A} L^{|P|-1} \prod_{\omega\in P}\delta_{\ell(\omega)}.
\end{equation}
In particular, if $c>0$, $u_a\ge0$, and $\delta_a\le e^{-cu_a}$ for every $a\in[N]$, then
\begin{equation}\label{eq:fc-weighted-lifting-exponent}
  \log\bar\rho(\Lambda+K) \le \max_{(P,\ell)\in\mathscr A} \left( (|P|-1)\log L -c\sum_{\omega\in P}u_{\ell(\omega)} \right).
\end{equation}
Here and below, $\log 0\coloneqq-\infty$.
\end{theorem}

\begin{proof}
Choose a linear isomorphism $S\colon E\to\R^m$ satisfying $S(\Lambda_h)=\Z^m$, apply $S\oplus\operatorname{id}_F$ to $K$ and all the $J_a$, and keep the same notation for their images.
The section inclusions and hole proportions are unchanged.
At the end we will apply the inverse map, which multiplies the covolume by $|\det(S^{-1})|=\det\Lambda_h$.

Choose a basis map $B\colon\R^d\to F$ satisfying $B(\Z^d)=\Lambda_v$.
Since $\Psi(\Z^d)\subseteq\Lambda_v$ has finite index, there is a nonsingular integer matrix $M\in M_d(\Z)$ such that $\Psi=B\circ M$.
Let $\T^m\coloneqq\R^m/\Z^m$, let $\pi\colon\R^m\to\T^m$ be the quotient map, let $\mu$ be Haar probability measure on $\T^m$, and set $U_a\coloneqq\T^m\setminus\pi(J_a)$.
Since $J_a$ is compact, $U_a$ is Borel, and $\mu(U_a)=\bar\rho(\Z^m+J_a)\le\delta_a$.

Apply \cref{thm:fc-weighted-amplification} to $\mathscr C$, obtaining shifts $z_1,\ldots,z_d\in\T^m$.
Write $z\coloneqq(z_1,\ldots,z_d)$.

Consider the homomorphism from $(\T^m)^d$ to itself given by
\[
  \chi\longmapsto M^{\mathsf T}\chi.
\]
This map is onto.
It need not be injective, but Smith normal form reduces surjectivity to the fact that multiplication by a nonzero integer is onto on $\T^m$.
Choose $\chi\coloneqq(\chi_1,\ldots,\chi_d)$ satisfying $M^{\mathsf T}\chi=z$, and choose representatives $\widehat\chi_j\in\R^m$.
For $\nu\in\Z^d$, write $\nu\cdot\widehat\chi\coloneqq\sum_{j\in[d]}\nu_j\widehat\chi_j$.

Define the \emph{shear lattice} in the normalized coordinates by
\[
  \Gamma \coloneqq \bigl\{ (k-\nu\cdot\widehat\chi,B(\nu))\colon k\in\Z^m,\ \nu\in\Z^d \bigr\}.
\]
Its generating matrix is block triangular, with diagonal blocks $I_m$ and $B$, so $\det\Gamma=|\det B|=\det\Lambda_v$.
Let $\mathcal D_h\coloneqq[0,1)^m$, and let $\mathcal D_v$ be a half-open fundamental domain for $\Lambda_v$.
To reduce a point $(x,w)$ modulo $\Gamma$, first choose $\nu\in\Z^d$, uniquely outside a null boundary set, such that $w-B(\nu)\in\mathcal D_v$, and then choose $k\in\Z^m$, again uniquely outside a null boundary set, such that $x+\nu\cdot\widehat\chi-k\in\mathcal D_h$.
Thus $\mathcal D_h\times\mathcal D_v$ is a fundamental domain for $\Gamma$, up to its null boundary.

Let
\[
  \mathcal H
  \coloneqq
  (\mathcal D_h\times\mathcal D_v)\setminus(\Gamma+K),
  \quad\text{and}\quad
  \mathcal H_w
  \coloneqq
  \{x\in\mathcal D_h\colon(x,w)\in\mathcal H\}.
\]
The set $\Gamma+K$ is closed because $K$ is compact and $\Gamma$ is discrete.
Hence $\mathcal H$ is Borel and its sections $\mathcal H_w$ are measurable.

Fix $w\in\mathcal D_v$, and choose $\lambda=B(\nu_0)$ and $(P,\ell)\in\mathscr A$ from \eqref{eq:fc-weighted-witness}.
For $\omega\in P$, put $\nu_\omega\coloneqq\nu_0+M\omega$.
The points of $\Gamma$ with vertical coordinate $\lambda+\Psi(\omega)=B(\nu_\omega)$ are
\[
  \left( k-\nu_\omega\cdot\widehat\chi,\, B(\nu_\omega) \right), \quad\text{for }k\in\Z^m.
\]
Suppose $x\in\mathcal H_w$.
If $\pi(x+\nu_\omega\cdot\widehat\chi)$ belonged to $\pi(J_{\ell(\omega)})$, there would be a $k\in\Z^m$ such that $x-k+\nu_\omega\cdot\widehat\chi\in J_{\ell(\omega)}$.
By \eqref{eq:fc-weighted-witness}, the difference between $(x,w)$ and the corresponding point of $\Gamma$ would then belong to $K$, contradicting $x\in\mathcal H_w$.
Therefore $\pi(x+\nu_\omega\cdot\widehat\chi)\in U_{\ell(\omega)}$ for every $\omega\in P$.

In $\T^m$, we have $(M\omega)\cdot\chi=\omega\cdot(M^{\mathsf T}\chi)=\omega\cdot z$.
It follows that the image of $\mathcal H_w$ under the measure-preserving map $x\mapsto\pi(x+\nu_0\cdot\widehat\chi)$ is contained in
\[
  \bigcap_{\omega\in P} \bigl(U_{\ell(\omega)}-\omega\cdot z\bigr).
\]
The restriction of $\pi$ to $\mathcal D_h$ preserves measure up to the null boundary.
Hence \cref{thm:fc-weighted-amplification} gives
\[
  \vol_m(\mathcal H_w) \le L^{|P|-1} \prod_{\omega\in P}\delta_{\ell(\omega)}.
\]
The bound is uniform in $w$, even though the witness need not be chosen measurably.
Since $\vol_d(\mathcal D_v)=\det\Lambda_v$, Tonelli's theorem now gives
\[
  \bar\rho(\Gamma+K)
  =
  \frac{1}{\det\Lambda_v}
  \int_{\mathcal D_v}\vol_m(\mathcal H_w)\,\dd w
  \le
  \max_{(P,\ell)\in\mathscr A}
  L^{|P|-1}
  \prod_{\omega\in P}\delta_{\ell(\omega)}.
\]
Finally, define $\Lambda\coloneqq(S^{-1}\oplus\operatorname{id}_F)(\Gamma)$.
The inverse map preserves the uncovered proportion and gives $\det\Lambda=\det\Lambda_h\det\Lambda_v$.
This proves \eqref{eq:fc-weighted-lifting-determinant} and \eqref{eq:fc-weighted-lifting-hole}; taking logarithms proves \eqref{eq:fc-weighted-lifting-exponent}.
\end{proof}

We finish with the form of entropy pruning used in the main proof.
It separates the vertical entropy estimates from the horizontal geometry.
A label map need not be injective: all sums and products indexed by Boolean vertices count repeated labels with their multiplicities.

\begin{lemma}\label{lem:fc-entropy-pruning}
Let $m,d,N$ be positive integers, let $\rho>0$, and let $u_1,\ldots,u_N>0$.
For $\varepsilon\in\{0,1\}^d$, let $a_\varepsilon>0$, and put $A\coloneqq\sum_\varepsilon a_\varepsilon$ and $q_\varepsilon\coloneqq a_\varepsilon/A$.
Suppose that, for some $\zeta\in[0,1/2]$, a label map $\ell\colon\{0,1\}^d\to[N]$ satisfies
\begin{equation}\label{eq:fc-pruning-load-comparison}
  (1-\zeta)\rho a_\varepsilon \le u_{\ell(\varepsilon)} \le (1+\zeta)\rho a_\varepsilon.
\end{equation}
For $T_0>0$, define the \emph{retained pattern} $P\coloneqq\{\varepsilon\colon\rho a_\varepsilon\ge2T_0\}$.
Then every retained label has load at least $T_0$.
If in addition
\begin{equation}\label{eq:fc-pruning-upper-window}
  (1+\zeta)\rho A\max_\varepsilon q_\varepsilon \le\eta_m,
\end{equation}
then every retained label has load at most $\eta_m$, and
\begin{equation}\label{eq:fc-pruning-retained-load}
  \sum_{\varepsilon\in P}u_{\ell(\varepsilon)} \ge (1-\zeta)\rho A \left( 1- \frac{H(q)} {\log(\rho A/(2T_0))} \right),
\end{equation}
provided $\rho A>2T_0$.
\end{lemma}

\begin{proof}
The lower and upper load bounds follow immediately from \eqref{eq:fc-pruning-load-comparison} and \eqref{eq:fc-pruning-upper-window}.
If $\varepsilon\notin P$, then $q_\varepsilon<2T_0/(\rho A)$.
For every probability vector $q$ and every $\theta\in(0,1)$,
\begin{equation}\label{eq:fc-entropy-small-atom-tail}
  \sum_{q_\varepsilon\le\theta}q_\varepsilon \le \frac{H(q)}{\log(1/\theta)},
\end{equation}
because each summand in the indicated set satisfies $-q_\varepsilon\log q_\varepsilon \ge q_\varepsilon\log(1/\theta)$.
Apply \eqref{eq:fc-entropy-small-atom-tail} with $\theta\coloneqq2T_0/(\rho A)$, and use the lower bound in \eqref{eq:fc-pruning-load-comparison} on the retained vertices.
This proves \eqref{eq:fc-pruning-retained-load}.
\end{proof}

\section{Proof of the universal upper bound}\label{sec:fc-upper-proof}

We prove the upper bound in four steps.
First, we normalize $K$, choose the vertical dimension, and discretize all relevant horizontal sections.
Second, the weighted full-cube construction supplies a vertical lattice and a Boolean family of Gaussian loads above every vertical translate.
Third, entropy pruning removes the small loads, after which the Rogers--Schmidt estimate gives one horizontal lattice for every surviving section label.
Finally, weighted lifting combines the horizontal and vertical lattices, and Rogers' completion lemma turns the resulting almost covering into an exact covering.

Put $\ell_n\coloneqq\log\log n$.
For a sufficiently large absolute constant $A$, choose
\begin{equation}\label{eq:fc-d-choice}
  d \coloneqq \left\lceil A\ell_n^{4/3}(\log\ell_n)^{2/3} \right\rceil.
\end{equation}
Set $m\coloneqq n-d$.

For orientation, the Gaussian width is $\tau=(2\log d)^{-1}$, the finite-field grid has size $Q\asymp d(\log d)^2$, and the vertical dimension satisfies $d\asymp\ell_n^{4/3}(\log\ell_n)^{2/3}$.
Later, the number of section labels satisfies $\log N_0=O(d\log n)$, the pruning threshold satisfies $T_0=O(d\log n)$, and the target density is $\rho=n\log n\,d^{5/2+o(1)}$.
The first two parameters define the vertical cube and lattice; the choice of $d$ provides the required min-entropy; $N_0$ and $T_0$ control the common-lattice union bound; and $\rho$ provides the total horizontal load.
Unless explicitly stated otherwise, every $o(1)$ term in this section is deterministic and uniform in the convex body $K$, the vertical translate $w$, and the Boolean vertex $\varepsilon$.

The next elementary estimates verify every range condition used below and show that the vertical factor $d^{5/2+o(1)}$ becomes the stated power of $\log\log n$.

\begin{lemma}\label{lem:fc-d-scale}
If $A$ is sufficiently large, then
\begin{equation}\label{eq:fc-entropy-margin}
  ds_\infty\ge4\ell_n,
\end{equation}
and
\[
  d=o(n),\qquad d\log n=o(n),\qquad 4^d d\log n=o(n\log n),\qquad d\le n^{c_1},\quad\text{and}\quad 2\sqrt2\,d^{3/2}\sqrt{\log d}=o(n^{c_4}).
\]
Moreover,
\begin{equation}\label{eq:fc-d-power}
  d^{5/2+o(1)} = (\log\log n)^{10/3+o(1)}.
\end{equation}
\end{lemma}

\begin{proof}
The one-coordinate min-entropy estimate from \cref{lem:fc-entropy-asymptotics} gives
\[
  ds_\infty = \left(\frac{2\pi-4}{\sqrt\pi}+o(1)\right) \frac{d^{3/4}}{\sqrt{\log d}}.
\]
The choice \eqref{eq:fc-d-choice} gives
\[
  \frac{d^{3/4}}{\sqrt{\log d}} = \left( \sqrt{\frac34}\,A^{3/4}+o(1) \right)\ell_n.
\]
Thus the last expression is $(c_A+o(1))\ell_n$, where $c_A\coloneqq((2\pi-4)/\sqrt\pi)\sqrt{3/4}\,A^{3/4}$.
Since $c_A\to\infty$ with $A$, taking $A$ sufficiently large proves \eqref{eq:fc-entropy-margin}.

The remaining growth conditions follow because $d$ is a fixed power of $\log\log n$, while $\log n=\exp(\ell_n)$ and $n=\exp(\exp(\ell_n))$.
In particular, $4^dd=\exp(O(d))=n^{o(1)}=o(n)$, which gives $4^dd\log n=o(n\log n)$.

Finally, raising \eqref{eq:fc-d-choice} to the power $5/2$ proves \eqref{eq:fc-d-power}, with all lower-order factors absorbed into the $o(1)$.
\end{proof}

\begin{proof}[Proof of \cref{thm:upper}]
Choose once and for all an absolute constant $A$ large enough for \cref{lem:fc-d-scale}, and then choose an absolute constant $B_0$ with $B_0\ge8$.
From this point on $A$ and $B_0$ are fixed in that order.
We take $n\ge n_0(A,B_0)$, increasing $n_0(A,B_0)$ whenever one of the uniform estimates below requires it.

We first choose a common affine normalization and vertical subspace.

Fix an arbitrary $n$-dimensional convex body.
By affine invariance, apply \cref{prop:fc-variable-section-discretization} with $R\coloneqq\sqrt{d/\tau}=\sqrt{2d\log d}$, and henceforth write $K$ for the affine image given there.
Fix the resulting orthogonal decomposition $\R^n=E\oplus F$, identify $F$ isometrically with $\R^d$, and retain the discretization bodies $J_1,\ldots,J_N$.
Indeed,
\[
  2dR=2\sqrt 2\,d^{3/2}\sqrt{\log d}<n^{c_4}
\]
for all sufficiently large $n$ by \cref{lem:fc-d-scale}, so the range condition of the proposition holds.
Thus the vertical construction below is made in the same subspace $F$ for which all the pointwise section comparisons and labels are available.

We now choose the target density and the corresponding vertical and horizontal determinants.

Apply the weighted vertical construction of \cref{prop:fc-weighted-vertical}.
Let $\Lambda_v\subseteq F$ be its scaled vertical lattice, set $\Psi\coloneqq\Psi_d$ and $D_v\coloneqq\det\Lambda_v$, and note that $\Psi(\Z^d)\subseteq\Lambda_v$ is a finite-index sublattice.
For a fixed $w\in F$, choose any $\lambda_w\in\Lambda_v$ supplied by the construction and form the vertical points
\[
  y_{w,\varepsilon} \coloneqq w-\lambda_w-\Psi(\varepsilon), \quad\text{for }\varepsilon\in\{0,1\}^d.
\]
All estimates below are uniform in $w$ and hold for every admissible choice of $\lambda_w$; no choice as $w$ varies will be required.

Using the notation of \cref{prop:fc-weighted-vertical}, set $\mathsf w_{w,\varepsilon}\coloneqq D_v\gamma_d(y_{w,\varepsilon})$, $\mathsf W_w\coloneqq\sum_\varepsilon\mathsf w_{w,\varepsilon}$, and $q_{w,\varepsilon}\coloneqq\mathsf w_{w,\varepsilon}/\mathsf W_w$, and write $q_w\coloneqq(q_{w,\varepsilon})_{\varepsilon\in\{0,1\}^d}$.
Let $\xi_d=C_\xi\log\log d/\log d$ be the deterministic sequence in \cref{prop:fc-weighted-vertical}.
Since $ds\asymp ds_\infty\asymp d^{3/4}/\sqrt{\log d}$, for $s_\star\in\{s,s_\infty\}$ we have
\[
  \frac{\sqrt{ds_\star\log d}+o(\log d)}{ds_\star}
  =
  O\left(\sqrt{\frac{\log d}{ds_\star}}\right)
  +o\left(\frac{\log d}{ds_\star}\right)
  =o\left(\frac{\log\log d}{\log d}\right).
\]
After increasing the absolute constant $C_\xi$ if necessary, the same sequence $\xi_d$ therefore absorbs the error terms in both \eqref{eq:fc-max-atom} and \eqref{eq:fc-weight-entropy}.
Thus, uniformly in $w$,
\begin{align}
  d^{-5/2-\xi_d} &\le \mathsf W_w\le d^{-5/2+\xi_d},\notag\\
  -\log\max_\varepsilon q_{w,\varepsilon} &\ge(1-\xi_d)ds_\infty,\notag\\
  H(q_w)&\le(1+\xi_d)ds. \label{eq:fc-vertical-estimates}
\end{align}

Set the \emph{target density} to be
\begin{equation}\label{eq:fc-target-density}
  \rho \coloneqq B_0n\log n\,d^{5/2+\xi_d}.
\end{equation}
Choose the \emph{horizontal determinant} $D_h\coloneqq\vol_n(K)/(\rho D_v)$.
The \emph{ideal horizontal load} of the point $y_{w,\varepsilon}$ is $\widehat u_{w,\varepsilon}\coloneqq\rho\mathsf w_{w,\varepsilon}$.
By \eqref{eq:fc-vertical-estimates},
\begin{equation}\label{eq:fc-ideal-total-load}
  \sum_\varepsilon\widehat u_{w,\varepsilon} =\rho\mathsf W_w \ge B_0n\log n.
\end{equation}
On the other hand, \eqref{eq:fc-entropy-margin} and \eqref{eq:fc-vertical-estimates}, together with the initial choice of $A$, imply
\begin{equation}\label{eq:fc-ideal-max-load}
  \max_\varepsilon\widehat u_{w,\varepsilon}
  \le
  B_0n\log n\,d^{2\xi_d}
  e^{-(1-\xi_d)ds_\infty}
  \le
  B_0n\log n\,e^{-(1-2\xi_d)ds_\infty}
  \le
  \frac14\eta_m
\end{equation}
for all $n\ge n_0(A,B_0)$; in the second inequality we used $2\log d\le ds_\infty$, which holds for large $d$.
Indeed, \eqref{eq:fc-entropy-margin} and $\xi_d\to0$ make the preceding quantity at most
\[
  B_0n(\log n)^{-3+8\xi_d},
\]
which is $o(n)$, whereas $\eta_m=b_0(n-d)-3\log(n-d)\sim b_0n$.

We next prune the small loads and assign the corresponding section labels.

Fix an absolute constant $C_{\mathrm{sec}}$ at least as large as the constant in \cref{prop:fc-variable-section-discretization} and such that $C_{\mathrm{sec}}\ge2(1+c_{\mathrm{RS}})$.
Let $N_0\coloneqq(C_{\mathrm{sec}}m^2)^d$ be the deterministic upper bound in the section discretization, and put $T_0\coloneqq4\log(C_{\mathrm{sec}}N_0)$.
Since $T_0=O(d\log n)=o(n\log n)$, \eqref{eq:fc-ideal-total-load} implies $\rho\mathsf W_w>2T_0$, uniformly in $w$, for all sufficiently large $n$.

All the points $y_{w,\varepsilon}$ lie in $B_R^d$.
Put $u_a\coloneqq\vol_m(J_a)/D_h$.
For a fixed admissible vertical witness and every $\varepsilon\in\{0,1\}^d$, choose any label $\alpha_w(\varepsilon)\in[N]$ furnished by \cref{prop:fc-variable-section-discretization}, so that \eqref{eq:fc-variable-label-inclusion} and \eqref{eq:fc-variable-label-comparison} hold for this same label.
By \eqref{eq:fc-variable-label-gaussian} and the definition of $D_h$, the normalization is
\begin{align*}
  u_{\alpha_w(\varepsilon)} = \frac{\vol_m(J_{\alpha_w(\varepsilon)})}{D_h}
  &= \left(1+O\left(\frac1m+n^{-c_3}\right)\right) \frac{\vol_n(K)\gamma_d(y_{w,\varepsilon})} {\vol_n(K)/(\rho D_v)}\\
  &= \left(1+O\left(\frac1m+n^{-c_3}\right)\right) \rho D_v\gamma_d(y_{w,\varepsilon})\\
  &=
  \left(1+O\left(\frac1m+n^{-c_3}\right)\right)
  \rho\mathsf w_{w,\varepsilon}
  =
  \left(1+O\left(\frac1m+n^{-c_3}\right)\right)
  \widehat u_{w,\varepsilon}.
\end{align*}
In particular, there is a deterministic sequence $\zeta_n\to0$, independent of $K,w,\varepsilon$, such that $0\le\zeta_n\le1/2$ and
\[
  (1-\zeta_n)\widehat u_{w,\varepsilon} \le u_{\alpha_w(\varepsilon)} \le(1+\zeta_n)\widehat u_{w,\varepsilon}.
\]

For each $w$, define the retained pattern
\[
  \mathcal P_w \coloneqq \{\varepsilon\in\{0,1\}^d\colon \widehat u_{w,\varepsilon}\ge2T_0\}.
\]
Equations \eqref{eq:fc-ideal-total-load} and \eqref{eq:fc-ideal-max-load} imply
\[
  (1+\zeta_n)\rho\mathsf W_w \max_\varepsilon q_{w,\varepsilon} \le\eta_m
\]
for $n\ge n_0(A,B_0)$.
Apply \cref{lem:fc-entropy-pruning} for each $w$ with
\[
  a_\varepsilon=\mathsf w_{w,\varepsilon},
  \qquad
  A=\mathsf W_w,
  \qquad
  q=q_w,
  \qquad
  \ell=\alpha_w,
  \quad\text{and}\quad
  \zeta=\zeta_n.
\]
Every label occurring in $\mathcal P_w$ then lies in the common load window
\begin{equation}\label{eq:fc-active-load-window}
  T_0\le u_a\le\eta_m.
\end{equation}

By \cref{lem:fc-entropy-asymptotics} and \eqref{eq:fc-d-choice}, there is an absolute constant $c>0$ such that, uniformly in $w$,
\begin{align*}
  H(q_w) \le(1+\xi_d)ds=O(\ell_n)
  \quad\text{and}\quad
  \log\frac{\rho\mathsf W_w}{2T_0} \ge \log\frac{c n}{d} = \log n-O(\log d) = (1-o(1))\log n.
\end{align*}
Thus $H(q_w)/\log(\rho\mathsf W_w/(2T_0))=O(\ell_n/\log n)=o(1)$.
The retained-load conclusion of \cref{lem:fc-entropy-pruning} therefore yields
\begin{equation}\label{eq:fc-retained-actual-load}
  \sum_{\varepsilon\in\mathcal P_w} u_{\alpha_w(\varepsilon)} \ge (1-o(1))B_0n\log n
\end{equation}
uniformly in $w$.

It remains to choose one horizontal lattice and complete the covering.

Let $\mathscr A_0$ be the set of all labeled retained patterns that can arise from the preceding construction: choose $w\in F$, an admissible vertical witness $\lambda_w$, and labels satisfying the section inclusions, and retain the resulting pair $(\mathcal P_w,\alpha_w|_{\mathcal P_w})$.
This is an existential definition and makes no global choice as $w$ varies.
The family $\mathscr A_0$ is finite because $\mathcal P_w$ is a subset of the finite Boolean cube and its labels lie in $[N]$.
By \eqref{eq:fc-retained-actual-load}, every pattern in $\mathscr A_0$ is nonempty for sufficiently large $n$, and the retained-load bound holds uniformly over the family.

Call a label \emph{active} if it appears in a member of $\mathscr A_0$, and let $\mathcal I_{\mathrm{act}}$ be the set of active labels.
It is finite; set $N_{\mathrm{act}}\coloneqq|\mathcal I_{\mathrm{act}}|\le N_0$.
Discard the unused discretization bodies and reindex the bodies in $\mathcal I_{\mathrm{act}}$, together with every map in $\mathscr A_0$, by $[N_{\mathrm{act}}]$.
Denote the resulting family of labeled patterns by $\mathscr A$.
Thus \eqref{eq:fc-active-load-window} holds for every label in this reindexed family.

Apply \cref{lem:fc-varying-common-lattice} to this reindexed family, with determinant $D_h$.
We have
\[
  T_0=4\log(C_{\mathrm{sec}}N_0) \ge 2\log\bigl(2(1+c_{\mathrm{RS}})N_{\mathrm{act}}\bigr).
\]
Together with \eqref{eq:fc-active-load-window}, these are the hypotheses of that lemma.
We obtain one horizontal lattice $\Lambda_h$ such that every active label satisfies
\begin{equation}\label{eq:fc-active-hole}
  \bar\rho(\Lambda_h+J_a) \le e^{-u_a/2}.
\end{equation}

Put $\mathscr C\coloneqq\operatorname{scl}(\mathscr A)$ and $L\coloneqq2|\mathscr C|$.
For every $w$, the vertical construction and the section discretization supply at least one member $(P,\alpha)\in\mathscr A$ and a $\lambda\in\Lambda_v$ satisfying the witness condition \eqref{eq:fc-weighted-witness}.
Thus the lifting theorem requires no measurable selection.

Moreover, $\mathscr C$ is contained in the family of all labeled subpatterns of the Boolean cube, and hence
\[
  \log L \le C2^d\log(N_0+1).
\]
For every $(P,\alpha)\in\mathscr A$, one has $|P|-1\le2^d$.
Since $\log(N_0+1)=O(d\log n)$, \cref{lem:fc-d-scale} therefore gives
\[
  (|P|-1)\log L \le C4^d\log(N_0+1) = O(4^d d\log n) = o(n\log n).
\]
Apply \cref{thm:fc-weighted-lifting} with $\delta_a=e^{-u_a/2}$ and $c=1/2$.
Together with \eqref{eq:fc-active-hole} and \eqref{eq:fc-retained-actual-load}, it produces a lattice $\Lambda\subseteq\R^n$ with $\det\Lambda=D_hD_v$ and
\begin{align*}
  \log\bar\rho(\Lambda+K)
  \le
  C4^d\log(N_0+1)
  -\frac12
  \min_{(P,\alpha)\in\mathscr A}
  \sum_{\varepsilon\in P}
  u_{\alpha(\varepsilon)}
  \le \left( -\frac12(1-o(1))B_0+o(1) \right)n\log n.
\end{align*}
Because $B_0\ge8$ was fixed before $n_0(A,B_0)$, the last expression is smaller than $-2n\log n$ after one final increase of $n_0(A,B_0)$.
Since $e^{-2n\log n}=n^{-2n}<(n^n+1)^{-1}$, it follows that
\[
  \bar\rho(\Lambda+K)< (n^n+1)^{-1}.
\]
By \cref{lem:fc-rogers-completion}, $\Lambda+(1+1/n)K=\R^n$.
Equivalently, $\Lambda'\coloneqq(1+1/n)^{-1}\Lambda$ satisfies $\Lambda'+K=\R^n$.
Since $\det\Lambda'=(1+1/n)^{-n}\det\Lambda$, the density of this exact lattice covering is
\[
  \frac{\vol_n(K)}{\det\Lambda'} = \left(1+\frac1n\right)^n \frac{\vol_n(K)}{D_hD_v} \le e\rho.
\]

It remains to express the auxiliary dimension in the form used in the theorem.
Since $\log d=(4/3)\log\ell_n+(2/3)\log\log\ell_n+O(1)$, we have
\[
  \frac{(5/2+\xi_d)\log d}{\log\ell_n}
  =
  \frac{10}{3}
  +\frac43\xi_d
  +\left(\frac53+\frac23\xi_d\right)
   \frac{\log\log\ell_n}{\log\ell_n}
  +O\left(\frac1{\log\ell_n}\right).
\]
For all sufficiently large $n$, define
\[
  \varepsilon_n \coloneqq \max\left\{ 0,\, \frac{(5/2+\xi_d)\log d}{\log\ell_n} -\frac{10}{3} \right\}.
\]
Set $\varepsilon_n\coloneqq0$ for the remaining finitely many $n$.
Then $\varepsilon_n\to0$, independently of $K$, and $d^{5/2+\xi_d}\le\ell_n^{10/3+\varepsilon_n}$ for all sufficiently large $n$.
Consequently, \eqref{eq:fc-target-density} gives
\[
  \theta_L(K) \le eB_0n\log n\, (\log\log n)^{10/3+\varepsilon_n}.
\]
Finally, affine invariance permits us to undo the normalization without changing the covering density.
The constants and the sequence $\varepsilon_n$ are independent of $K$, so taking the supremum over all $n$-dimensional convex bodies proves \cref{thm:upper} with $C=eB_0$.
\end{proof}

\part{Worst-case lower bounds}

The lower-bound argument is technically independent of Part~I.
We therefore introduce its spherical, Poisson, and flat-torus notation locally.

\section{Tools for the lower bound}\label{lb:sec:prelim}

The proof uses three auxiliary tools.
Spherical-cap estimates control intersections among the sets of Poisson directions that cut prescribed residual vectors.
A packing argument on a flat torus produces many targets with separated residual directions.
Finally, a continuous Janson inequality converts bounds for the resulting Poisson witnesses and their dependency diagrams into a small covering probability for a fixed lattice.

We write $B_2^n$ for the Euclidean unit ball, $\Sph^{n-1}$ for its boundary, $\kappa_n\coloneqq\vol_n(B_2^n)$, and $\sigma$ for the rotation-invariant probability measure on $\Sph^{n-1}$.
Throughout Part~II, the symbols $c$ and $C$ denote positive absolute constants whose values may change from one occurrence to the next.
Every other constant is fixed within the statement or proof in which it is introduced, and the constants in the hierarchy of \cref{lb:sec:clusters} remain fixed thereafter.
Any dependence on fixed parameters is stated explicitly.

For $r>1$, let
\[
  p_n(r)\coloneqq\sigma\{u\in\Sph^{n-1}\colon |u_1|>1/r\}.
\]
For nonzero vectors $x$ and $z$ in $\R^n$, define their projective angle by
\[
  \angle_{\mathrm{proj}}(x,z)
  \coloneqq\min\{\angle(x,z),\pi-\angle(x,z)\}.
\]

\subsection{Spherical caps}

The following estimates will be used repeatedly.
The uniformity over a shell of width $O(\log n/n)$ is essential later.

\begin{lemma}\label{lb:lem:caps}
Fix $a>1$ and positive constants $A$ and $\beta$, and put $L\coloneqq a^2-1$.
The following assertions hold for all sufficiently large $n$, where the lower threshold for $n$ may depend on the fixed parameters displayed in the relevant assertion.

\begin{enumerate}[label=\textup{(\roman*)}]
\item There is a constant $C_{\mathrm{rad}}=C_{\mathrm{rad}}(a,A,\beta)>0$ such that, if $|t|\le A\log n+\beta$ and $r=ae^{t/n}$, then
\[
  C_{\mathrm{rad}}^{-1}e^{t/L}
  \le \frac{p_n(r)}{p_n(a)}
  \le C_{\mathrm{rad}}e^{t/L}.
\]

\item For every compact interval $I\subseteq(0,1)$ and every $H_0>0$, there is a constant $H_1=H_1(I,H_0)>0$ such that, for every $\tau\in I$ and $0\le h\le H_0\log n/n$,
\[
  \frac{\sigma\{u\colon u_1>\tau+h\}}
       {\sigma\{u\colon u_1>\tau\}}
  \ge \exp(-H_1nh-H_1).
\]

\item Let $r_0\coloneqq ae^{-A\log n/n}$ and $b\coloneqq ae^{\beta/n}$. For
$r,s\in[r_0-n^{-4},b+n^{-4}]$, define
\[
  C(x)\coloneqq\{u\in\Sph^{n-1}\colon |\ip{u}{x}|>1\}.
\]
For every fixed $J>0$, there is a constant $D=D(a,A,\beta,J)>0$ such that, whenever $|x|=r$, $|z|=s$, and $\angle_{\mathrm{proj}}(x,z)\ge D\sqrt{(\log n)/n}$, we have
\[
  \sigma(C(x)\cap C(z))
  \le n^{-J}\min\{\sigma(C(x)),\sigma(C(z))\}.
\]
The same conclusion holds when the threshold $1$ in both caps is replaced by $1+\eta$, uniformly for $0\le\eta\le n^{-4}$.
\end{enumerate}
\end{lemma}

\begin{proof}
The first coordinate of a uniform point of $\Sph^{n-1}$ has density
\[
  f_n(s)\coloneqq c_n(1-s^2)^{(n-3)/2},\quad\text{where}\quad
  c_n\coloneqq\frac{\Gamma(n/2)}{\sqrt\pi\,\Gamma((n-1)/2)}\asymp\sqrt n.
\]
If $\tau$ lies in a fixed compact subinterval of $(0,1)$, then
\begin{equation}\label{lb:eq:cap-endpoint}
  \int_\tau^1 f_n(s)\dd s
  \asymp n^{-1/2}(1-\tau^2)^{(n-1)/2}.
\end{equation}
Indeed,
\[
  \int_\tau^1 s(1-s^2)^{(n-3)/2}\dd s
  =\frac{(1-\tau^2)^{(n-1)/2}}{n-1}.
\]
On the interval of integration, the factor $s$ lies between $\tau$ and $1$.
Removing this factor changes the integral by a bounded factor, and $c_n\asymp\sqrt n$ gives \eqref{lb:eq:cap-endpoint}.

For part~\textup{(i)}, apply \eqref{lb:eq:cap-endpoint} with $\tau=1/r$.
When $r=ae^{t/n}$, Taylor expansion gives
\[
  \frac{n-1}{2}
  \bigl[\log(1-r^{-2})-\log(1-a^{-2})\bigr]
  =\frac{t}{a^2-1}+O\left(\frac{t^2+|t|}{n}\right).
\]
The error is $o(1)$ uniformly for $|t|\le A\log n+\beta$, and the remaining factors in \eqref{lb:eq:cap-endpoint} have bounded ratios.
This proves part~\textup{(i)}.

For part~\textup{(ii)}, put $p(\zeta)\coloneqq\int_\zeta^1f_n(s)\dd s$.
The endpoint estimate implies
\[
  -\frac{\mathrm d}{\mathrm d\zeta}\log p(\zeta)\le Cn
\]
uniformly on the fixed threshold interval.
Integrating from $\tau$ to $\tau+h$ yields
\[
  \log\frac{p(\tau+h)}{p(\tau)}\ge-Cnh.
\]
For sufficiently large $n$, the point $\tau+h$ remains in a slightly larger compact subinterval of $(0,1)$.
Part~\textup{(ii)} follows after adjusting the constant.

It remains to prove part~\textup{(iii)}.
Split each absolute cap into two oriented caps.
The shell assumptions place all oriented thresholds in a fixed compact interval $[\tau_0,\tau_1]\subsetneq(0,1)$ and give
\begin{equation}\label{lb:eq:threshold-difference}
  |\tau-\xi|\le C_2\frac{\log n}{n}.
\end{equation}
Let $e$ and $f$ be unit vectors with angle $\psi$, and define
\[
  q_*(\tau,\xi,\psi)
  \coloneqq\min\{|v|^2\colon \ip{v}{e}\ge\tau
                   \text{ and }\ip{v}{f}\ge\xi\}.
\]
We first show that there are constants $c_3,C_3>0$ such that
\begin{equation}\label{lb:eq:cap-angle-large}
  \psi^2\ge C_3\frac{\log n}{n}
\end{equation}
implies
\begin{equation}\label{lb:eq:minnorm-gain}
  q_*(\tau,\xi,\psi)
  \ge\max\{\tau^2,\xi^2\}+c_3\psi^2,
\end{equation}
unless the two constraints have empty intersection with the unit ball.

Fix a small constant $\psi_0>0$.
Suppose first that $0<\psi\le\psi_0$ and, by symmetry, that $\tau\ge\xi$.
In coordinates $e=(1,0)$ and $f=(\cos\psi,\sin\psi)$, we have
\[
  \xi-\tau\cos\psi
  =\tau(1-\cos\psi)-(\tau-\xi).
\]
By \eqref{lb:eq:threshold-difference} and \eqref{lb:eq:cap-angle-large}, a sufficiently large choice of $C_3$ makes the right side at least $c\psi^2$.
Both constraints are therefore active at the minimizer, and
\[
  q_*(\tau,\xi,\psi)
  =\tau^2+\frac{(\xi-\tau\cos\psi)^2}{\sin^2\psi}
  \ge\tau^2+c\psi^2.
\]
The case $\xi\ge\tau$ is symmetric.

Now suppose that $\psi\ge\psi_0$.
If the oriented caps intersect the unit ball, then
\[
  \tau+\xi\le |e+f|=2\cos(\psi/2),
\]
and hence $\psi\le\pi-\delta_0$ for some fixed $\delta_0>0$.
At equal thresholds, the minimum norm exceeds the square of the threshold by a positive amount uniformly for $\psi\in[\psi_0,\pi-\delta_0]$.
Continuity gives the same conclusion when $|\tau-\xi|$ is sufficiently small, and \eqref{lb:eq:threshold-difference} supplies this smallness for large $n$.
After decreasing $c_3$, we again obtain \eqref{lb:eq:minnorm-gain}.

Let $P$ denote orthogonal projection onto $\operatorname{span}\{e,f\}$.
If $u$ is uniform on $\Sph^{n-1}$, then
\[
  \Prb\left[|Pu|^2\ge q\right]=(1-q)^{(n-2)/2}\quad\text{for}\quad 0\le q\le1.
\]
Together with \eqref{lb:eq:cap-endpoint} and \eqref{lb:eq:minnorm-gain}, this gives
\begin{equation}\label{lb:eq:coarse-two-cap}
  \frac{\sigma\{u\colon \ip{u}{e}>\tau
                       \text{ and }\ip{u}{f}>\xi\}}
       {\min\{\sigma\{u\colon\ip{u}{e}>\tau\},
              \sigma\{u\colon\ip{u}{f}>\xi\}\}}
  \le n^{C_4}e^{-c_4n\psi^2}.
\end{equation}
Indeed, the numerator is at most $(1-q_*)^{(n-2)/2}$, while the smaller one-cap probability has threshold $\max\{\tau,\xi\}$ and is bounded below by \eqref{lb:eq:cap-endpoint}.

For two absolute caps, two sign choices have angle equal to the projective angle, and the other two have angle at least $\pi/2$.
Apply \eqref{lb:eq:coarse-two-cap} to all four pairs.
Taking $D$ sufficiently large gives the factor $n^{-J}$.
Replacing the threshold $1$ by $1+\eta$ changes each oriented threshold by $O(n^{-4})$, so the preceding estimates remain uniform for $0\le\eta\le n^{-4}$.
\end{proof}

\subsection{Unfolding on a flat torus}

For a lattice $\Lambda$, let $\T_\Lambda\coloneqq\R^n/\Lambda$, equipped with normalized Haar measure.
If $y\in\T_\Lambda$, an expression involving the coset $y-\Lambda$ is evaluated using any representative of $y$; the result is independent of that choice.
In particular, for every nonnegative measurable function $f$, the sum below is a well-defined function on $\T_\Lambda$, and we use the unfolding identity
\begin{equation}\label{lb:eq:unfold}
  \int_{\T_\Lambda}\sum_{\lambda\in\Lambda}f(y-\lambda)\dd y
  =\frac1{\det\Lambda}\int_{\R^n}f(x)\dd x.
\end{equation}
To obtain this identity, integrate over a fundamental domain and translate its lattice copies.

The next packing lemma does not require the quotient map to be injective on a small Euclidean ball.

\begin{lemma}\label{lb:lem:torus-packing}
Let $G\subseteq\T_\Lambda$ have Haar measure at least $\gamma>0$, and let $\delta>0$.
There are points $y_1,\ldots,y_T\in G$ such that $\dist_{\T_\Lambda}(y_\alpha-y_\beta,0)>\delta$ and $\dist_{\T_\Lambda}(y_\alpha+y_\beta,0)>\delta$ for all distinct $\alpha$ and $\beta$, and such that $T\ge \gamma\det\Lambda/(2\kappa_n\delta^n)$.
\end{lemma}

\begin{proof}
Choose the points greedily.
After selecting a point $y$, remove from $G$ the images of the two closed balls $y+\delta B_2^n$ and $-y+\delta B_2^n$.
The selected points satisfy both separation conditions.
Since every $\delta$-separated subset of the compact torus is finite, the process ends with a maximal list $y_1,\ldots,y_T$.

The image of a Euclidean ball of radius $\delta$ has normalized Haar measure at most $\kappa_n\delta^n/\det\Lambda$; identifications under the quotient map can only decrease its measure.
By maximality, the $2T$ removed images cover $G$.
Hence
\[
  \gamma\le \frac{2T\kappa_n\delta^n}{\det\Lambda},
\]
which proves the lemma.
\end{proof}

\subsection{A continuous Janson inequality}

Here a witness is an injective tuple of Poisson points chosen from prescribed measurable slots.
The diagram sums below record the ways in which two witnesses can share Poisson points.

Let $(\Omega,\mu)$ be a finite nonatomic standard Borel measure space, and let $\Pi$ be a Poisson process with intensity measure $\mu$.
Let $\mathcal A$ be a finite index set.
For each $\alpha\in\mathcal A$, let $I_\alpha$ be a finite nonempty set.
For every $i\in I_\alpha$, let $E_{\alpha i}\subseteq\Omega$ be measurable with $0<\mu(E_{\alpha i})<\infty$, and put
\[
  \lambda_{\alpha i}\coloneqq\mu(E_{\alpha i})
  \quad\text{and}\quad
  \nu_\alpha\coloneqq\prod_{i\in I_\alpha}\lambda_{\alpha i}.
\]
Let $Z_\alpha$ count the injective tuples $(u_i)_{i\in I_\alpha}$ of points of $\Pi$ such that $u_i\in E_{\alpha i}$ for every $i\in I_\alpha$.
The Campbell--Mecke formula~\cite{Mecke1967} gives
\[
  \E Z_\alpha
  =\int_{\Omega_{\ne}^{I_\alpha}}
    \prod_{i\in I_\alpha}\mathbf 1_{E_{\alpha i}}(u_i)\,
    \dd\mu^{\otimes I_\alpha}((u_i)_i)
  =\prod_{i\in I_\alpha}\lambda_{\alpha i}
  =\nu_\alpha,
\]
where $\Omega_{\ne}^{I_\alpha}$ is the set of tuples with pairwise distinct coordinates.
Its complement has product measure zero because $\mu$ is nonatomic.
Thus $\nu_\alpha$ is the mean number of witnesses of type $\alpha$.

A partial injection $M\colon I_\alpha\dashrightarrow I_\beta$ is an injective map from a subset of $I_\alpha$ into $I_\beta$.
For a nonempty partial injection $M$, define
\begin{equation}\label{lb:eq:diagram-weight}
  W_{\alpha\beta}(M)\coloneqq
  \nu_\alpha\nu_\beta
  \prod_{(i,j)\in M}
  \frac{\mu(E_{\alpha i}\cap E_{\beta j})}
       {\lambda_{\alpha i}\lambda_{\beta j}}.
\end{equation}
Let $D_{\alpha\beta}$ be the sum of \eqref{lb:eq:diagram-weight} over all nonempty partial injections from $I_\alpha$ to $I_\beta$.

We first record the finite indexed form of Janson's inequality that will be used in the approximation.

\begin{lemma}\label{lb:lem:finite-janson}
Let $(X_v)_{v\in V}$ be independent Bernoulli variables on a finite set $V$, and let $(H_h)_{h\in\mathcal H}$ be a finite indexed family of nonempty subsets of $V$.
Repetitions among the underlying subsets are allowed.
Put
\[
  I_h\coloneqq\prod_{v\in H_h}X_v,\qquad
  Y\coloneqq\sum_{h\in\mathcal H}I_h,\quad\text{and}\quad
  m\coloneqq\E Y,
\]
and define the ordered dependency sum, including the diagonal, by
\[
  \Delta^*\coloneqq\sum_{(h,h')\colon H_h\cap H_{h'}\ne\varnothing}
  \E(I_hI_{h'}).
\]
If $m>0$, then
\begin{equation}\label{lb:eq:finite-indexed-janson}
  \Prb\left[Y=0\right]\le\exp\left(-\frac{m^2}{2\Delta^*}\right).
\end{equation}
\end{lemma}

\begin{proof}
For distinct underlying subsets, \eqref{lb:eq:finite-indexed-janson} is the lower-tail form of Janson's inequality with $t=m$; see Janson, {\L}uczak, and Ruci\'nski~\cite[Theorem~2.14]{JLR}.
With our ordered convention, $\Delta^*$ is the quantity denoted there by $\bar\Delta$.

To remove the distinctness assumption, fix $0<\varepsilon<1$.
For each index $h$, introduce a private Bernoulli variable $B_h$ with parameter $1-\varepsilon$.
Assume that all private variables are independent of one another and of the original family, and replace $I_h$ by $I_h B_h$.
The augmented underlying sets are distinct.
Their mean and ordered dependency sum are
\[
  m_\varepsilon=(1-\varepsilon)m
\]
and
\[
  \Delta_\varepsilon^*
  =(1-\varepsilon)\sum_h\E I_h
   +(1-\varepsilon)^2
    \sum_{\substack{h\ne h'\\H_h\cap H_{h'}\ne\varnothing}}
    \E(I_hI_{h'}).
\]
The probability that the augmented family has no present member differs from $\Prb\left[Y=0\right]$ by at most $|\mathcal H|\varepsilon$.
Since $m_\varepsilon\to m$ and $\Delta_\varepsilon^*\to\Delta^*$, letting $\varepsilon\downarrow0$ proves the indexed form.
\end{proof}

\begin{lemma}\label{lb:lem:poisson-janson}
There is an absolute constant $c>0$ such that
\begin{equation}\label{lb:eq:poisson-janson}
  \Prb\left[Z_\alpha=0\text{ for every }\alpha\right]
  \le
  \exp\left(
  -c\frac{(\sum_\alpha\nu_\alpha)^2}
           {\sum_{\alpha,\beta}D_{\alpha\beta}}
  \right).
\end{equation}
In particular, $D_{\alpha\alpha}$ contains the full diagonal matching, whose weight is $\nu_\alpha$.
\end{lemma}

\begin{proof}
We approximate the Poisson process by independent cell-occupancy variables.
All sets may first be modified on null sets, so we may assume that they are Borel.

\medskip\noindent\emph{Nested partitions.}
The standard Borel structure has a countable family of Borel sets that separates points.
Start with the Boolean atoms generated by the finitely many sets $E_{\alpha i}$.
At stage $k$, refine the preceding partition by the first $k$ point-separating Borel sets.
Using nonatomicity, subdivide every positive-measure cell into finitely many measurable pieces of measure at most $2^{-k}$.
These choices can be made recursively, so they give nested finite partitions $(\mathcal P_k)_{k\ge1}$.
After discarding null cells, we have
\[
  \delta_k\coloneqq\max_{C\in\mathcal P_k}\mu(C)\longrightarrow0.
\]
Every cell is contained in or disjoint from each $E_{\alpha i}$, and every two distinct points of $\Omega$ are eventually separated by the partitions.

For $C\in\mathcal P_k$, let $X_C$ indicate that $\Pi$ meets $C$.
The variables $(X_C)_{C\in\mathcal P_k}$ are independent, and
\begin{equation}\label{lb:eq:cell-occupancy}
  q_C\coloneqq\Prb\left[X_C=1\right]=1-e^{-\mu(C)}
  =\mu(C)(1+O(\delta_k))
\end{equation}
uniformly in $C$.

An indexed hyperedge consists of a label $\alpha$ and an assignment $(C_i)_{i\in I_\alpha}$ of pairwise distinct cells satisfying $C_i\subseteq E_{\alpha i}$ for every $i$.
Its underlying vertex set is $\{C_i\colon i\in I_\alpha\}$, and its indicator is $\prod_i X_{C_i}$.
Different labels and different assignments remain distinct indices even when they have the same underlying cell set.
Let $Y_k$ count the present indexed hyperedges, let $m_k\coloneqq\E Y_k$, and let $\Delta_k^*$ be the ordered dependency sum from \cref{lb:lem:finite-janson}.
Then
\begin{equation}\label{lb:eq:finite-janson-step}
  \Prb\left[Y_k=0\right]\le
  \exp\left(-\frac{m_k^2}{2\Delta_k^*}\right)
\end{equation}
whenever $m_k>0$.

\medskip\noindent\emph{First moment.}
Fix $\alpha$ and put $s_\alpha\coloneqq|I_\alpha|$.
By \eqref{lb:eq:cell-occupancy},
\begin{equation}\label{lb:eq:first-moment-cells}
  \E Y_{k,\alpha}
  =\bigl(1+O(s_\alpha\delta_k)\bigr)
  \sum_{\substack{(C_i)_{i\in I_\alpha}\text{ pairwise distinct}\\
                   C_i\subseteq E_{\alpha i}}}
  \prod_{i\in I_\alpha}\mu(C_i).
\end{equation}
Without the distinctness restriction, the sum equals $\nu_\alpha$.
For fixed distinct $i$ and $j$, assignments with $C_i=C_j$ have total product mass at most
\[
  \sum_{C\subseteq E_{\alpha i}\cap E_{\alpha j}}\mu(C)^2
  \prod_{\ell\ne i,j}\lambda_{\alpha\ell}
  \le
  \delta_k\mu(E_{\alpha i}\cap E_{\alpha j})
  \prod_{\ell\ne i,j}\lambda_{\alpha\ell}.
\]
A union bound over the finitely many pairs of slots shows that
\begin{equation}\label{lb:eq:janson-first-limit}
  m_k\longrightarrow\sum_\alpha\nu_\alpha.
\end{equation}
In particular, $m_k>0$ for all sufficiently large $k$.

\medskip\noindent\emph{Dependency sum.}
Fix $\alpha$, $\beta$, and a nonempty partial injection $M\colon I_\alpha\dashrightarrow I_\beta$.
Form the quotient of the disjoint union $I_\alpha\sqcup I_\beta$ by identifying $i$ and $j$ whenever $(i,j)\in M$.
For each quotient class $Q$, let $F_Q$ be the intersection of the one or two slot sets indexed by $Q$.
If the quotient classes choose cells independently and unintended equalities are ignored, their total product mass is
\[
  \prod_Q\mu(F_Q)
  =\nu_\alpha\nu_\beta
   \prod_{(i,j)\in M}
   \frac{\mu(E_{\alpha i}\cap E_{\beta j})}
        {\lambda_{\alpha i}\lambda_{\beta j}}
  =W_{\alpha\beta}(M).
\]
If two distinct quotient classes $Q$ and $Q'$ use the same cell, then the total product mass is at most
\[
  \delta_k\mu(F_Q\cap F_{Q'})
  \prod_{Q''\ne Q,Q'}\mu(F_{Q''}).
\]
There are finitely many pairs of quotient classes, so requiring the quotient classes to use pairwise distinct cells changes the total by $O(\delta_k)$.

For two internally injective cell assignments, the complete set of cross equalities is uniquely a partial injection: each cell on the left equals at most one cell on the right, and conversely.
The two indexed hyperedges intersect exactly when this partial injection is nonempty.
Hence the ordered dependency sum decomposes without omission or repetition over the nonempty partial injections.
Replacing each factor $\mu(C)$ by $q_C$ changes each finite sum by another $O(\delta_k)$, by \eqref{lb:eq:cell-occupancy}.
Therefore
\begin{equation}\label{lb:eq:janson-second-limit}
  \Delta_k^*\longrightarrow
  \sum_{\alpha,\beta\in\mathcal A}D_{\alpha\beta}.
\end{equation}

\medskip\noindent\emph{Limiting event.}
Because the partitions are nested, the events $\{Y_k>0\}$ increase with $k$.
A witness in pairwise distinct coarse cells remains a witness after one occupied subcell is chosen inside each coarse cell.
Conversely, the Poisson process is finite and simple almost surely because $\mu(\Omega)<\infty$ and $\mu$ is nonatomic.
Every finite list of distinct Poisson points is eventually separated by the partitions.
Consequently,
\[
  \bigcup_k\{Y_k>0\}
  =\left\{\sum_\alpha Z_\alpha>0\right\}
  \quad\text{almost surely}.
\]
Thus
\[
  \Prb\left[Y_k=0\right]\downarrow
  \Prb\left[Z_\alpha=0\text{ for every }\alpha\right].
\]
Letting $k\to\infty$ in \eqref{lb:eq:finite-janson-step}, and using \eqref{lb:eq:janson-first-limit} and \eqref{lb:eq:janson-second-limit}, proves \eqref{lb:eq:poisson-janson} with $c=1/2$.

For $\alpha=\beta$, the full matching $i\mapsto i$ has weight $\nu_\alpha^2\prod_i\lambda_{\alpha i}^{-1}=\nu_\alpha$.
\end{proof}

For distinct $\alpha$ and $\beta$, put
\begin{equation}\label{lb:eq:rho-cross}
  \rho_{\alpha\beta}\coloneqq
  \sum_{\substack{i\in I_\alpha\\j\in I_\beta}}
  \frac{\mu(E_{\alpha i}\cap E_{\beta j})}
       {\lambda_{\alpha i}\lambda_{\beta j}}.
\end{equation}
Dropping the matching constraint and expanding the resulting product gives
\begin{equation}\label{lb:eq:cross-diagram}
  D_{\alpha\beta}
  \le\nu_\alpha\nu_\beta
  \bigl(e^{\rho_{\alpha\beta}}-1\bigr).
\end{equation}
This estimate is the interface used later to convert cross-target cap overlaps into the dependency term in the continuous Janson inequality.

\section{The random convex body}\label{lb:sec:body}

The construction has three components.
The cube fixes the John normalization, the outer Euclidean ball gives a uniform radius bound, and the random slabs create many potential holes while preserving a constant fraction of the volume with positive probability.

For the remainder of the lower-bound argument, fix
\[
  a\coloneqq2,\qquad L\coloneqq a^2-1=3,\quad\text{and}\quad \beta\coloneqq1,
\]
and put $b\coloneqq ae^{\beta/n}$.
Let $\Pi$ be a Poisson process on $\Sph^{n-1}$ with intensity measure
\[
  \mu\coloneqq M\sigma,\quad\text{where}\quad M\coloneqq\frac{L}{p_n(a)}.
\]
This normalization gives $Mp_n(a)=L$, so the cap associated with a residual vector of norm $a$ has constant intensity.
It balances the random cuts against survival of a constant fraction of the volume.

Define
\begin{equation}\label{lb:eq:random-body}
  K_\Pi\coloneqq bB_2^n\cap[-1,1]^n\cap
  \bigcap_{u\in\Pi}\{x\in\R^n\colon |\ip{u}{x}|\le1\}.
\end{equation}
The set $K_\Pi$ is an origin-symmetric convex body, and
\[
  B_2^n\subseteq K_\Pi\subseteq[-1,1]^n\cap bB_2^n.
\]

\begin{proposition}\label{lb:prop:volume}
There are absolute constants $c_v,c_p>0$ such that, for all sufficiently large $n$,
\[
  \Prb\left[\vol_n(K_\Pi)\ge c_v\kappa_n a^n\right]\ge c_p.
\]
For every realization of $\Pi$, the John ellipsoid of $K_\Pi$ is $B_2^n$.
\end{proposition}

\begin{proof}
For a vector $x$ with $|x|=r>1$, let
\[
  C_x\coloneqq\{u\colon |\ip{u}{x}|>1\}
  \quad\text{and}\quad
  c(x)\coloneqq\mu(C_x)=Mp_n(r).
\]
The Poisson void probability gives
\[
  \Prb\left[\Pi\cap C_x=\varnothing\right]=e^{-c(x)}.
\]
Restrict polar integration to the shell $r\in[ae^{-1/n},a]$.
By \cref{lb:lem:caps}\textup{(i)}, the quantity $c(x)$ is bounded above by an absolute constant throughout this shell.
If $v$ is uniform on $\Sph^{n-1}$, then, for every such $r$,
\[
  \Prb\left[rv\notin[-1,1]^n\right]
  \le n\Prb\left[|v_1|>1/a\right]
  \le ne^{-c n}=o(1),
\]
where the last inequality follows from \eqref{lb:eq:cap-endpoint}.
Fubini's theorem therefore gives
\begin{align*}
  \E\vol_n(K_\Pi)
  &\ge (1-o(1))e^{-C}
     \vol_n\bigl(aB_2^n\setminus ae^{-1/n}B_2^n\bigr)\\
  &\ge c\kappa_n a^n.
\end{align*}

We also have $\vol_n(K_\Pi)\le\kappa_n b^n=e^\beta\kappa_n a^n$.
Let $X\coloneqq\vol_n(K_\Pi)$ and $t\coloneqq(c/2)\kappa_n a^n$.
Since
\[
  X\le t+e^\beta\kappa_n a^n\mathbf 1_{\{X\ge t\}},
\]
the expectation bound implies
\[
  \Prb\left[X\ge\frac c2\kappa_n a^n\right]
  \ge \frac{c}{2e^\beta}.
\]
This proves the first assertion after renaming the constants.

It remains to identify the John ellipsoid.
The contact points $\pm e_i$ of the cube, each with weight $1/2$, satisfy
\[
  \sum_{i=1}^n\sum_{\varepsilon\in\{-1,1\}}
  \frac12\varepsilon e_i=0
  \quad\text{and}\quad
  \sum_{i=1}^n\sum_{\varepsilon\in\{-1,1\}}
  \frac12(\varepsilon e_i)(\varepsilon e_i)^{\mathsf T}=I_n.
\]
By John's contact criterion~\cite{John1948,Ball1992}, the unit ball is the maximum-volume ellipsoid contained in $[-1,1]^n$.
Every ellipsoid contained in $K_\Pi$ is contained in the cube, while $B_2^n\subseteq K_\Pi$.
Hence $B_2^n$ is the John ellipsoid of $K_\Pi$.
\end{proof}

For the stability argument, we also use an enlarged body.
For $0\le\eta\le n^{-4}$, let
\begin{equation}\label{lb:eq:enlarged-body}
  \widehat K_{\Pi,\eta}
  \coloneqq(b+\eta)B_2^n\cap
  \bigcap_{u\in\Pi}\{x\colon |\ip{u}{x}|\le1+\eta\}.
\end{equation}
The cube is omitted, and hence $K_\Pi\subseteq\widehat K_{\Pi,\eta}$.
All fixed-lattice estimates below are uniform for $0\le\eta\le n^{-4}$.

\section{Fixed lattices and separated targets}\label{lb:sec:fixed}

We first identify a positive-measure set of good torus targets: they have no short residual vector, and both the number of relevant residuals and their aggregate cap cost are $O(R)$.
We then pack this set to obtain exponentially many targets whose residual directions are separated across distinct targets.

Fix a lattice $\Lambda$, and put
\begin{equation}\label{lb:eq:R-def}
  R\coloneqq\frac{\kappa_n a^n}{\det\Lambda}
  \quad\text{and}\quad
  F_n\coloneqq n\log n.
\end{equation}
We work in the range
\begin{equation}\label{lb:eq:R-range}
  c_*\le R\le c_0F_n,
\end{equation}
where $c_*\coloneqq e^{-\beta}/4$ and the upper-range constant $c_0>0$ will be chosen only after the cross-target exponent has been fixed.
The bodies considered below have outer radius $b+\eta$, with $0\le\eta\le n^{-4}$.
If such a body covers with lattice $\Lambda$, then integration of the covering multiplicity gives
\[
  \det\Lambda\le\kappa_n(b+\eta)^n.
\]
Consequently,
\[
  R\ge\left(\frac{a}{b+\eta}\right)^n
  =e^{-\beta+o(1)}\ge\frac12e^{-\beta}>c_*.
\]
Thus the lower bound in \eqref{lb:eq:R-range} excludes no lattice that can cover one of these bodies.

Fix $A\coloneqq10$, and let
\[
  r_0\coloneqq a\exp(-A\log n/n).
\]
For $0\le\eta\le n^{-4}$, put $b_\eta\coloneqq b+\eta$ and define
\[
  C_x^\eta\coloneqq\{u\colon |\ip{u}{x}|>1+\eta\}
  \quad\text{and}\quad
  c_\eta(x)\coloneqq\mu(C_x^\eta).
\]
For $y\in\T_\Lambda$, let
\[
  X_y^\eta\coloneqq(y-\Lambda)\cap
  \bigl(b_\eta B_2^n\setminus(r_0-\eta)B_2^n\bigr).
\]
Set
\[
  m(y)\coloneqq|X_y^\eta|,\qquad
  q(t)\coloneqq-\log(1-e^{-t}),\quad\text{and}\quad
  Q(y)\coloneqq\sum_{x\in X_y^\eta}q(c_\eta(x)).
\]
The quantity $1-e^{-c_\eta(x)}$ is the probability that the cap $C_x^\eta$ is hit by the Poisson process.
Thus $Q(y)$ is the total logarithmic cost of hitting all residual caps if their dependencies are ignored.

\begin{lemma}\label{lb:lem:good-targets}
There are absolute constants $C$, $C_A$, and $\gamma>0$ such that, for all sufficiently large $n$, the following holds uniformly over every lattice satisfying \eqref{lb:eq:R-range} and every $0\le\eta\le n^{-4}$.
The set $G_\Lambda^\eta$ of targets $y$ that satisfy
\begin{equation}\label{lb:eq:good-properties}
  (y-\Lambda)\cap(r_0-\eta)B_2^n=\varnothing,
  \qquad m(y)\le CR,
  \quad\text{and}\quad Q(y)\le CR
\end{equation}
has Haar measure at least $\gamma$.
Moreover, if $y$ is good and $x\in X_y^\eta$, then
\begin{equation}\label{lb:eq:cap-intensity-range}
  n^{-C_A}\le c_\eta(x)\le C_A.
\end{equation}
\end{lemma}

\begin{proof}
The unfolding identity gives
\begin{align}
  \E_y m(y)
  &=\frac{\kappa_n(b_\eta^n-(r_0-\eta)^n)}{\det\Lambda}
    \le C_mR,\label{lb:eq:mean-m}\\
  \Prb_y\left[(y-\Lambda)\cap(r_0-\eta)B_2^n\ne\varnothing\right]
  &\le\frac{\kappa_n(r_0-\eta)^n}{\det\Lambda}
    \le 2Rn^{-A}=o(1).\label{lb:eq:inner-bad}
\end{align}
All constants in this proof are uniform in $\Lambda$ and $\eta$.

Write $r=ae^{t/n}$ in the relevant shell.
Then
\[
  c_\eta(r)=M p_n\!\left(\frac{r}{1+\eta}\right)
  \quad\text{and}\quad
  n\log\frac{r/(1+\eta)}a
  =t-n\log(1+\eta)=t+o(1).
\]
\Cref{lb:lem:caps}\textup{(i)} therefore gives
\begin{equation}\label{lb:eq:c-radial}
  C^{-1}e^{t/L}\le c_\eta(r)\le Ce^{t/L}
  \quad\text{for}\quad -A\log n-O(1)\le t\le\beta+O(1).
\end{equation}
This proves \eqref{lb:eq:cap-intensity-range} after enlarging $C_A=C_A(a,A,\beta)$.
Since $a$, $A$, and $\beta$ are fixed, $C_A$ is absolute.

For $0<s\le C_A$, the function $q(s)=-\log(1-e^{-s})$ satisfies $q(s)\le C(1+|\log s|)$.
Applying the unfolding identity again and using the substitution $r=ae^{t/n}$, we obtain
\begin{align*}
  \E_yQ(y)
  =\frac{n\kappa_n}{\det\Lambda}
    \int_{r_0-\eta}^{b_\eta}r^{n-1}q(c_\eta(r))\dd r
  \le C_QR\int_{-A\log n-O(1)}^{\beta+O(1)}
    e^t(1+|t|)\dd t
  \le C_Q'R.
\end{align*}
Choose $C\ge8\max\{C_m,C_Q'\}$.
Markov's inequality shows that each of the events $m(y)>CR$ and $Q(y)>CR$ has probability at most $1/8$.
For sufficiently large $n$, the probability in \eqref{lb:eq:inner-bad} is at most $1/4$.
Hence the three conditions in \eqref{lb:eq:good-properties} hold on a set of Haar measure at least $1/2$.
We may take $\gamma\coloneqq1/2$.
\end{proof}

If a good target $y$ satisfies $m(y)=0$, then every residual vector $y-\lambda$ has norm greater than $b_\eta$.
Therefore $y$ is a deterministic hole in $\widehat K_{\Pi,\eta}+\Lambda$.
In the witness construction below, we may consequently restrict attention to good targets with $m(y)\ge1$.

We now choose targets whose residual directions are separated across different targets.
Let
\[
  \alpha_{\mathrm{sep}}\coloneqq D_{\mathrm{sep}}\sqrt{\frac{\log n}{n}}
  \quad\text{and}\quad
  \delta\coloneqq3\bigl(2\eta+b-r_0+b_\eta\alpha_{\mathrm{sep}}\bigr),
\]
where $D_{\mathrm{sep}}$ is for now an arbitrary fixed constant.
It will be chosen in \cref{lb:sec:clusters}, after the desired cross-target exponent $J$ is fixed.
Suppose that $x\in X_y^\eta$ and $z\in X_{y'}^\eta$ have projective angle less than $\alpha_{\mathrm{sep}}$.
For some sign $\varepsilon\in\{-1,1\}$,
\[
  |x-\varepsilon z|
  \le 2\eta+b-r_0+b_\eta\alpha_{\mathrm{sep}}<\delta.
\]
The vector $x-\varepsilon z$ represents $y-\varepsilon y'$ on the torus.
The contrapositive of this observation, together with \cref{lb:lem:torus-packing}, gives the next lemma.

\begin{lemma}\label{lb:lem:separated-targets}
For every fixed $D_{\mathrm{sep}}>0$ and all sufficiently large $n$, there are good targets $y_1,\ldots,y_T$ such that residual directions from distinct targets have projective angle at least $\alpha_{\mathrm{sep}}$, and
\begin{equation}\label{lb:eq:number-targets}
  T\ge c\frac{\det\Lambda}{\kappa_n\delta^n}
  =\exp\bigl((1/2-o(1))n\log n\bigr).
\end{equation}
The estimate is uniform over all lattice shapes satisfying \eqref{lb:eq:R-range} and all $0\le\eta\le n^{-4}$; the $o(1)$ term may depend on the fixed value of $D_{\mathrm{sep}}$.
\end{lemma}

\begin{proof}
Apply \cref{lb:lem:torus-packing} to $G_\Lambda^\eta$.
Since $\det\Lambda=\kappa_n a^n/R$ and $\delta\asymp\sqrt{\log n/n}$ uniformly in the stated range,
\[
  \log\frac{a^n}{R\delta^n}
  =\frac n2\log n-\frac n2\log\log n+O(n)-O(\log n).
\]
This proves \eqref{lb:eq:number-targets}.
\end{proof}

\section{Poisson witnesses and the common-lens hierarchy}\label{lb:sec:clusters}

For a fixed good target, each residual vector in the outer shell determines a cap of Poisson directions that excludes that residual from the random body.
We group nearby residual caps into blocks and place a common measurable slot inside every cap in a block.
An injective Poisson point in each slot then certifies that the target is uncovered.
The hierarchy below controls dependencies among the slots of one target, while the separation from \cref{lb:lem:separated-targets} controls dependencies between different targets.
We finally apply the continuous Janson inequality to many targets for which the products of the slot measures are comparable.

We make the constant choices as follows.
The body and shell parameters $a$, $A$, and $\beta$ are fixed first.
The geometric estimates in \cref{lb:lem:cap-affinity,lb:lem:metric-lens} determine their absolute constants, including $A_0$.
We then fix $\varepsilon$, take $c_{\mathrm{sh}}$ and $B_0$ from \cref{lb:lem:row-bound}, and choose a fixed $B\ge B_0$; this determines the anchor exponent $C_{\mathrm{anc}}$.
These choices are independent of the tail exponent $J$.
After fixing $J>0$, we choose $D_{\mathrm{sep}}=D_{\mathrm{sep}}(J)$ sufficiently large for the cross-target estimate and then choose $c_0=c_0(J)$ sufficiently small.
No constant depends on $n$, the shape of the lattice, or $\eta$.

Fix a good target $y$ with $m(y)\ge1$, put $m\coloneqq m(y)$, and list its residuals as
\[
  X_y^\eta=\{x_1,\ldots,x_m\}.
\]
For $i\in[m]$, put
\[
  C_i\coloneqq C_{x_i}^\eta
  \quad\text{and}\quad
  c_i\coloneqq\mu(C_i)=c_\eta(x_i).
\]
The good-target estimates give
\begin{equation}\label{lb:eq:hier-basic}
  m\le CR,\qquad m\le n^2,\qquad
  \sum_{i=1}^m q(c_i)\le CR,\quad\text{and}\quad
  n^{-C_A}\le c_i\le C_A.
\end{equation}
The inequality $m\le n^2$ follows from $m\le CR$ and $R\le c_0 n\log n$ for sufficiently large $n$.

Put
\[
  \bar c_i\coloneqq\min\{c_i,1/2\}
  \quad\text{and}\quad
  s_i\coloneqq-\log\bar c_i.
\]
The truncation keeps the slot scale at most $1/2$, while $s_i$ records the rarity of the cap.
The metric below combines angular energy with mismatch between these logarithmic intensities.

View $e_i\coloneqq x_i/|x_i|$ as a point of real projective space.
Define
\begin{equation}\label{lb:eq:parameter-metric}
  \mathfrak d(i,j)\coloneqq\sqrt n\,\angle_{\mathrm{proj}}(e_i,e_j)
         +\sqrt{|s_i-s_j|}.
\end{equation}
The projective angle is the geodesic metric on real projective space.
The second summand is also a metric because $\sqrt{x+y}\le\sqrt x+\sqrt y$ for nonnegative $x$ and $y$.
Hence $\mathfrak d$ is a metric, and
\begin{equation}\label{lb:eq:metric-energy}
  n\angle_{\mathrm{proj}}(e_i,e_j)^2+|s_i-s_j|
  \le \mathfrak d(i,j)^2
  \le2n\angle_{\mathrm{proj}}(e_i,e_j)^2+2|s_i-s_j|.
\end{equation}
For every nonempty set $S\subseteq[m]$, let
\begin{equation}\label{lb:eq:cluster-radius}
  r(S)\coloneqq\min_{a\in S}\max_{i\in S}\mathfrak d(a,i)^2.
\end{equation}

The next proposition is the deterministic core of the fixed-target argument.
It partitions the residual caps into blocks and assigns a common-lens slot to each block.
Later, an injective choice of one Poisson point from every slot will certify that the target is a hole.
Estimate \eqref{lb:eq:hier-activity-cost} bounds the negative logarithm of the product of the slot measures, \eqref{lb:eq:hier-loop-cost} controls diagonal identifications, and \eqref{lb:eq:row-product} controls off-diagonal identifications after the normalization in \eqref{lb:eq:normalized-overlap}.
The last two estimates feed into the self-diagram bound in \cref{lb:lem:self-diagram}.

\begin{proposition}\label{lb:prop:hierarchical-slots}
For all sufficiently large $n$, every good target $y$ with $m(y)\ge1$ has a partition $\cP_y$ of $[m(y)]$.
For every block $S\in\cP_y$, there is a measurable slot
\[
  E_S\subseteq\bigcap_{i\in S}C_i,
  \quad\text{where}\quad 0<\lambda_S\coloneqq\mu(E_S)\le\frac12,
\]
such that
\begin{align}
  \sum_{S\in\cP_y}-\log\lambda_S&\le CR,
    \label{lb:eq:hier-activity-cost}\\
  \sum_{S\in\cP_y}\log(1+\lambda_S^{-1})&\le CR.
    \label{lb:eq:hier-loop-cost}
\end{align}
For distinct blocks $S$ and $T$, put
\begin{equation}\label{lb:eq:normalized-overlap}
  w_{ST}\coloneqq
  \frac{\mu(E_S\cap E_T)}
  {\sqrt{\lambda_S(1+\lambda_S)\lambda_T(1+\lambda_T)}}.
\end{equation}
Then
\begin{equation}\label{lb:eq:row-product}
  \sum_{S\in\cP_y}
  \log\left(1+\sum_{T\in\cP_y\setminus\{S\}}w_{ST}\right)
  \le CR.
\end{equation}
Every slot has an anchor $a(S)\in S$ such that
\begin{equation}\label{lb:eq:slot-anchor}
  E_S\subseteq C_{a(S)}
  \quad\text{and}\quad
  \lambda_S\ge n^{-C_{\mathrm{anc}}}c_{a(S)}.
\end{equation}
All constants are uniform in the lattice, the good target, and $0\le\eta\le n^{-4}$.
\end{proposition}

\subsection{A sharp affinity estimate for caps}

The coarse two-cap estimate in \cref{lb:lem:caps} has a polynomial prefactor, which is too expensive for the row sum below.
The metric \eqref{lb:eq:parameter-metric} also records radial intensity and yields a constant-prefactor estimate.

\begin{lemma}\label{lb:lem:cap-affinity}
For all sufficiently large $n$ and all caps $C_i$ and $C_j$ belonging to one good target,
\begin{equation}\label{lb:eq:cap-affinity}
  \frac{\mu(C_i\cap C_j)}{\sqrt{\bar c_i\bar c_j}}
  \le C\exp\{-c \mathfrak d(i,j)^2\}.
\end{equation}
The constants are uniform for $0\le\eta\le n^{-4}$.
\end{lemma}

\begin{proof}
We first prove an oriented-cap estimate.
Let
\[
  H(e,\tau)\coloneqq\{u\colon\ip{u}{e}>\tau\}
  \quad\text{and}\quad
  H(f,\xi)\coloneqq\{u\colon\ip{u}{f}>\xi\},
\]
where $\tau$ and $\xi$ lie in a fixed compact interval $[\tau_0,\tau_1]\subsetneq(0,1)$.
Put
\[
  p(\zeta)\coloneqq\sigma(H(e,\zeta)),\qquad
  \Delta\coloneqq|\log p(\tau)-\log p(\xi)|,
  \quad\text{and}\quad t\coloneqq n\psi^2,
\]
where $\psi$ is the angle between $e$ and $f$.

\medskip\noindent\emph{One-cap estimates.}
Let
\[
  c_n\coloneqq\frac{\Gamma(n/2)}{\sqrt\pi\,\Gamma((n-1)/2)}.
\]
Inserting the factor $x$ into the one-cap integral gives
\begin{equation}\label{lb:eq:aff-onecap}
  \frac{c_n}{n-1}(1-\zeta^2)^{(n-1)/2}
  \le p(\zeta)
  \le\frac{c_n}{(n-1)\zeta}(1-\zeta^2)^{(n-1)/2}.
\end{equation}
Since $p'(\zeta)=-c_n(1-\zeta^2)^{(n-3)/2}$, we obtain
\begin{equation}\label{lb:eq:aff-hazard}
  \frac{(n-1)\zeta}{1-\zeta^2}
  \le-\frac{\mathrm d}{\mathrm d\zeta}\log p(\zeta)
  \le\frac{n-1}{1-\zeta^2}.
\end{equation}
Consequently,
\begin{equation}\label{lb:eq:aff-radial}
  cn|\tau-\xi|\le\Delta\le Cn|\tau-\xi|
\end{equation}
uniformly on the threshold interval.

We claim that
\begin{equation}\label{lb:eq:oriented-affinity}
  \frac{\sigma(H(e,\tau)\cap H(f,\xi))}
       {\sqrt{p(\tau)p(\xi)}}
  \le C e^{-c(t+\Delta)}.
\end{equation}
The intersection is contained in the smaller cap, and hence the left side is at most $e^{-\Delta/2}$.
Choose a sufficiently large absolute constant $K$.
If $t\le K(\Delta+1)$, then this trivial estimate proves \eqref{lb:eq:oriented-affinity} after adjusting $c$ and $C$.
We may therefore assume that
\begin{equation}\label{lb:eq:aff-separated-regime}
  t>K(\Delta+1).
\end{equation}
By \eqref{lb:eq:aff-radial},
\begin{equation}\label{lb:eq:threshold-vs-angle}
  |\tau-\xi|\le C\frac{\Delta}{n}
  \le\frac{C}{K}\psi^2.
\end{equation}

\medskip\noindent\emph{Two-cap optimization.}
Use orthonormal coordinates on $\operatorname{span}\{e,f\}$ in which
\[
  e=(1,0)
  \quad\text{and}\quad
  f=(\cos\psi,\sin\psi).
\]
The intersection of the two boundary lines is
\begin{equation}\label{lb:eq:v-star}
  v_*\coloneqq
  \left(\tau,\frac{\xi-\tau\cos\psi}{\sin\psi}\right).
\end{equation}
Writing $v_*=\omega_e e+\omega_f f$ gives
\begin{equation}\label{lb:eq:lagrange-multipliers}
  \omega_e\coloneqq\frac{\tau-\xi\cos\psi}{\sin^2\psi}
  \quad\text{and}\quad
  \omega_f\coloneqq\frac{\xi-\tau\cos\psi}{\sin^2\psi}.
\end{equation}
We show that, whenever the cap intersection is nonempty, there is a uniform constant $c>0$ such that
\begin{equation}\label{lb:eq:quadratic-expansion}
  |v|^2\ge |v_*|^2+c(\delta_e+\delta_f)
\end{equation}
for every feasible $v$, where
\[
  \delta_e\coloneqq\ip{v}{e}-\tau\ge0
  \quad\text{and}\quad
  \delta_f\coloneqq\ip{v}{f}-\xi\ge0,
\]
and that
\begin{equation}\label{lb:eq:vstar-gain}
  |v_*|^2\ge\max\{\tau^2,\xi^2\}+c\min\{\psi^2,1\}.
\end{equation}

Fix a small constant $\psi_0>0$.
If $0<\psi\le\psi_0$, then
\[
  \omega_e=\frac{\tau}{1+\cos\psi}
    +\frac{(\tau-\xi)\cos\psi}{\sin^2\psi}
  \quad\text{and}\quad
  \omega_f=\frac{\xi}{1+\cos\psi}
    +\frac{(\xi-\tau)\cos\psi}{\sin^2\psi}.
\]
By \eqref{lb:eq:threshold-vs-angle}, choosing $K$ sufficiently large gives $\omega_e,\omega_f\ge c_{\mathrm{lag}}>0$.
Since $v-v_*$ has inner products $\delta_e$ and $\delta_f$ with $e$ and $f$, respectively,
\[
  |v|^2
  =|v_*|^2+2\omega_e\delta_e+2\omega_f\delta_f+|v-v_*|^2,
\]
which proves \eqref{lb:eq:quadratic-expansion}.
If $\tau\ge\xi$, then
\[
  \xi-\tau\cos\psi
  =\tau(1-\cos\psi)-(\tau-\xi)\ge c\psi^2,
\]
and hence
\[
  |v_*|^2-\tau^2
  =\frac{(\xi-\tau\cos\psi)^2}{\sin^2\psi}
  \ge c\psi^2.
\]
The case $\xi\ge\tau$ is symmetric, so \eqref{lb:eq:vstar-gain} follows in the small-angle range.

Suppose now that $\psi\ge\psi_0$ and the cap intersection is nonempty.
For a point $u$ in the intersection,
\[
  2\tau_0\le\ip{u}{e+f}\le |e+f|=2\cos(\psi/2).
\]
Thus $\psi\le\pi-\delta_0$ for some fixed $\delta_0>0$.
When $\tau=\xi$, the quantities in \eqref{lb:eq:lagrange-multipliers} are positive and $|v_*|^2-\tau^2=\tau^2\tan^2(\psi/2)$ is bounded below on $[\psi_0,\pi-\delta_0]$.
The compactness of this parameter range, together with \eqref{lb:eq:threshold-vs-angle}, gives uniform positive lower bounds for $\omega_e$, $\omega_f$, and the norm gain when $K$ is sufficiently large.
This proves \eqref{lb:eq:quadratic-expansion} and \eqref{lb:eq:vstar-gain} in the remaining range.

\medskip\noindent\emph{Two-dimensional integration.}
The first two orthonormal coordinates of a uniform point of $\Sph^{n-1}$ have density
\[
  c_{n,2}(1-|v|^2)_+^{(n-4)/2},\quad\text{where}\quad
  c_{n,2}\coloneqq\frac{\Gamma(n/2)}{\pi\Gamma((n-2)/2)}
  =\frac{n-2}{2\pi}.
\]
The map $v\mapsto(\delta_e,\delta_f)$ has Jacobian determinant $\sin\psi$.
If $|v_*|\ge1$, then the cap intersection is empty.
Otherwise, let $Q\coloneqq1-|v_*|^2$.
By \eqref{lb:eq:quadratic-expansion},
\begin{align*}
  \sigma(H(e,\tau)\cap H(f,\xi))
  \le\frac{c_{n,2}}{\sin\psi}
    \int_0^\infty\!\int_0^\infty
    (Q-c(\delta_e+\delta_f))_+^{(n-4)/2}\,\dd\delta_e\,\dd\delta_f
  =\frac{c_{n,2}}{\sin\psi}
    \frac{4Q^{n/2}}{c^2n(n-2)}.
\end{align*}
Therefore
\begin{equation}\label{lb:eq:two-dimensional-endpoint}
  \sigma(H(e,\tau)\cap H(f,\xi))
  \le\frac{C}{n\sin\psi}(1-|v_*|^2)^{n/2}.
\end{equation}
Let $\zeta\coloneqq\max\{\tau,\xi\}$.
The norm gain \eqref{lb:eq:vstar-gain} and compactness of the threshold interval imply
\[
  \left(\frac{1-|v_*|^2}{1-\zeta^2}\right)^{n/2}
  \le e^{-cn\psi^2}.
\]
Dividing \eqref{lb:eq:two-dimensional-endpoint} by the smaller one-cap probability and using \eqref{lb:eq:aff-onecap}, we obtain
\[
  \frac{\sigma(H(e,\tau)\cap H(f,\xi))}
       {\min\{p(\tau),p(\xi)\}}
  \le\frac{C}{\sqrt n\sin\psi}e^{-ct}.
\]
In the small-angle range, \eqref{lb:eq:aff-separated-regime} gives $\sqrt n\sin\psi\asymp\sqrt t\ge c\sqrt K$.
In the large-angle range, $\sin\psi$ is bounded below.
Thus the last expression is at most $Ce^{-c't}$.
Multiplying by
\[
  \sqrt{\frac{\min\{p(\tau),p(\xi)\}}
              {\max\{p(\tau),p(\xi)\}}}=e^{-\Delta/2}
\]
proves \eqref{lb:eq:oriented-affinity}.

\medskip\noindent\emph{Absolute caps.}
Let $\tau_i\coloneqq(1+\eta)/|x_i|$.
Then $c_i=2Mp(\tau_i)$.
Among the four pairs of oriented halves of $C_i$ and $C_j$, two have angle $\phi=\angle_{\mathrm{proj}}(e_i,e_j)$, and the other two have angle $\pi-\phi\ge\pi/2$.
Summing \eqref{lb:eq:oriented-affinity} over these four pairs gives
\begin{equation}\label{lb:eq:absolute-affinity-untruncated}
  \frac{\mu(C_i\cap C_j)}{\sqrt{c_ic_j}}
  \le C\exp\{-c(n\phi^2+|\log c_i-\log c_j|)\}.
\end{equation}
Since $c_i\le C_A$, replacing $c_i$ by the previously defined $\bar c_i$ changes the denominator by at most a fixed factor.
Furthermore,
\begin{equation}\label{lb:eq:truncated-log-comparison}
  |\log c_i-\log c_j|\ge |s_i-s_j|-C.
\end{equation}
If neither intensity is truncated, this is equality.
If both are truncated, the right side is zero.
If exactly one is truncated, the two sides differ by at most $\log(2C_A)$.
Combining \eqref{lb:eq:absolute-affinity-untruncated}, \eqref{lb:eq:truncated-log-comparison}, and \eqref{lb:eq:metric-energy} proves \eqref{lb:eq:cap-affinity}.
\end{proof}

\subsection{Scale-sensitive common lenses}

Pairwise affinity is not enough to assign an entire cluster of caps to one Poisson point.
The next lemma gives a quantitative common intersection for caps that are close in the parameter metric.

\begin{lemma}\label{lb:lem:metric-lens}
There is an absolute constant $A_0$ with the following property.
Fix $B>0$, and let $n$ be sufficiently large in terms of $B$.
If $S\subseteq[m]$, $|S|=N\ge2$, and $r(S)\le B\log n$, then there exist an anchor $a(S)\in S$ and a measurable set $F_S\subseteq\bigcap_{i\in S}C_i$ such that
\begin{equation}\label{lb:eq:metric-lens}
  \mu(F_S)\ge
  \exp\{-A_0(1+r(S)+\log N)\}\bar c_{a(S)}.
\end{equation}
The anchor may be chosen to minimize $c_i$ over $i\in S$.
\end{lemma}

\begin{proof}
Choose $j_0\in S$ that attains the minimum in \eqref{lb:eq:cluster-radius}.
Since $r(S)\le B\log n=o(n)$, every projective direction in $S$ has, for sufficiently large $n$, a unique representative making angle at most $\sqrt{r(S)/n}<\pi/4$ with a fixed representative of $e_{j_0}$.
Use these representatives below.

Let $i_0\coloneqq a(S)$ minimize $c_i$ over $S$.
By the triangle inequality in projective space, every $e_i$ makes angle at most $\theta\coloneqq2\sqrt{r(S)/n}$ with $e_{i_0}$.
Hence
\[
  e_i=\cos\phi_i\,e_{i_0}+\sin\phi_i\,w_i,
  \quad\text{where}\quad w_i\perp e_{i_0},\quad |w_i|=1,
  \quad\text{and}\quad 0\le\phi_i\le\theta.
\]
Let $S_+\coloneqq\{i\in S\colon\phi_i>0\}$.
No tangential constraint is needed when $\phi_i=0$.

Put $\tau_i\coloneqq(1+\eta)/|x_i|$.
Since $c_i=2Mp(\tau_i)$ and $p$ is decreasing, the choice of the anchor gives
\[
  \tau\coloneqq\max_{i\in S}\tau_i=\tau_{i_0}.
\]
Choose sufficiently large absolute constants $C_s$ and $C$, in this order, and set
\[
  s\coloneqq C_s\sqrt{\frac{\log(2N)}n}
  \quad\text{and}\quad
  h\coloneqq C(\theta s+\theta^2).
\]
Because $N\le m\le n^2$ and $r(S)\le B\log n$, we have $s=o(1)$ and $h=O_B(\log n/n)$.
Thus $\tau+h$ remains in a fixed compact subinterval of $(0,1)$ for all sufficiently large $n$.

Condition a uniform vector $u\in\Sph^{n-1}$ on $t=\ip{u}{e_{i_0}}$, and write
\[
  u=t e_{i_0}+\sqrt{1-t^2}\,w,
  \quad\text{where}\quad w\in\Sph^{n-2}.
\]
For every fixed $w_i$, spherical concentration gives
\[
  \Prb_w\left[\ip{w}{w_i}<-s\right]\le e^{-cns^2}.
\]
Choose $C_s$ so that $Ne^{-cns^2}\le1/2$.
A union bound yields
\begin{equation}\label{lb:eq:tangential-good}
  \Prb_w\left[\ip{w}{w_i}\ge-s
  \text{ for every }i\in S_+\right]\ge\frac12.
\end{equation}
Define
\begin{equation}\label{lb:eq:explicit-common-lens}
  F_S\coloneqq\left\{
  u=t e_{i_0}+\sqrt{1-t^2}\,w\colon
  t>\tau+h\text{ and }
  \ip{w}{w_i}\ge-s\text{ for every }i\in S_+
  \right\}.
\end{equation}
For $u\in F_S$ and $i\in S$,
\begin{align*}
  \ip{u}{e_i}
  =t\cos\phi_i
    +\sqrt{1-t^2}\sin\phi_i\ip{w}{w_i}
  \ge t\cos\phi_i-s\sin\phi_i
  \ge\tau+h-C(\phi_i^2+s\phi_i)
  >\tau_i.
\end{align*}
The last inequality follows from the choice of $C$ in the definition of $h$.
Thus $F_S$ lies in one oriented half of every $C_i$, and hence $F_S\subseteq\bigcap_{i\in S}C_i$.

Integrating first in $w$ and using \eqref{lb:eq:tangential-good}, we obtain
\[
  \sigma(F_S)\ge\frac12p(\tau+h).
\]
The oriented anchor cap has measure $p(\tau)$, and $c_{i_0}=2Mp(\tau)$.
Therefore
\begin{equation}\label{lb:eq:lens-measure-before-hazard}
  \mu(F_S)\ge\frac{c_{i_0}}{4}\frac{p(\tau+h)}{p(\tau)}.
\end{equation}
The hazard estimate \eqref{lb:eq:aff-hazard} gives
\[
  \frac{p(\tau+h)}{p(\tau)}\ge e^{-Cnh}.
\]
Moreover,
\[
  nh\le C\bigl(r(S)+\sqrt{r(S)\log N}\bigr)
  \le C'\bigl(1+r(S)+\log N\bigr).
\]
Combining the last two estimates with \eqref{lb:eq:lens-measure-before-hazard} and $\bar c_{i_0}\le c_{i_0}$ proves \eqref{lb:eq:metric-lens}.
The resulting constant $A_0$ is absolute and does not depend on $B$.
\end{proof}

Increase $A_0$ once, if necessary, so that $A_0\ge1$; this choice is independent of $\varepsilon$.
Then fix a sufficiently small absolute constant $\varepsilon>0$.
Define
\begin{equation}\label{lb:eq:lens-penalty}
  \ell(S)\coloneqq A_0(1+r(S)+\log|S|)\quad\text{when }|S|\ge2,
  \quad\text{and}\quad \ell(\{i\})\coloneqq0.
\end{equation}
Let
\begin{equation}\label{lb:eq:R-star}
  R_*\coloneqq B\log n,
\end{equation}
where $B$ is the fixed constant from the hierarchy at the beginning of the section.
A nonempty set $S\subseteq[m]$ is \emph{economical} if it is a singleton or if
\begin{equation}\label{lb:eq:economical}
  r(S)\le R_*
  \quad\text{and}\quad
  \ell(S)\le\varepsilon|S|.
\end{equation}

Start with the partition into singletons.
Whenever at least two current blocks have an economical union, replace those blocks by their union.
Each merge decreases the number of blocks, so the procedure terminates.
Let $\cP_y$ be a terminal partition.
Every block of $\cP_y$ is economical, and
\begin{equation}\label{lb:eq:terminal}
  \text{no union of at least two blocks of $\cP_y$ is economical.}
\end{equation}

For a nonsingleton block $S\in\cP_y$, apply \cref{lb:lem:metric-lens}.
By nonatomicity, choose a measurable subset $E_S\subseteq F_S$ with
\begin{equation}\label{lb:eq:exact-lambda}
  \lambda_S\coloneqq\mu(E_S)=e^{-\ell(S)}\bar c_{a(S)}.
\end{equation}
For a singleton $S=\{i\}$, choose $E_S\subseteq C_i$ with $\mu(E_S)=\bar c_i$, and put $a(S)\coloneqq i$.
Thus \eqref{lb:eq:exact-lambda} holds for every block.

\subsection{Terminal local density and weighted rows}

Terminality prevents too many low-penalty blocks from having anchors near one point.
The resulting local-count estimate will be converted into the row-sum bound needed for the self-diagram calculation.

\begin{lemma}\label{lb:lem:terminal-density}
Let $T_{\mathrm{loc}},H\ge0$, let $N_0\ge1$, let $z\in[m]$, and let $S_1,\ldots,S_{N_0}$ be distinct blocks of $\cP_y$ such that $\mathfrak d(a(S_j),z)^2\le T_{\mathrm{loc}}$ and $\ell(S_j)\le H$ for every $j\in[N_0]$.
If
\begin{equation}\label{lb:eq:local-range}
  8T_{\mathrm{loc}}+8H/A_0\le R_*,
\end{equation}
then
\begin{equation}\label{lb:eq:local-count}
  N_0\le C(1+T_{\mathrm{loc}}+H).
\end{equation}
\end{lemma}

\begin{proof}
Let $U\coloneqq S_1\cup\cdots\cup S_{N_0}$ and $M_U\coloneqq|U|$.
For each $j$, choose $i_j\in S_j$ that attains the minimum in the definition of $r(S_j)$.
If $i\in S_j$, then the triangle inequality gives
\[
  \mathfrak d(i,a(S_j))
  \le \mathfrak d(i,i_j)+\mathfrak d(i_j,a(S_j))
  \le2\sqrt{r(S_j)}
  \le2\sqrt{H/A_0}.
\]
Using $a(S_1)$ as a center and applying the triangle inequality through $z$ and $a(S_j)$, we obtain
\[
  \mathfrak d(a(S_1),i)\le2\sqrt{T_{\mathrm{loc}}}+2\sqrt{H/A_0}.
\]
Consequently,
\[
  r(U)\le8T_{\mathrm{loc}}+8H/A_0\le R_*.
\]
The conclusion is immediate when $N_0=1$.
Suppose that $N_0\ge2$.
Then
\[
  \ell(U)
  \le A_0+8A_0T_{\mathrm{loc}}+8H+A_0\log M_U.
\]
Choose $C$ so large that
\[
  A_0+8A_0T_{\mathrm{loc}}+8H\le\frac\varepsilon2M_U
\]
whenever $M_U\ge N_0>C(1+T_{\mathrm{loc}}+H)$, and enlarge $C$ so that $A_0\log u\le\varepsilon u/2$ for every $u\ge C$.
If $N_0>C(1+T_{\mathrm{loc}}+H)$, then $\ell(U)\le\varepsilon|U|$.
Thus $U$ is an economical union of at least two terminal blocks, contrary to \eqref{lb:eq:terminal}.
This contradiction proves \eqref{lb:eq:local-count}.
\end{proof}

For every block $S\in\cP_y$, put
\begin{equation}\label{lb:eq:L-S}
  L_S\coloneqq-\log\lambda_S=s_{a(S)}+\ell(S).
\end{equation}
By \cref{lb:lem:cap-affinity}, slot containment, and \eqref{lb:eq:exact-lambda},
\begin{equation}\label{lb:eq:angular-weight}
  w_{ST}\le
  C\exp\left\{\frac{\ell(S)+\ell(T)}2
  -c \mathfrak d(a(S),a(T))^2\right\}.
\end{equation}
The elementary bound $\mu(E_S\cap E_T)\le\min\{\lambda_S,\lambda_T\}$ also gives
\begin{equation}\label{lb:eq:intensity-weight}
  w_{ST}\le e^{-|L_S-L_T|/2}\le1.
\end{equation}

\begin{lemma}\label{lb:lem:row-bound}
There are absolute constants $c_{\mathrm{sh}}>0$ and $B_0>0$ such that, for every fixed $B\ge B_0$ and all sufficiently large $n$, every $S\in\cP_y$ with $L_S\le c_{\mathrm{sh}}R_*$ satisfies
\begin{equation}\label{lb:eq:row-bound}
  \sum_{T\in\cP_y\setminus\{S\}}w_{ST}\le C(1+L_S).
\end{equation}
\end{lemma}

\begin{proof}
Every positive weight $w_{ST}$ lies in a unique interval $\bigl(e^{-(k+1)},e^{-k}\bigr]$.
For $k\ge0$, let $\mathcal N_k(S)$ be the set of blocks $T\ne S$ whose weight lies in this interval.
By \eqref{lb:eq:intensity-weight},
\[
  |L_S-L_T|\le2(k+1)
  \quad\text{and}\quad
  \ell(T)\le L_T\le L_S+2(k+1).
\]
The lower bound $w_{ST}>e^{-(k+1)}$, together with \eqref{lb:eq:angular-weight}, implies
\[
  \mathfrak d(a(S),a(T))^2\le C_2(1+L_S+k).
\]
Apply \cref{lb:lem:terminal-density} with
\[
  z=a(S),\qquad T_{\mathrm{loc}}=C_2(1+L_S+k),
  \quad\text{and}\quad H=L_S+2(k+1).
\]
There is a sufficiently small absolute constant $c_1>0$ such that the range condition \eqref{lb:eq:local-range} holds whenever $L_S+k+1\le c_1R_*$.
Hence
\begin{equation}\label{lb:eq:band-count}
  |\mathcal N_k(S)|\le C_3(1+L_S+k)
\end{equation}
in this range.

Choose $c_{\mathrm{sh}}\le c_1/4$ and let $K\coloneqq\lfloor c_1R_*/2\rfloor$.
Then \eqref{lb:eq:band-count} holds for $0\le k\le K$.
All remaining weights are at most $e^{-K}$.
Since $|\cP_y|\le m\le n^2$, there is an absolute $B_0>0$ such that every fixed $B\ge B_0$ ensures that $me^{-K}\le1$ for all sufficiently large $n$.
Therefore
\[
  \sum_{T\ne S}w_{ST}
  \le C\sum_{k\ge0}e^{-k}(1+L_S+k)+1
  \le C'(1+L_S).\qedhere
\]
\end{proof}

\begin{proof}[Proof of \cref{lb:prop:hierarchical-slots}]
For $u>0$,
\[
  -\log\min\{u,1/2\}\le C(1+q(u)).
\]
The anchors of distinct blocks are distinct.
Hence
\begin{equation}\label{lb:eq:anchor-cost}
  \sum_{S\in\cP_y}s_{a(S)}
  \le\sum_{i=1}^m s_i
  \le C\left(m+\sum_iq(c_i)\right)
  \le CR.
\end{equation}
Every terminal block is economical, so
\begin{equation}\label{lb:eq:penalty-cost}
  \sum_{S\in\cP_y}\ell(S)
  \le\varepsilon\sum_{S\in\cP_y}|S|
  =\varepsilon m.
\end{equation}
Equations \eqref{lb:eq:L-S}, \eqref{lb:eq:anchor-cost}, and \eqref{lb:eq:penalty-cost} prove \eqref{lb:eq:hier-activity-cost}.
Also,
\[
  \log(1+\lambda_S^{-1})\le L_S+\log 2,
\]
and hence \eqref{lb:eq:hier-loop-cost} follows.

To prove \eqref{lb:eq:row-product}, call a block shallow when $L_S\le c_{\mathrm{sh}}R_*$ and deep otherwise.
\Cref{lb:lem:row-bound} bounds the logarithmic row cost of a shallow block by $C(1+L_S)$.
For a deep block, \eqref{lb:eq:intensity-weight} and $|\cP_y|\le n^2$ give
\[
  \log\left(1+\sum_{T\ne S}w_{ST}\right)
  \le C\log n\le C_B L_S.
\]
Since $|\cP_y|\le m\le CR$ and $\sum_{S\in\cP_y}L_S\le CR$, summing over all blocks proves \eqref{lb:eq:row-product}.

Finally, every block satisfies $r(S)\le R_*$ and $|S|\le n^2$, so $\ell(S)\le C_B\log n$.
Since $\bar c_{a(S)}\ge c_{a(S)}/(2C_A)$, equation \eqref{lb:eq:exact-lambda} gives \eqref{lb:eq:slot-anchor} after enlarging the absolute exponent $C_{\mathrm{anc}}$.
\end{proof}

\subsection{Self and cross diagrams}

For a good target $y$ with $m(y)\ge1$, let
\[
  \nu_y\coloneqq\prod_{S\in\cP_y}\lambda_S,
\]
and let $Z_y$ count the injective tuples with one Poisson point in every final slot.
If $Z_y>0$, then every residual cap belonging to $y$ is hit.
More precisely, goodness excludes residuals of norm less than $r_0-\eta$.
A residual of norm between $r_0-\eta$ and $b+\eta$ is cut by the point chosen for the slot containing its cap, while a residual of norm greater than $b+\eta$ lies outside the outer ball of $\widehat K_{\Pi,\eta}$.
Therefore
\begin{equation}\label{lb:eq:witness-hole}
  Z_y>0\quad\Longrightarrow\quad
  y\notin\widehat K_{\Pi,\eta}+\Lambda.
\end{equation}
By \eqref{lb:eq:hier-activity-cost},
\begin{equation}\label{lb:eq:activity-range}
  e^{-CR}\le\nu_y\le1.
\end{equation}
Let $D_{yy}$ denote the self-diagram sum from \cref{lb:lem:poisson-janson}.

The calculation below uses \eqref{lb:eq:hier-loop-cost} for diagonal loops and \eqref{lb:eq:row-product}, after \eqref{lb:eq:normalized-overlap}, for the off-diagonal partial matching.

\begin{lemma}\label{lb:lem:self-diagram}
For every good target $y$ with $m(y)\ge1$,
\begin{equation}\label{lb:eq:self-diagram-bound}
  \frac{D_{yy}}{\nu_y^2}\le e^{CR}.
\end{equation}
\end{lemma}

\begin{proof}
For every block $S\in\cP_y$, put $\delta_S\coloneqq1+\lambda_S^{-1}$.
For distinct blocks $S$ and $T$, put
\[
  h_{ST}\coloneqq\frac{\mu(E_S\cap E_T)}{\lambda_S\lambda_T}.
\]
We first fix the off-diagonal identifications between two copies of the witness and then sum over diagonal loops on the compatible, as-yet-unmatched slots.
Consider two copies of $\cP_y$.
Fix an off-diagonal partial matching $P$ between them, and let $\deg_P(S)\in\{0,1,2\}$ record whether the block $S$ is used on the left, on the right, or on both sides.
Summing all diagonal loops compatible with $P$ gives
\[
  \left(\prod_S \delta_S\right)
  \left(\prod_{(S,T)\in P}h_{ST}\right)
  \left(\prod_{S\colon\deg_P(S)>0}\delta_S^{-1}\right).
\]
Whenever $\deg_P(S)\in\{1,2\}$, we have $\delta_S^{-1}\le\delta_S^{-\deg_P(S)/2}$.
The preceding expression is at most
\[
  \left(\prod_S \delta_S\right)
  \prod_{(S,T)\in P}\frac{h_{ST}}{\sqrt{\delta_S\delta_T}}
  =\left(\prod_S \delta_S\right)
   \prod_{(S,T)\in P}w_{ST}.
\]
Summing over all off-diagonal partial matchings and adding the empty matching only enlarges the sum.
Hence
\[
  \frac{D_{yy}}{\nu_y^2}
  \le
  \prod_{S\in\cP_y}(1+\lambda_S^{-1})
  \sum_{\substack{P\text{ an off-diagonal}\\
                       \text{partial matching}}}
  \prod_{(S,T)\in P}w_{ST}.
\]
The first product is at most $e^{CR}$ by \eqref{lb:eq:hier-loop-cost}.
If we drop injectivity on the right copy, then
\[
  \sum_P\prod_{(S,T)\in P}w_{ST}
  \le\prod_{S\in\cP_y}
  \left(1+\sum_{T\ne S}w_{ST}\right)
  \le e^{CR}
\]
by \eqref{lb:eq:row-product}.
Absorbing the two absolute constants in the exponent proves the lemma.
\end{proof}

Fix $J>0$, and now choose the separation constant $D_{\mathrm{sep}}=D_{\mathrm{sep}}(J)$ promised above.
Index the separated targets by $\alpha$.
For the target $y_\alpha$, write $x_{\alpha,i}$ for its residuals and $\cP_\alpha$, $E_{\alpha S}$, $\lambda_{\alpha S}$, and $a_\alpha(S)$ for its partition, slots, slot measures, and anchors.
In the estimate below, $S\in\cP_\alpha$ and $T\in\cP_\beta$.
The constants $A_0$, $B$, and the anchor exponent $C_{\mathrm{anc}}$ are already fixed.
If slots belong to distinct targets, then slot containment, \eqref{lb:eq:cap-intensity-range}, and \cref{lb:lem:caps}\textup{(iii)} give
\begin{align}
  \frac{\mu(E_{\alpha S}\cap E_{\beta T})}
       {\lambda_{\alpha S}\lambda_{\beta T}}
  \le
  \frac{\mu(C_{x_{\alpha,a_\alpha(S)}}^\eta\cap
             C_{x_{\beta,a_\beta(T)}}^\eta)}
       {n^{-2C_{\mathrm{anc}}}c_\eta(x_{\alpha,a_\alpha(S)})
                    c_\eta(x_{\beta,a_\beta(T)})}
  \le n^{-J-5}.\label{lb:eq:cross-slot-overlap}
\end{align}
Here $D_{\mathrm{sep}}$ is chosen sufficiently large.
To verify the last inequality, apply \cref{lb:lem:caps}\textup{(iii)} with exponent $J+2C_{\mathrm{anc}}+C_A+5$.
Since each target has at most $n^2$ slots,
\begin{equation}\label{lb:eq:cross-rho-small}
  \rho_{\alpha\beta}\le n^{-J}
  \quad\text{for}\quad \alpha\ne\beta.
\end{equation}

This completes the choice of $D_{\mathrm{sep}}$.
All constants entering the self-diagram estimate are already fixed; the only remaining choice is the upper-range constant $c_0$, which is made in the proof of \cref{lb:prop:fixed-lattice}.
Thus there is no circular dependence.

\begin{proposition}\label{lb:prop:fixed-lattice}
For every fixed $J>0$, there are constants $c=c(J)>0$ and $c_0=c_0(J)>0$ such that, for all sufficiently large $n$, the following holds.
If $c_*\le R\le c_0F_n$, then, uniformly for $0\le\eta\le n^{-4}$,
\begin{equation}\label{lb:eq:fixed-lattice-probability}
  \Prb\left[\widehat K_{\Pi,\eta}+\Lambda=\R^n\right]
  \le e^{-cn^J}.
\end{equation}
\end{proposition}

\begin{proof}
Fix $J$, and choose $D_{\mathrm{sep}}$ as above.
A good target with $m(y)=0$ is a deterministic hole, so we may assume that every target selected below has at least one residual.

Let $y_1,\ldots,y_T$ be the separated targets from \cref{lb:lem:separated-targets}.
For a constant $C_T=C_T(J)>0$,
\begin{equation}\label{lb:eq:explicit-target-count}
  \log T\ge
  \frac n2\log n-\frac n2\log\log n-C_T n
\end{equation}
for all sufficiently large $n$.
By \eqref{lb:eq:activity-range}, all values $-\log\nu_y$ lie in an interval of length at most $CR$.
Partition this interval into subintervals of length $\log 2$.
One bin contains at least
\[
  T_{\mathrm{bin}}\ge\frac{T}{C(R+1)}
\]
targets whose activities differ by a factor of less than two.
Since $R\le c_0 n\log n$, taking $c_0$ sufficiently small and then $n$ sufficiently large gives
\begin{equation}\label{lb:eq:target-bin-lower}
  \log T_{\mathrm{bin}}\ge\frac13n\log n.
\end{equation}

Let $\nu$ be the smallest activity in this bin, and select $N$ targets from it.
Their activities satisfy $\nu\le\nu_\alpha<2\nu$.
Let $C_{\mathrm{sd}}$ be the constant in \cref{lb:lem:self-diagram}.
Then
\[
  \sum_\alpha D_{\alpha\alpha}
  \le4N\nu^2e^{C_{\mathrm{sd}}R}.
\]
For distinct targets, \eqref{lb:eq:cross-diagram} and \eqref{lb:eq:cross-rho-small}, together with $e^x-1\le2x$ for $0\le x\le1$, give
\[
  \sum_{\alpha\ne\beta}D_{\alpha\beta}
  \le C N^2\nu^2n^{-J}.
\]
Also, $(\sum_\alpha\nu_\alpha)^2\ge N^2\nu^2$.
Applying \cref{lb:lem:poisson-janson}, we obtain
\begin{equation}\label{lb:eq:janson-exponent-fixed}
  \Prb\left[Z_\alpha=0\text{ for every selected }\alpha\right]
  \le
  \exp\left[-c\min\left\{
  Ne^{-C_{\mathrm{sd}}R},n^J\right\}\right].
\end{equation}

Choose
\begin{equation}\label{lb:eq:chosen-target-number}
  N\coloneqq\left\lceil e^{C_{\mathrm{sd}}R}n^J\right\rceil.
\end{equation}
This choice is available when $c_0$ is sufficiently small.
Indeed,
\[
  \log N
  \le C_{\mathrm{sd}}c_0 n\log n+J\log n+O(1).
\]
Choose $c_0$ so that $C_{\mathrm{sd}}c_0<1/4$.
For sufficiently large $n$, the right side is less than $(1/3)n\log n$, and \eqref{lb:eq:target-bin-lower} gives $N\le T_{\mathrm{bin}}$.
Substitution into \eqref{lb:eq:janson-exponent-fixed} yields the bound $e^{-cn^J}$.

If $\widehat K_{\Pi,\eta}+\Lambda$ covers, then no selected target has a witness, by \eqref{lb:eq:witness-hole}.
The covering event is therefore contained in the event on the left side of \eqref{lb:eq:janson-exponent-fixed}.
This proves \eqref{lb:eq:fixed-lattice-probability}.
\end{proof}

\section{Uniformity over all lattices}\label{lb:sec:net}

We now pass from one fixed lattice to all lattices by a finite-net argument.
Any lattice covering $K_\Pi$ has covering radius at most $b$ and hence admits a short basis.
The determinant range then gives a quantitative lower bound on the least singular value of a basis matrix, allowing us to round it to a finite matrix grid.
A hole for the enlarged body is stable under the resulting near-identity linear perturbation.
We take $J=6$ in \cref{lb:prop:fixed-lattice}; the resulting probability $e^{-cn^6}$ dominates the grid size $\exp(O(n^3\log n))$.

For a lattice $\Lambda$, let $\covrad(\Lambda)$ denote its Euclidean covering radius.
For $1\le i\le n$, let $\lambda_i(\Lambda)$ denote its $i$th successive minimum with respect to $B_2^n$.
All unqualified matrix norms in this section are Euclidean operator norms; $\|\cdot\|_{\mathrm F}$ denotes the Frobenius norm.

\begin{lemma}\label{lb:lem:short-basis}
If $\covrad(\Lambda)\le b$, then $\Lambda$ has a basis $g_1,\ldots,g_n$ such that
\[
  |g_i|\le Cnb\quad\text{for}\quad 1\le i\le n.
\]
\end{lemma}

\begin{proof}
Jarn{\'{\i}}k's inequality~\cite{Jarnik1941} gives $\lambda_n(\Lambda)\le2\covrad(\Lambda)$; a modern formulation appears in~\cite[Eq.~(12)]{CelayaHenk2023}.
The Mahler basis estimate~\cite[Eq.~(2.5)]{FreyerHenk2022} gives a basis satisfying
\[
  |g_i|\le\max\{1,i/2\}\lambda_i(\Lambda).
\]
Since $\lambda_i(\Lambda)\le\lambda_n(\Lambda)$, the desired bound follows.
\end{proof}

The enlarged body from \eqref{lb:eq:enlarged-body} provides the buffer needed to transfer a hole between nearby lattices.

\begin{lemma}\label{lb:lem:stable-hole}
Let $0\le\eta\le n^{-4}$.
Let $G$ be an invertible matrix, let $\Lambda\coloneqq G\Z^n$, and let $\mathcal T$ be an invertible linear map.
Suppose that $y\notin\widehat K_{\Pi,\eta}+\Lambda$.
If
\[
  \|\mathcal T-I\|\le\frac{\eta}{4(b+\eta)},
\]
then $\mathcal T y\notin K_\Pi+\mathcal T\Lambda$.
\end{lemma}

\begin{proof}
Fix $z\in\Z^n$, and put $x\coloneqq y-Gz$.
Since $x\notin\widehat K_{\Pi,\eta}$, either $|x|>b+\eta$, or some $u\in\Pi$ satisfies $|\ip{u}{x}|>1+\eta$.

In the first case,
\[
  |\mathcal T x|\ge(1-\|\mathcal T-I\|)|x|
  >(1-\|\mathcal T-I\|)(b+\eta)
  \ge b+\frac{3\eta}{4}
  \ge b.
\]
In the second case, if $|x|\le b+\eta$, then
\[
  |\ip{u}{\mathcal T x}|
  \ge|\ip{u}{x}|-\|(\mathcal T-I)x\|
  >1+\eta-\frac\eta4
  =1+\frac{3\eta}{4}
  \ge1.
\]
When $|x|>b+\eta$, the first estimate applies.
Thus $\mathcal T x\notin K_\Pi$ for every $z\in\Z^n$.
Since $\mathcal T x=\mathcal T y-\mathcal T Gz$, the point $\mathcal T y$ is not covered by $K_\Pi+\mathcal T\Lambda$.
\end{proof}

Let $c_{\mathrm{fix}}$ be the upper-range constant supplied by \cref{lb:prop:fixed-lattice} for $J=6$.
Fix $c_{\mathrm{net}}>0$ so small that $2c_{\mathrm{net}}\le c_{\mathrm{fix}}$, and put
\begin{equation}\label{lb:eq:Dn}
  D_n\coloneqq\frac{\kappa_n a^n}{c_{\mathrm{net}}F_n}.
\end{equation}

\begin{proposition}\label{lb:prop:uniform-exclusion}
With probability $1-o(1)$, no lattice $\Lambda$ with $\det\Lambda\ge D_n$ covers $K_\Pi$.
\end{proposition}

\begin{proof}
Suppose that $K_\Pi+\Lambda$ covers.
Since $K_\Pi\subseteq bB_2^n$, the translates of $bB_2^n$ also cover.
Hence
\begin{equation}\label{lb:eq:covering-lattice-range}
  \covrad(\Lambda)\le b
  \quad\text{and}\quad
  \det\Lambda\le\kappa_n b^n.
\end{equation}
The determinant bound follows by integrating the covering multiplicity over a fundamental domain.

Now consider a lattice satisfying \eqref{lb:eq:covering-lattice-range} and $\det\Lambda\ge D_n$.
By \cref{lb:lem:short-basis}, it has a basis matrix $G'$ whose columns have norm at most $Cnb$.
Consequently,
\[
  \|G'\|\le\|G'\|_{\mathrm F}\le Cbn^{3/2}.
\]
Let $s_1(G')\ge\cdots\ge s_n(G')=s_{\min}(G')$ be the singular values.
Since
\[
  |\det G'|=\prod_{i=1}^ns_i(G')
  \le s_{\min}(G')\|G'\|^{n-1},
\]
we have
\begin{equation}\label{lb:eq:smin}
  s_{\min}(G')
  \ge\frac{D_n}{(Cbn^{3/2})^{n-1}}
  \eqqcolon\underline s_n.
\end{equation}
Stirling's formula gives $-\log\kappa_n=(n/2)\log n+O(n)$, and hence
\begin{equation}\label{lb:eq:smin-log}
  -\log\underline s_n=O(n\log n).
\end{equation}

Set $\eta\coloneqq n^{-5}$, and define the mesh
\begin{equation}\label{lb:eq:matrix-mesh}
  \xi_n\coloneqq\frac{\eta\underline s_n}{100n(b+\eta)}.
\end{equation}
Let $\Gamma_n\coloneqq\xi_n\Z\cap[-2Cnb,2Cnb]$ and let $\mathcal B_n^{\mathrm{grid}}\coloneqq\Gamma_n^{n\times n}$.
Retain the invertible matrices $G\in\mathcal B_n^{\mathrm{grid}}$ satisfying
\begin{equation}\label{lb:eq:grid-det-range}
  \frac{D_n}{2}\le|\det G|\le2\kappa_n b^n,
\end{equation}
and denote this finite set by $\mathcal G_n$.

One entry has at most $1+4Cnb/\xi_n$ possible values.
By \eqref{lb:eq:smin-log} and $\eta=n^{-5}$,
\[
  \log\left(1+\frac{4Cnb}{\xi_n}\right)=O(n\log n).
\]
Since there are $n^2$ entries,
\begin{equation}\label{lb:eq:net-size}
  \log|\mathcal G_n|=O(n^3\log n).
\end{equation}

Choose, once and for all, a nearest-point rounding map from $[-Cnb,Cnb]$ to $\Gamma_n$.
Apply it entrywise to $G'$ and call the rounded matrix $G$.
The larger interval in the definition of $\Gamma_n$ ensures that the rounded entries remain in the grid box.
Also,
\begin{equation}\label{lb:eq:rounding-error}
  \|G-G'\|\le\|G-G'\|_{\mathrm F}\le n\xi_n.
\end{equation}
By Weyl's inequality and \eqref{lb:eq:smin},
\[
  s_{\min}(G)
  \ge\underline s_n-n\xi_n
  \ge\frac12\underline s_n
\]
for sufficiently large $n$.
In particular, $G$ is invertible.

Put
\[
  E\coloneqq(G-G')(G')^{-1}.
\]
Equation \eqref{lb:eq:rounding-error} gives
\begin{equation}\label{lb:eq:E-bound}
  \|E\|\le\frac{n\xi_n}{\underline s_n}
  =\frac{\eta}{100(b+\eta)}<\frac12.
\end{equation}
Since $G=(I+E)G'$, every singular value of $I+E$ lies in $[1-\|E\|,1+\|E\|]$.
Therefore
\begin{equation}\label{lb:eq:determinant-stability}
  \left|\log\frac{|\det G|}{|\det G'|}\right|
  =|\log|\det(I+E)||
  \le2n\|E\|=O(n^{-4}).
\end{equation}
The inequalities $D_n\le|\det G'|\le\kappa_n b^n$ now imply \eqref{lb:eq:grid-det-range} for sufficiently large $n$.
Thus $G\in\mathcal G_n$.

The linear map from the grid lattice to the original lattice is
\[
  \mathcal T\coloneqq G'G^{-1}=(I+E)^{-1}.
\]
By \eqref{lb:eq:E-bound},
\begin{equation}\label{lb:eq:T-close}
  \|\mathcal T-I\|
  \le\frac{\|E\|}{1-\|E\|}
  \le2\|E\|
  \le\frac{\eta}{4(b+\eta)}.
\end{equation}

For every retained grid lattice $G\Z^n$, its parameter $R$ satisfies
\[
  \frac{e^{-\beta}}2
  =\frac{a^n}{2b^n}
  \le R
  \le\frac{2\kappa_n a^n}{D_n}
  =2c_{\mathrm{net}}F_n
  \le c_{\mathrm{fix}}F_n.
\]
In particular, $R\ge c_*$, so \cref{lb:prop:fixed-lattice} applies uniformly with $\eta=n^{-5}$.
By \eqref{lb:eq:net-size} and a union bound,
\begin{equation}\label{lb:eq:no-grid-cover}
  \Prb\left[\widehat K_{\Pi,\eta}+G\Z^n
  \text{ covers for some }G\in\mathcal G_n\right]
  \le\exp\bigl(O(n^3\log n)-cn^6\bigr)=o(1).
\end{equation}

On the complementary event, every retained grid lattice has a hole for $\widehat K_{\Pi,\eta}$.
Suppose to the contrary that an actual lattice $\Lambda\coloneqq G'\Z^n$ with $\det\Lambda\ge D_n$ covers $K_\Pi$.
The preceding rounding argument supplies a matrix $G\in\mathcal G_n$.
Choose $y\notin\widehat K_{\Pi,\eta}+G\Z^n$.
By \eqref{lb:eq:T-close} and \cref{lb:lem:stable-hole},
\[
  \mathcal T y\notin K_\Pi+\mathcal T G\Z^n=K_\Pi+G'\Z^n,
\]
contrary to the assumed covering.
This contradiction proves the proposition.
\end{proof}

\section{Proofs of the main results}\label{lb:sec:lower-main-proof}

\begin{proof}[Proof of \cref{thm:lower}]
Let $V$ be the event $\vol_n(K_\Pi)\ge c_v\kappa_n a^n$. 
By \cref{lb:prop:volume}, $\Prb\left[V\right]\ge c_p>0$.
Let $U$ be the event that no lattice of determinant at least $D_n$ covers $K_\Pi$.
By \cref{lb:prop:uniform-exclusion}, $\Prb\left[U\right]=1-o(1)$.
Therefore
\[
  \Prb\left[U\cap V\right]
  \ge\Prb\left[V\right]-\Prb\left[U^{\mathrm c}\right]
  \ge c_p-o(1)>0
\]
for sufficiently large $n$.
Fix a realization in $U\cap V$, and denote the resulting body by $K_n$.

Every lattice $\Lambda$ that covers $K_n$ satisfies $\det\Lambda<D_n$.
Hence
\[
  \frac{\vol_n(K_n)}{\det\Lambda}
  \ge\frac{c_v\kappa_n a^n}{D_n}
  =c_v c_{\mathrm{net}} n\log n.
\]
Since $K_n$ is origin-symmetric by construction, taking the infimum over all covering lattices proves \cref{thm:lower} with $c=c_v c_{\mathrm{net}}$.
By \cref{lb:prop:volume}, the John ellipsoid of $K_n$ is $B_2^n$.
\end{proof}

\begin{proof}[Proof of \cref{thm:main}]
The lower inequality follows from \cref{thm:lower}.
Every origin-symmetric convex body is a convex body, so the definitions give $\Theta_L^{\mathrm{sym}}(n)\le\Theta_L^{\mathrm{conv}}(n)$.
The upper inequality follows from \cref{thm:upper}.
Taking $c$ from \cref{thm:lower} and $C$ and $\varepsilon_n$ from \cref{thm:upper} proves \cref{thm:main}.
\end{proof}

\section*{Acknowledgments}

H.L. was supported by the National Natural Science Foundation of China (12501487), the China Scholarship Council, and the Institute for Basic Science (IBS-R029-C4). X.L. was supported by the Excellent Young Talents Program (Overseas) of the National Natural Science Foundation of China.

\section*{Declaration on the use of AI}

The authors used generative AI tools to assist in discussing proof strategies, checking proofs, and improving exposition. 

\bibliographystyle{abbrv}
\bibliography{../lattice_coverings_log_one}

\end{document}